\documentclass[12pt,twosided]{amsbook}
\usepackage{palatino, mathpazo}
\usepackage[a4paper,
margin=25mm]{geometry}
 \usepackage[OT2,T1]{fontenc}
\usepackage[utf8]{inputenc}
\usepackage{tikz-cd}
\usepackage{microtype,amsfonts,amssymb,amsthm,mathrsfs}
\usepackage[colorlinks=false]{hyperref}
\usepackage{graphicx}
\usepackage{wrapfig}
\usepackage{stmaryrd}
\usepackage[toc,page]{appendix}
\usepackage[colorinlistoftodos]{todonotes}
\usepackage{amssymb}
\usepackage{amsmath}
\usepackage{csquotes}
\usepackage{bm}
\usepackage{verbatim} 
\usepackage{latexsym}
\usepackage[all]{xy}
\usepackage{fancyhdr}
\usepackage{indentfirst}
\usepackage{color}
\usepackage{graphicx}
\usepackage{newlfont}
\usepackage{amscd}
\usepackage{amsmath}
\usepackage[colorinlistoftodos]{todonotes}
\usepackage{latexsym}
\usepackage{tikz-cd}
\usepackage{enumitem}
\theoremstyle{plain}
\newtheorem{thm}{theorem}[chapter]
\newtheorem{prop}[thm]{Proposition}
\newtheorem*{claim}{Claim}
\newtheorem{cor}[thm]{Corollary}
\newtheorem{lem}[thm]{Lemma}
\newtheorem{fact}[thm]{Fact}
\theoremstyle{definition}
\newtheorem{defin}{Definition}[chapter]

\newtheorem{exam}{Example}[chapter]
\theoremstyle{remark}
\newtheorem{rmk}{Remark}[chapter]

\newcommand{\E}{\mathcal E}

\newcommand{\C}{\mathbb C}
\newcommand{\D}{\mathcal D}

\newcommand{\R}{\mathbb{R}}

\newcommand{\N}{\mathbb N}

\newcommand{\OO}{\mathcal O}

\newcommand{\h}{\operatorname{Hom}}

\newcommand{\End}{\operatorname{End}}
\newcommand{\Spec}{\operatorname{Spec}}
\newcommand{\db}{\bar{\partial}}

\DeclareMathOperator\dom{Dom}
\DeclareMathOperator\supp{supp}
\DeclareMathOperator\im{im}
\newcommand{\loc}{\mathrm{loc}}

\newcommand{\p}{\partial}

\newcommand{\Imp}{\mathrm{Im}}
\newcommand{\comp}{\mathrm{comp}}
\newcommand{\ii}{\sqrt{-1}}

\newcommand{\cl}{\mathrm{cl}}

\numberwithin{section}{chapter}

\newcommand{\la}{\langle}
\newcommand{\ra}{\rangle}
\renewcommand{\h}{\widehat}

\numberwithin{equation}{chapter}
\begin{document}
\title{Asymptotic of Bergman Kernel}
\author{
{\Large \rm Yu-Chi Hou}\\
\vspace{7cm}
{\large \it A thesis presented for the degree of Master of Science}\\ 
{\large \it of Graduate Institute of Mathematics,}\\
{\large \it National Taiwan University, Taipei}\\
{\Large {\it supervised by Professor}   Chin-Yu Hsiao}\\
\vfill
{\large \it February 2022}
}
\begin{abstract}
In this thesis, we give a new proof on the pointwise asymptotic expansion for Bergman kernel of a hermitian holomorphic line bundle on the points where the curvature of the line bundle is positive and satisfy  local spectral gap condition. The main point is to introduce a suitable semi-classical symbol space and related symbolic calculus inspired from recent work of Hsiao and Savale. Particularly, we establish the existence of pointwise asymptotic expansion on the positive part for certain semi-positive line bundles.
\\ \hspace{\fill} 
\end{abstract}
\maketitle

\frontmatter
\tableofcontents
\mainmatter

\chapter{Introduction}
Bergman kernel, first introduced by Stefan Bergman \cite{Bergman} in 1922 for domains in $\C^n$, is a fundamental object in several complex variables and complex geometry, and has been the impetus of many researches in various fields. The main goal of this thesis is to study the asymptotic of Bergman kernel of holomorphic sections for high tensor powers of a holomorphic line bundle over a complex manifold. Let us first define Bergman kernel for a hermitian holomorphic line bundle over a complex manifold.

Let $L$ be a holomorphic line bundle over a complex manifold $X$ with $\dim_\C X=n$. If we endow a positive, smooth $(1,1)$-form $\omega$ on $X$, which induces a Riemannian volume form $d\nu_X=\omega_n:=\frac{\omega^n}{n!}$, and a hermitian metric $h^L$ on $L$ given by local weight $\phi$, then they give rise to a scalar product on $C^\infty_c(X,L)$, the space of smooth global sections for $L$ with compact supports. We then complete  $C^\infty_c(X,L)$ with respect to the scalar product to get a Hilbert space $L^2_{\omega,\phi}(X,L)$. The orthogonal projection $\Pi:L^2_{\omega,\phi}(X,L)\to \mathcal{H}^0(X,L)$ onto the subspace of $L^2$-integrable holomorphic sections of $L$ is called the \textbf{Bergman projection}, and its Schwartz kernel $K(z,w)$ is called the \textbf{Bergman kernel}. It is well-known that $K(z,w)$ is a smooth kernel (cf. Lemma \ref{lem:Bergman kernel is smoothing}), i.e., $K(z,w)\in C^\infty(X\times X,L^k\boxtimes L^{-k})$ and satisfies the \textbf{reproducing property}:
\begin{equation}
u(z)=(\Pi u)(z)=\int_X K(z,w)(u(w))d\nu_X(w),\quad u\in L^2_{\omega,\phi}(X,L)\cap H^0(X,L).
\end{equation} 

In general, it is difficult to calculate the Bergman kernel explicitly. However, when we replace $L$ by the $k$-th tensor power $L^k:=L^{\otimes k}$ and the hermitian metric $\phi$ by $k\phi$, the large $k$-behavior of the Bergman kernel is rather tractable and has important applications such as approximation of K\"ahler metrics by Fubini-Study metrics via the Kodaira map (\cite{Tian}, \cite{Bouche3}, and \cite{Ruan}), existence of canonical K\"ahler metrics (eg. \cite{Donaldson2001}, \cite{DS}, and \cite{CDS1}--\cite{CDS3}), Berezin--Toeplitz quantization (eg. \cite{Schlichenmaier1}, \cite{Schlichenmaier2}, and \cite{MamaBT}), and in physics \cite{Klevtsov}. 

If $L$ is  positive and $X$ is compact, then a well-known asymptotic formula asserts that there exists smooth functions $b_r(x)\in C^\infty(X)$, for $r\in \N_0$, such that for any $N,l\in \N$, there exists a constant $C:=C_{N,l}>0$ independent of $k$ satisfying
\begin{equation}\label{Zelditch--Catlin}
\|K_{k}(z,z)-\sum_{r=0}^N k^{n-r}b_r(x)\|_{C^l(X)}\leq Ck^{n-(N+1)}, \quad k\gg1.
\end{equation}
The existence of formula \eqref{Zelditch--Catlin} have been worked out in various generalities and through a variety of methods over the last thirty years. The leading asymptotic was first proved independently by Tian (1990, \cite{Tian}) using H\"ormander's $L^2$-estimates and by Bouche (1990, \cite{Bouche3}) using heat kernel. The full asymptotic was later developed independently by Catlin (1999, \cite{Catlin}) and Zelditch (1998, \cite{Zelditch}) using a result in CR geometry due to Boutet de Monvel and Sj\"orstrand (1975, \cite{BMS}). Later, Dai, Liu, and Ma (2006, \cite{DLM}) and Ma and Marienscu (2006, \cite{MM}) obtained both diagonal and off-diagonal expansions for generalized Bergman kernels for $\mathrm{spin}^c$-Dirac operators on compact symplectic manifolds based on the analytic localization technique due to Bismut and Lebeau. We refer the book of Ma and Marinescu \cite{Mama} and the references therein for this approach. 

Also, Berman, Berndtsson, and Sj\"ortrand (2008, \cite{BBS}) gave an alternative proof of \eqref{Zelditch--Catlin} by constructing approximate Bergman kernel inspired from analytic microlocal analysis and Hezari et al. (2016, \cite{Seto2}) (cf. also \cite{Seto}) gave a proof in similar spirit by constructing approximate kernel from Bargmann-Fock kernels. We also mention a recent result due to Rouby, Sj\"orstrand, and Ng\d{o}c (2020, \cite{Siorstrand2}) which proved the asymptotic expansion with exponentially decaying remainder when the metric $h^L$ is real-analytic by using Berezin quantization and analytic microlocal analysis.

If one drops the positive curvature assumption for $L$ and assume instead that the curvature is non-degenerate with constant signature $(n_+,n_-)$, then Berman and Sj\"ostrand (2007, \cite{BS}) showed the similar asymptotic expansions holds for orthogonal projection onto the space of harmonic $(0,q)$-form $\mathscr{H}^{0,q}(X,L^k)$ of Kodaira Laplacian (cf. section \ref{Gaffney Extension of Kodaira Laplacian}) if $X$ is compact and $q=n_-$. Independently, Ma and Marinescu (2006, \cite{MM2}) proved the analogous results in the setting of $\mathrm{spin}^c$-Dirac operators on compact symplectic manifolds. In \cite{HsiaoMa}, Hsiao and Marinescu (2014) proved that the spectral function for Kodaira Laplacian always admits local asymptotic expansion for any hermitian holomorphic line bundle on the non-degenerate points of the curvature, and they deduce local asymptotic of Bergman kernel under the additional spectral gap condition (cf. Definition \ref{defin:spectral gap condition}).
 
The philosophy behind the large $k$ asymptotic of Bergman kernel is closely tied to the notion of \textbf{positivity} and \textbf{ampleness} in algebraic geometry and complex geometry.  If $X$ is compact and $L$ is \textbf{positive} (i.e., the curvature form $R^L$ is positive definite), then classical \textbf{Kodaira vanishing Theorem} (cf. \cite[Theorem 1.5.6]{Mama}) implies that $H^q(X,L^k)=0$ for sufficiently large $k$, for any $q>0$. This is the key step in the tradition proof of \textbf{Kodaira embedding Theorem}, which asserts that 
if $L$ is positive line bundle and $X$ is compact, then any orthonormal basis $\{f_j\}_{j=1}^{d_k}$ of $H^0(X,L^k)$ gives rise to a holomorphic embedding
\begin{equation}
\begin{aligned}
\Phi_{L^k}:X&\to \mathbb{CP}^{d_k-1}\\
x&\mapsto [f_1(x):\cdots:f_{d_k}(x)]
\end{aligned}
\end{equation}
In fact, the full asymptotic of the Bergman kernel $K_k(z,w)$ also leads to a direct proof for Kodaira embedding Theorem (cf. \cite{Hsiao}). Now, combining Kodaira vanishing Theorem with classical Riemann--Roch formula
\[
\sum_{q=0}^n (-1)^q\dim H^q(X,L^k)=\frac{k^n}{n!}\int_{X} \left(\frac{\ii}{2\pi}R^L\right)^n+o(k^n)
\]
shows that $\dim H^0(X,L^k)=\frac{k^n}{n!}\int_X  \left(\frac{\ii}{2\pi}R^L\right)^n+o(k^n)$. 

If $L$ is only assumed to be semi-positive and positive at some point of $X$, called \textbf{quasi-positive}, then Siu (1984, \cite{Siu} ) showed the following asymptotic vanishing Theorem \[
\dim H^{0,q}(X,L^k)=o(k^n), \quad \forall q>0, k\gg1.
\]
By Riemann--Roch formula above again, we again obtain that 
\[
\dim H^0(X,L^k)=\frac{k^n}{n!}\int_X  \left(\frac{\ii}{2\pi}R^L\right)^n+o(k^n),
\] 
and thus $L$ is big. Siu then used this to prove a bimeromorphic analogue of Kodaira embedding Theorem, known as \textbf{Grauert--Riemenschneider conjecture}, which asserts that a compact, connected, complex manifold $X$ is \textbf{Moischezon}\footnote{A complex manifold $X$ is Moischezon if the fields of meromorphic functions on $X$ has transcendental degree $\dim X$ over $\C$.} if it admits a holomorphic hermitian line bundle $(L,h^L)\to X$ so that it is quasi-positive.

In the present thesis, we establish a pointwise asymptotic expansion for Bergman kernel. Instead of requiring $L$ to be globally positive, the expansion exists on the positively curved points which lie in an open set satisfying spectral gap condition (cf. Definition \ref{defin:spectral gap condition}). First, we define the approximate Bergman kernel on $\C^n$ and develop the related semi-classical $L^2$-estimates. Next, inspired by recent work of Hsiao and Savale \cite{HsiaoNikhill} in the CR setting, we define a semi-classical symbol spaces on Euclidean spaces and develop their symbolic calculus. Particularly, we establish the notion of asymptotic sum and compositions of their quantizations in such symbol spaces. Combining with symbolic calculus developed earlier, we obtain an asymptotic expansion for approximate Bergman kernel near the origin. Finally, under the assumption of spectral gap condition,  we prove that global Bergman kernel $K^{(q)}_k$ coincides with approximate kernel up to error term $O(k^{-\infty})$ (cf. Definition \ref{defin:k-negligible}) near the point where $L$ is positive.
\section{Statement of Main Results and Applications}\label{Statement of Main Results and Applications}
We now formulate the main result presented in the thesis. We refer to chapter \ref{db-Harmonic Forms and Bergman Kernel} for relevant notations and terminology used here.
Let $(X,\omega)$ be a hermitian manifold of complex dimension $n$, where $\omega$ is a smooth, positive $(1,1)$-form on $X$, inducing the hermitian structure on $X$. We denote $\la\cdot|\cdot\ra_\omega$ by the hermitian metric on $T^{1,0}X$ induced by $\omega$.  A canonical Riemannian volume form $d\nu_X$ for $(X,\omega)$ is given by $d\nu_X=\omega_n:=\frac{\omega^n}{n!}$. Let $L$ be a holomorphic line bundle on $X$ and set $L^k:=L^{\otimes k}$, for $k\in \N$. For any hermitian metric $h$ on $L$, we can define the Chern connection $\nabla$ on $L$ with respect to $h$ with curvature $R^L(h)=(\nabla^L)^2\in A^{1,1}(X)$. We identify $R^L(h)$ with the curvature operator $\dot{R}^L(h)\in C^\infty(X,\End(T^{1,0}X))$ by
\begin{equation}\label{curvature operator}
\sqrt{-1}R^L(h)(x)(v\wedge \overline{w})=\la \dot{R}^L(h)(x)v|w\ra_\omega,
\end{equation}
for any $x\in X$, $v,w\in T^{1,0}_xX$. We denote $n_+(x)$,$n_-(x)$,$n_0(x)$ by the number of positive, negative, and zero eigenvalues of $\dot{R}^L(h)$ at $x$. For $q=0,\dots,n$, we let
\begin{equation}
X(q):=\{x\in X:n_+(x)=n-q,n_-(x)=q,n_0(x)=0\}.
\end{equation}
Notice that $X(q)$ is an open set of $X$, for each $q\in \{0,1\dots,n\}$.

Locally, if $s$ is a holomorphic trivialization of $L$ over an open set $U\subset X$, then the hermitian metric $h$ is determined by
\[
|s|^2_h=e^{-2\phi},\quad \phi\in C^\infty(U,\R),
\]
called the \textbf{local weight} of $h$. On the $k$-th tensor power $L^k$ of $L$, $h$ induces a natural hermitian metric $h^k$ on $L^k$ with local weight $k\phi$. Let $\la \cdot|\cdot\ra_{k\phi}$ be the pointwise scalar product on the bundle $L^k$ (cf. \eqref{pointwise norm}) and $(\cdot|\cdot)_{\omega,k\phi}$ be the inner product on the space $C^\infty_c(X,L^k)$ of compact supported smooth sections of $L^k$, induced by $\omega$ and $h^k$ (cf. \eqref{L2-norm}). We denote $|\cdot|_{k\phi}$ and $\|\cdot\|_{\omega,k\phi}$ be the pointwise and $L^2$-norm associated to $\omega$ and $h^k$, and let $L^2(X,L^k)$ be the completion of $C^\infty_c(X,L^k)$ with respect to $\|\cdot\|_{\omega,k\phi}$. 

Let $\db:C^\infty(X,L^k)\to A^{0,1}(X,L^k)$ be the Cauchy--Riemann operator acting on smooth sections of $L^k$, $\db^*$ be the formal adjoint of $\db$ with respect to $(\cdot|\cdot)_{\omega,k\phi}$, and $\square_{\omega,k\phi}^{(0)}:=\db^*\db$ be the Kodaira Laplacian acting on $C^\infty(X,L^k)$ (cf. \eqref{Kodaira--Hodge Laplacian}). We denote by $\square_k^{(0)}$ by Gaffney extension of the Kodaira Laplacian (cf. Definition \ref{Gaffney extension}).  Let $\mathscr{H}^{0}(X,L^k)$ be the kernel of $\square_k^{(0)}$ and let $\Pi^{(0)}_k:L^2(X,L^k)\to \mathscr{H}^{0}(X,L^k)$ be the Bergman kernel for $L^k$-sections.

To state our result, we first define local spectral gap property.
\begin{defin}\label{defin:spectral gap condition}
For $d\in \R$ and an open set $D\subset X$,  we say that $\square_{\omega,k\phi}^{(0)}$ has \textbf{local spectral gap condition} of order $d$ on $D$ if there exists $C>0$ and $k_0\in \N$ such that for any $u\in C^\infty_c(D,L^k)$, if $k\geq k_0$,  
 \begin{equation}\label{spectral gap condition}
 \left\|(I-\Pi_{k}^{(0)})u \right\|_{\omega,k\phi}\leq \frac{1}{Ck^d} \left\|\square_{\omega,k\phi}^{(0)} u\right\|_{\omega,k\phi},
 \end{equation}
 where $\Pi_{k}^{(0)}$ is the Bergman projection from $L^2(X,L^k)\to \mathscr{H}^{0}(X,L^k)$.
\end{defin}
Next, we introduce the key ingredients in our approach. Namely, a kind of semi-classical symbol space inspired from the recent work of \cite{HsiaoNikhill}.
\begin{defin}\label{defin:semi-classical symbol space} For $m\in \R$, a function $a(x,y,k)$ with parameter $k\in \N$ is in $\widehat{S}^m(\R^d\times \R^d)$ if
\begin{itemize}
\item[(i)] $a(x,y,k)\in C^\infty(\R^d\times \R^d)$, for each $k\in \N$, and
\item[(ii)] for any$(\alpha,\beta)\in \N_0^{2d}$, there exists $l=l(\alpha,\beta)\in \N$ and $k_0\in \N$ such that for any $N\in \N$,  there exists a constant $C=C_{\alpha,\beta,N}(a)>0$ independent of $k$ satisfying
\begin{equation}\label{semi-classical symbol space}
\left|\p_x^\alpha\p_y^\beta a(x,y,k)\right|\leq Ck^{m+\frac{|\alpha|+|\beta|}{2}}\frac{(1+|\sqrt{k}x|+|\sqrt{k}y|)^{l}}{(1+|\sqrt{k}(x-y)|)^N},\end{equation}
for any  $(x,y)\in \R^d\times \R^d$, any $k\geq k_0$. 
\end{itemize}
Furthermore, we say $a\in \h{S}^m_{\cl}(\R^d\times \R^d)$ if $a\in \h{S}^m(\R^d\times \R^d)$ and there exists a sequence $a_j\in \h{S}(\R^d\times \R^d)$ (cf. Definition \ref{defin:Generalized Schwartz Space} for the definition of $\h{S}(\R^d\times \R^d)$), for $j\in \N_0$, so that
\begin{equation}\label{asymptotic sum of classical symbol}
a(x,y,k)-\sum_{j=0}^{N-1} k^{m-\frac{j}{2}}a_j(\sqrt{k}x,\sqrt{k}y)\in \h{S}^{m-\frac{N}{2}}(\R^d\times \R^d),\quad \forall N\in \N.
\end{equation}
\end{defin}
Also, let $s$ be a local holomorphic trivialization of $L$ over an open set $U\subset X$, we can make the following identification:
\[
A^{0,q}(U,L^k)\to A^{0,q}(U),\quad u=s^k\otimes \alpha\mapsto \widetilde{\alpha}:=\alpha e^{-k\phi}
\]
so that for any $u,v\in L^2_{0,q}(U,L^k)\cap A^{0,q}(U,L^k)$,
\[
(u|v)_{\omega,k\phi}=\int_U \la \alpha|\beta\ra_\omega e^{-2k\phi}\omega_n=\int_U \la e^{-k\phi}\alpha|e^{-k\phi}\beta\ra_\omega\omega_n=:(\widetilde{\alpha}|\widetilde{\beta})_\omega.
\]
Let $L^2_{0,q}(U,\omega)$ be the completion of $L^2_{0,q}(U)$ with respect to $(\cdot|\cdot)_{\omega}$ defined above. Clearly, above identification extends to an isometry $L^2_{0,q}(U,L^k)\cong L^2_{0,q}(U,\omega)$. We define the \textbf{localized Bergman kernel} with respect to $s$ by
\[
\Pi_{k,s}^{(q)}\alpha=e^{-k\phi}s^{-k}\Pi_k^{(q)}(e^{k\phi}\alpha\otimes s^k),
\]
where $s^{-k}$ is the dual section of $s^k$ so that $s^{-k}(s^k)\equiv 1$ on $U$. We denote $K_{k,s}$ by the Schwartz kernel of the localized Bergman kernel $\Pi_{k,s}$, called the \textbf{localized Bergman kernel}.

We now can state the main result for the thesis.
\begin{thm}[cf. Theorem \ref{thm:Pointwise Asymptotic of Bergman Kernel Function}]\label{Pointwise Asymptotic of Bergman Kernel Function}
Suppose $X(0)\neq \emptyset$, say $x\in X(0)$. For any $D\subset X(0)$ of $x$ satisfying the spectral gap condition (cf. \eqref{defin:spectral gap condition}), there exists a trivializing open set $U\Subset D$ and a holomorphic coordinate $z$ on $U$ centered at $x$ so that on $U$,  we have 
\[
\rho(z)K_{k,s}(z,w)\chi_k(w)\in \h{S}^n_{\cl}(\C^n\times \C^n),
\]
where $\rho\in C^\infty_c(U)$, $\chi_k(z):=\chi(8k^{1/2-\epsilon}z)$, $\chi\in C^\infty_c(\C^n)$ satisfying
\[
\supp\chi\subset B_1(0),\quad \chi=1\text{ on }B_{1/2}(0),\quad \rho=1\text{ near } 0,
\]
and $\epsilon\in (0,\frac{1}{6})$.
\end{thm}
From \eqref{asymptotic sum of classical symbol}, there exists a sequence $a_j\in \h{S}(\C^n\times \C^n)$ such that for any $N\in \N$,
\begin{equation}\label{asymptotic expansion for Bergman kernel}
\rho(z)K_{k,s}(z,w)\chi_k(w)-\sum_{j=0}^N k^{n-j/2}a_j(\sqrt{k}z,\sqrt{k}w)\in \h{S}^{n-\frac{N+1}{2}}(\C^n\times \C^n).
\end{equation}
Furthermore, we can calculate the first coefficients of the expansion.
\begin{thm}\label{the first coefficient}
In the setting of Theorem \ref{Pointwise Asymptotic of Bergman Kernel Function}, under the choice of holomorphic coordinates and holomorphic trivializations on $U$ (cf. Lemma \ref{Canonical Coordinate}) such that 
\begin{align*}
\phi(z)&=\phi_0(z)+O(|z|^3),\quad \phi_0(z):=\sum_{i=1}^n \lambda_{i,x}|z^i|^2,\quad \lambda_{i,x}>0,\\
\omega&=\omega_0(z)+O(|z|),\quad \omega_0(z):=\frac{\sqrt{-1}}{2}\sum^n_{j=1}dz^j\wedge d\bar{z}^j,
\end{align*}
the first coefficients in \eqref{asymptotic expansion for Bergman kernel} is given by
\[
a_0(z,w)=\frac{2^n\lambda_{1,x}\cdots \lambda_{n,x}}{\pi^n}e^{\sum_{j=1}^n\lambda_{j,x}(2z^j\overline{w}^j-|z^j|^2-|w^j|^2)}.
\]
Here, $4\lambda_{1,x},\dots,4\lambda_{n,x}$ is the eigenvalues of curvature operator $\dot{R}^L(h)$ at $x\in X$ (cf. \eqref{curvature operator}).
\end{thm}
By \eqref{semi-classical symbol space}, this means that given $N\in \N$, there exists $l=l(N)\in \N$ such that for any $M>0$, there exists a constant $C=C(N,M)>0$ satisfying
\[
\left|\rho(z)K_{k,s}(z,0)-\sum_{j=0}^N k^{n-j/2}a_j(z,0)\right|\leq Ck^{n-\frac{N+1}{2}}(1+|z|)^{l-M},\quad |z|<\frac{1}{2}k^{\epsilon-1/2}.
\]
If we further put $z=0$, then we obtain a pointwise asymptotic for 
\[
K_{k,s}(z)\sim \sum_{j=0}^\infty k^{n-j/2}a_j(z)
\]
 in the sense that for any $N\in \N$,
\[
\left|K_{k,s}(0,0)-\sum_{j=0}^N k^{n-j/2}a_j(0,0)\right|\leq Ck^{n-\frac{N+1}{2}}.
\]
This establishes the local pointwise asymptotic of Bergman kernel function on $X(0)$ with local spectral gap condition.

We now give a digression on spectral gap condition given in Definition \ref{defin:spectral gap condition}. First of all, it is clear that if \eqref{spectral gap condition} holds on $D$, it holds on any open subset of $D$. Now, we say $\square_{\omega,k\phi}^{(0)}$ satisfies \textbf{global spectral gap condition} if \eqref{spectral gap condition} holds for $D=X$. When $X$ is compact, by Hodge theorem, global spectral gap condition is equivalent to 
\[
\lambda_1(X,L^k):=\inf\{\lambda\in \Spec\square_{\omega,k\phi}^{(0)}:\lambda\neq 0\}\geq Ck^d.
\]
In other words,  it is equivalent to $\Spec \square_{\omega,k\phi}^{(0)}\subset \{0\}\cup [Ck^d,\infty)$.  We now give some known examples for spectral gap condition.
\begin{exam}[cf. \cite{Mama}, Theorem 1.5.5]
Given a compact complex manifold $X$, a positive line bundle $L$ with respect to a hermitian metric $h$, by Nakano inequality,  there exists constants $C_0,C_1>0$ such that for any $k\in \N$,
\[
\Spec\square_{\omega,k\phi}^{(0)}\subset \{0\}\cup (C_0k-C_1,\infty).
\] 
Hence, $\square_{\omega,k\phi}^{(0)}$ satisfies global spectral gap condition of order $1$.
\end{exam}
\begin{exam}  In  \cite{Siu}, Siu conjectured the following ''eigenvalue conjecture'': if $X$ is compact and $L$ is quasi-positive, then
\begin{equation}\label{Siu}
\inf_{k\in \N}\lambda_1(X,L^k)>0.
\end{equation}
Particularly, \eqref{Siu} implies that $\square_{\omega,k\phi}^{(0)}$ satisfies global spectral gap of order $N$, for some $N>0$. However, Donnelly \cite{Donnelly} demonstrated that Siu's conjecture is false in general. Moreover, let $S\to X$ be the unit circle bundle of $L$, which is a CR manifold, Donnelly also showed that \eqref{Siu} is true if the tangential Cauchy--Riemann operator $\db_b$ has closed range. From this, one can deduce that if $L$ is a positive line bundle with semi-positive metric, then \eqref{Siu} is true (with respect to the semi-positive metric). This particularly implies that \eqref{Siu} is true for any quasi-positive line bundle on compact Riemann surfaces.
\end{exam}
\begin{exam}
Let $(L,h^L)$ be a semi-positive holomorphic line bundle over a compact hermitian manifold $(X,\omega)$ with $\dim_\C X=n$. If we arrange the eigenvalue of $\dot{R}^L$ at $x$ as $0\leq \mu_1(x)\leq \cdots \leq \mu_n(x)$, then $\mu_1(x)$ is a continuous function on $X$. Bouche \cite{Bouche2} showed that if $\int_X \mu_1^{-6n}d\nu_X<\infty$, then $\lambda_1(X,L^k)\geq k^{\frac{10n+1}{12n+1}}$. Hence, Bouche condition implies that $\square_{\omega,k\phi}^{(0)}$ satisfies global spectral gap condition of order $s=\frac{10n+1}{12n+1}$. 
\end{exam}
Now, we consider non-compact examples for spectral gap condition.
\begin{exam}
Let $(L,h^L)$ be a semi-positive holomorphic line bundle over a complete K\"ahler manifold $(X,\omega)$ with $\dim_\C=X=n$, $\omega$ is a K\"ahler metric which is not necessarily complete. Then Demailly's $L^2$-estimate \cite{Demailly3} implies the following. If $g\in L^2_{0,1}(X,K_X\otimes L)$ satisfying $\db g=0$ and $\int_X |g|_{R^L}^2d\nu_X<\infty$, where $|g|_{R^L}(x):=\inf_{g'\in\bigwedge^{n,1}T^*X\otimes L}\frac{\langle\ii R^L \Lambda g'.g'\rangle}{\langle g,g'\rangle^2(x)}$, then there exists $f\in L^2(X,L\otimes K_M)$ with $\db f=g$ and
\[
\int_X |f|^2_{h^L}d\nu_X\leq \int_X |g|_{R^L}^2d\nu_X.
\]
From Demailly's result, Hsiao and Marinescu in \cite{HsiaoMa} proved that for any precompact open set $D\Subset X(0)$, $\square_{\omega,k\phi}^{(0)}$ has spectral gap of order $1$ on $D$.
\end{exam}
\begin{exam} Let $(X,\omega)$ be a compact hermitian manifold. Assume $(L,h^L)\to X$ is a smooth quasi-positive line bundle. Then by the solution of Grauert--Riemenschneider conjecture (cf. \cite[Chapter2]{Mama}), we know that $X$ is a Moischezon manifold and $L$ is a big line bundle.  From \cite[Lemma 2.3.6]{Mama}, $L$ admits a singular hermitian metric $h^L_{\mathrm{sing}}$ which is smooth outside an analytic set $\Sigma$ and whose curvature is strictly positive current. Hsiao and Marinescu in \cite[Lemma 8.1, Theorem 8.2]{HsiaoMa} proved that for any open set $D\Subset X(0)\cap (X\setminus \Sigma)$, $\square_{\omega,k\phi}^{(0)}$ for the open manifold $X\setminus \Sigma$ has spectral gap of order $N=-\sup_{x\in D}2(\phi(x)-\phi_{\mathrm{sing}}(x))$, where $\phi$ and $\phi_{\mathrm{sing}}$ are local weights of $h^L$ and $h^L_{\mathrm{sing}}$, respectively.
\end{exam}
\begin{rmk} We remark that there exists compact example for which spectral gap condition does not hold globally on whole space. In  \cite[Theroem 10.1]{HsiaoMa}, Hsiao and Marinescu demonstrated that there exists a compact complex manifold $X$ and a holomorphic line bundle $L\to X$ with hermitian metric $h^L$ such that for any $N>0$, $\lim_{k\to\infty}k^N\lambda_1(X,L^k)=0$. Hence, such $X$ cannot satisfy any global spectral gap condition. 
\end{rmk}
These examples illustrates that Theorem \ref{Pointwise Asymptotic of Bergman Kernel Function} asserts the existence of pointwise asymptotic of Bergman kernel in more general situation.

The present thesis is organized as following. In chapter \ref{db-Harmonic Forms and Bergman Kernel}, we briefly recall some preliminaries in complex geometry, including Kodaira Laplacian, its Gaffney extension, Bergman kernels.  In chapter \ref{chapter:Asymptotic Expansion of Bergman Kernel}, we prove Theorem \ref{Pointwise Asymptotic of Bergman Kernel Function}. We first set up localized version of Bergman kernel in section \ref{Notation and Set-Up}. Next, we extend these local operators and metrics to whole $\C^n$ , define the approximate kernel, and prove the required semi-classical $L^2$-estimate.  One of the key in our approach is a semi-classical symbolic calculus, which is developed in section \ref{symbolic calculus and Asymptotic Sum}. Using results in previous sections, we prove the asymptotic expansion for approximate kernel near the origin in section \ref{Asymptotic Expansion of Approximate Kernel}. In the section \ref{sec:Localization of Bergman Kernel via Spectral Gap}, we then show that under the spectral gap condition (cf. Definition \ref{defin:spectral gap condition}) near a point of $X(0)$,  Bergman kernel and approximate Bergman kernel are identical modulo $k$-negligible operator (cf. Definition \ref{defin:k-negligible}), and thus the pointwise asymptotic follows from the result.

\chapter{Preliminaries on \texorpdfstring{$\db$}{}-Harmonic Forms and Bergman Kernel}
\label{db-Harmonic Forms and Bergman Kernel}
In this chapter, we introduce Bergman kernel. The first three sections introduce some relevant notations and terminology. In section \ref{Gaffney Extension of Kodaira Laplacian}, we define a suitable extension of Kodaira Laplacians $\square^{(q)}_{\omega,k\phi}$ due to Gaffney for the study of $\db$-harmonic forms on non-compact case. Finallly, in section \ref{Bergman Kernel} , we introduce Bergman kernel and discuss some fundamental properties of them. 
\section{Some Standard Notations}
We denote $\N:=\{1,2,\dots\}$  by the set of natural numbers and $\N_0:=\N\cup\{0\}$. We adopt the following two multi-indices notations. For a multi-index $\alpha=(\alpha_1,\dots,\alpha_n)\in \N_0^n$, we denote $|\alpha|:=\sum_{i=1}^n\alpha_i$. We adopt standard notations such as $\alpha!=\alpha_1!\cdots\alpha_n!$, $x^\alpha$, and $\p_x^\alpha$. On the other hand, a $n$-tuple $J=(j_1,\dots,j_q)\in \{1,\dots,n\}^q$ is called a \textbf{strictly increasing multi-index} if $1\leq j_1<\cdots <j_q\leq n$. For a differential $q$-form $\alpha$, the local expression in local coordinate $x=(x^1,\dots,x^n)$ is given by
\[
\alpha=\sum_{|I|=q}\nolimits^\prime \alpha_Idx^I,
\]
where $\sum'_{|I|=q}$ means that summation is over strictly increasing multi-indices $I$ of length $q$. Also, we denote $dm$ by the standard Lebesgue measure on Euclidean spaces and $B_r(z)$ by the open ball with radius $r>0$ and center $z\in \C^n$  

Let $X$ be a complex manifold. We introduce some standard notations of various function spaces. For any open subset $U\subset X$, we denote $\OO_X(U)$ by the space of holomorphic function on $U$. In case of $X=\C^n$, we denote $\OO_{\C^n}(U)$ by $\OO(U)$. We also denote $C^\infty(U)$ and $C^\infty_c(U)$ by the space of smooth functions and the test functions on $U$, respectively. If $E\to X$ is a complex vector bundle, we denote $C^\infty(U,E)$ and $C^\infty_c(U,E)$ by the space of smooth sections and its subspace whose elements having compact supports on an open subset $U\subset X$. Similarly,  we denote $\mathcal{D}'(U,E)$ and $\mathcal{E}'(U,E)$ by the space of distribution sections of $E$ over $U$ and its subspace whose elements having compact supports. For $t\in \R$, we denote $W^{t}(U,E)$ by the Sobolev space\footnote{The usual notation for $L^2$-Sobolev space is $H^s$. However, to avoid the confusion with cohomology group, we denote it by $W^s$.} of order $t$ of sections of $E$ over $U$, 
 \begin{align*}
&W^t_{\loc}(U,E):=\{u\in \mathcal{D}'(U,E):\phi u\in W^t(U,E),\forall \phi\in C^\infty_c(U)\},\text{ and }\\
&W^t_{\comp}(U,E):=W^t_{\loc}(U,E)\cap \E'(U,E).
\end{align*}
Finally, if $F$ is a holomorphic vector bundle, we denote $H^0(U,F)$ by the space of holomorphic sections of $F$ over $U$.
\section{Hermitian Metrics}
For a holomorphic line bundle $L\to X$ with hermitian metric $h$, given a holomorphic trivialization $\{(U_i,s_i)\}$ of $L\to X$, where $\{U_i\}$ is an open covering of $X$ and $s_i:U_i\to L$ is a non-zero holomorphic section of $L$, we define the \textbf{local weight} $\phi_j\in C^\infty(U_j,\R)$ of $h^L$ with respect to $s_i$ by $\phi_j:=\frac{-1}{2}\log |s_j|^2_h$, or equivalently, $|s_j|^2_h=e^{-2\phi_j}$. For $x\in U_i\cap U_j$, assume $s_i=\sigma_{ij}s_j$. Here, $\sigma_{ij}\in \OO_X(U_{i}\cap U_j)$ is nowhere vanishing. We then have
\begin{equation}\label{transition}
e^{-2\phi_i}=|s_i|_{h^L}^2=\langle \sigma_{ij}s_j|\sigma_{ij}s_j\rangle=|\sigma_{ij}|^2|s_j|^2=|\sigma_{ij}|^2e^{-2\phi_j}, \quad \text{ on }U_i\cap U_j.
\end{equation}
In other words, $\phi_i=\phi_j-\frac{1}{2}\log |\sigma_{ij}|^2$. Hence a hermitian metric $h$ is equivalent to a collection of local weight $\{\phi_i\}$ satisfying \eqref{transition} with respect to a holomorphic trivialization $(U_i,s_i)$ of $L$. Hence, for $u,v\in C^\infty(U_i,L)$, we write $u=f_is_i$, $v=g_is_i$, for some $f_i,g_i\in C^\infty(U_i,\C)$, $\langle s|t\rangle_{h}=f_i\overline{g_i}e^{-2\phi_i}$. \par

For $L^k:=L^{\otimes k}$, a natural hermitian metric $h^k$ on $L^k$ induced from $h$ is given by $h^k(s^k)=h(s)^k$, for any $s\in C^\infty(U,L^k)$, for any open subset $U\subset X$. We then see that the canonical hermitian metric $h^k$ on $L^k$ is given by the local weight $\{k\phi_i\}$ if $\{\phi_i\}$ is a family of local weight for $(L,h^L)$ with respect to a holomorphic trivialization. From now on, we denote $\langle u|v\rangle_h$ by $\langle u|v\rangle_\phi$ and we interchangeably use local weight $\phi$ to represent the hermitian metric $h$ on a holomorphic line bundle $L$.

For a complex manifold $X$, we have a natural almost complex structure $J:TX\to TX$ from multiplication by $\sqrt{-1}$. Hence, $J$ induces an eigenspace decomposition $TX\otimes _\R\C=T^{1,0}X\oplus T^{0,1}X$,
where $T^{1,0}X$ is the $\sqrt{-1}$-eigenspace of $J$ and $T^{0,1}X$ is the $-\sqrt{-1}$-eigenspace of $J$. Also, $J$ induces an almost complex structure on $T^*X$. Hence, we also have the eigenspace decomposition for complexified cotangent bundle
\[
T^*X\otimes_\R\C=\bigwedge\nolimits^{1,0}X\oplus \bigwedge\nolimits^{0,1}X.
\]
Moreover, this extends to exterior algebra of complexified cotangent bundle:
 \[
 \bigwedge\nolimits^rT^*X\otimes_\R\C=\bigoplus_{p+q=r} \bigwedge\nolimits^{p,q}X,
 \]
 where $\bigwedge^{p,q}X$ is locally spanned by $dz^I\wedge d\bar{z}^J$, for any strictly increasing multi-indices $I\in \{1,\dots,n\}^p$, $J\in \{1,\dots,n\}^q$. We denote $A^r(U)=C^\infty(U,\bigwedge^rT^*X)$ and $A^{p,q}(U)=C^\infty(U,\bigwedge^{p,q}X)$ by the space of smooth $r$-forms and smooth $(p,q)$-forms on $U$, respectively.

Recall that a \textbf{hermitian form} on a complex manifold $X$ is a smooth $(1,1)$-form $\omega\in A^{1,1}(M)$ such that in a local holomorphic coordinate $\mathbf{z}=(z^1,\dots,z^n)$ on a chart  $U$ of $X$,
\begin{equation}
\omega\big|_U=\frac{\ii}{2}\sum_{i,j=1}^n H_{ij}dz^i\wedge d\bar{z}^j,
\end{equation}
where $H(x)=(H_{ij}(x))_{i,j=1}^n$ is a hermitian matrix for any $x\in U$. It is well-known that a hermitian form $\omega$ is equivalent to a Riemannian metric $g$ on the underlying real manifold $X$ which the complex structure is an isometry. We then extend $g$ to a hermitian metric $g_\C$ on $TX\otimes_\R\C$ by 
\[
g_\C(v\otimes \lambda, w\otimes \mu):=\lambda \overline{\mu}g(v,w),
\] 
where $v,w\in T_xX$ and $\lambda,\mu\in \C$. Thus, we can define a pointwise hermitian inner product $\langle\cdot|\cdot\rangle_\omega$ on $A^{p,q}(X)$ induced from $g_\C$. Now, let $\alpha,\beta\in A^{p,q}(X,L^k):=C^\infty(X,\bigwedge^{p,q}X\otimes L^k)$ be two $L^k$-valued $(p,q)$-forms. Under a choice of trivialization $s:U\to L$ of $L$, we can write $\alpha=f\otimes s^k$, $\beta=g\otimes s^k$. We define
\begin{equation}\label{pointwise norm}
\langle \alpha|\beta\rangle_{\omega,\phi}:=\langle f|g\rangle_\omega e^{-2\phi},
\end{equation}
where $\phi$ is the local weight of $h$ associated to $s$. We then define a $L^2$-hermitian inner product $(\cdot|\cdot)$ on $A^{p,q}_c(X,L^k)$, the space of compact supported $(p,q)$-forms valued in $L^k$ by
\begin{equation}\label{L2-norm}
(\alpha|\beta)_{\omega,k\phi}:=\int_X \langle \alpha|\beta\rangle_{\omega,k\phi}d\nu_X.
\end{equation}
We write $\|\alpha\|^2_{\omega,k\phi}:=(\alpha|\alpha)_{\omega,k\phi}$ and denote $L^2_{p,q}(X,L^k)$ by the completion of $A^{p,q}_c(X,L^k)$ with respect to the norm $\|\cdot\|_{\omega,k\phi}$. We sometimes denote $L^2_{p,q}(X,L^k)$ by $L^2_{\omega,k\phi}(\bigwedge^{p,q}X\otimes L^k)$ if we wish to stress the choice of $\omega$ and $k\phi$.

\section{Differential Operators and Curvatures}\label{Differential Operators and Curvatures}
Recall that given a holomorphic vector bundle $E$ over a complex manifold $X$, we can always define the \textbf{Cauchy--Riemann operator} $\db^E:A^{p,q}(X,E)\to A^{p,q+1}(X,E)$ satisfying $(\db^E)^2=0$ and the Leibniz rule $\db^E(fs)=\db(f)\wedge s+f\db^Es$, for any $s\in A^{p,q}(X,E)$, any $f\in C^\infty(X)$. It is well-known that we can always choose a canonical connection on $E$ which is compatible with a given hermitian metric $h^E$ on $E$.
\begin{fact}[cf. \cite{GH}, p.73] Let $E$ be any holomorphic vector bundle over a complex manifold $X$, $h^E$ be a hermitian metric on $E$. Then there exists a unique connection $\nabla^E$, called the \textbf{Chern connection} on $E$ such that
\begin{itemize}
\item[(i)] for any  $s\in C^\infty(X,E)$, $\nabla^{0,1} s=\db^Es$, where $\nabla^{0,1}s$ means taking the $(0,1)$-part of $\nabla s$.
\item[(ii)] for any $s_1,s_2\in C^\infty(X,E)$, $d(h^E(s_1,s_2))=h^E(\nabla s_1,s_2)+h^E(s_1,\nabla s_2)$.
\end{itemize}
In fact, the connection $1$-form for $\nabla^E$ is given by $\Theta=(\p H)H^{-1}$, where $H$ is the matrix representation of $h^E$ with respect to a a holomorphic frame $s=(s_1,\dots,s_r)$ of $E$.
\end{fact}
Recall that the curvature form $R^E(h^E)\in A^2(X)$ with respect to the Chern connection $\nabla$ is given by $\nabla^2s=R^E(h^E)s$. Since $\Theta=(\p H)H^{-1}$ and $\p(H^{-1})=-H^{-1}(\p H)H^{-1}$,
\begin{align*}
R^E&=d\Theta-\Theta\wedge \Theta=(\p+\bar{\p})(\p H H^{-1})- (\p H)H^{-1}\wedge (\p H)H^{-1}\\
&=-\p H\wedge \p (H^{-1})+(\bar{\p}\p H)H^{-1}-\p H\wedge \bar{\p} H^{-1}- (\p H)H^{-1}\wedge (\p H)H^{-1}\\
&=(\bar{\p}\p H)H^{-1}+(\bar{\p} H)H^{-1}\wedge(\p H)H^{-1}=\db^E\Theta\in A^{1,1}(\mathrm{End}E).
\end{align*}
Particularly, let $L\to X$ be a holomorphic line bundle, $h$ be a hermitian metric on it. If $(s,U)$ is holomorphic trivialization of $L$ over $U$ and $\phi$ is the local weight of $h$ determined by $s$. In this case, the curvature form is locally given by
\begin{equation}
\label{curvature of line bundle}
R^L(h)=-\p\bar{\p}\log e^{-2\phi}=2\p\db \phi=2\sum_{j,l=1}^n\frac{\p^2\phi}{\p z^j\p \bar{z}^l}dz^j\wedge d\bar{z}^l.
\end{equation}
In particular, $\ii R^L(h)$ is a closed, real $(1,1)$-form on $X$. We define the \textbf{curvature operator} $\dot{R}^L\in \End(\bigwedge^{1,0}X)$ as in \eqref{curvature operator}.

We end this section by proving a Lemma which asserts that the local weight and base metric can always have a normal form.
\begin{lem}\label{Canonical Coordinate} Let $X$ be a complex manifold with hermitian form $\omega$, $L$ be a holomorphic line bundle on $X$ with a hermitian metric $h$. Fix a point $x\in X$, we can choose a local complex coordinate $(z^1,\dots,z^n)$ on an open neighborhood $U\subset X$ of $x$ and a holomorphic trivializing section $s\in H^0(U,L)$ such that
\begin{itemize}
\item[(i)] $z^i(x)=0$, for $i=1,\dots,n$,
\item[(ii)] $\omega(z)=\frac{\ii}{2}\sum_{i,j=1}^nH_{ij}(z)dz^i\wedge d\bar{z}_j$ with $H_{ij}(0)=\delta_{ij}$, and
\item[(iii)] $|s(z)|^2_h=e^{-2\phi(z)}$ with local weight 
\[
\phi(z)=\sum_{i=1}^n \lambda_{i,x}|z^i|^2+O(|z|^3),
\] 
where $4\lambda_{1,x},\dots,4\lambda_{n,x}$ are eigenvalues of $\dot{R}^L(x)$.
\end{itemize}
We usually denote $\phi_0(z)=\sum_{i=1}^n \lambda_{i,x}|z^i|^2$.
\end{lem}
\begin{proof}
Let $(w_1,\dots,w_n)$ be any complex coordinate of $X$ on an open chart $V$. First, for each $p\in V$, we choose $\widetilde{z}^i(p):=w^i(p)-w^i(x)$, for $i=1,\dots,n$. Hence, $\widetilde{z}^i(x)=0$, for any $i=1,\dots,n$. Next, suppose the hermitian form $\omega\big|_V$ is given by
\[
\omega(\widehat{z})=\frac{\ii}{2}\sum_{i,j=1}^n \widehat{H}_{ij}(\widetilde{z})d\widetilde{z}^i\wedge d\overline{\widetilde{z}}^j.
\]
Since $\widehat{H}(0):=(\widehat{H}_{ij}(0))$ is hermitian and positive definite, there exists unitary matrix $U=(U_{ij})$ such that $U\widehat{H}(0)U^{-1}=I_n$. We then pick $z^i=\sum_{j=1}^nU_{ij}\widehat{z}_j$. In this case, $\frac{\p }{\p z^i}=\sum_{j=1}^nU_{ij}\frac{\p }{\p \widehat{z}_j}$. Therefore, $\omega(z)=\frac{\ii}{2}\sum_{i,j=1}^nH_{ij}(\widehat{z})dz^i\wedge d\bar{z}_j$ and $H_{ij}(0)=\delta_{ij}$. In sum, we have found a holomorphic coordinate $(z_1,\dots,z_n)$ on an open chart $V$ of $x$ such that  (i) and (ii) hold.

For (iii), we first shrink $V$ to an open neighborhood $U$ of $x$ such that $L\big|_U$ is trivialized by a holomorphic section $\widehat{s}:U\to L$. Suppose $|\widehat{s}|_h^2=e^{-2\widehat{\phi}(z)}$ with $\widehat{\phi}\in C^\infty(U,\R)$. We then expand $\widehat{\phi}$ in Taylor series up to the second order with respect to above local holomorphic coordinate $(z_1,\dots,z_n)$:
\begin{align*}
\widehat{\phi}(z)&=\widehat{\phi}(0)+\sum_{i=1}^n \left[\frac{\p \widehat{\phi}}{\p z^i}(0)z^i+\frac{\p \widehat{\phi}}{\p \bar{z}^i}(0)\bar{z}^i\right]\\
&+\sum_{i,j=1}^n \left[\frac{\p^2 \widehat{\phi}}{\p z^i\p z^j}(0)z^iz^j+\frac{\p^2 \widehat{\phi}}{\p z^i\bar{z}^j}(0)z^i\bar{z}^j+\frac{\p^2\widehat{\phi}}{\p \bar{z}^i\bar{z}^j}(0)\bar{z}^i\bar{z}^j\right]+O(|z|^3).
\end{align*}
We then define a new holomorphic trivializing section $s:U\to L$ by $s(z)=e^{G(z)}\widehat{s}(z)$, where
\[
G(z)=\widehat{\phi}(0)+2\sum_{i=1}^n \frac{\p \widehat{\phi}}{\p z^i}(0)z^i+\sum_{i,j=1}^n\frac{\p^2\widehat{\phi}}{\p z^i\p z^j}(0)z^iz^j.
\]
Then the local weight $\phi$ of $h$ with respect to this new section $s$ is given by
\[
e^{-2\phi}=|s|_h^2=e^{G(z)+\overline{G(z)}}e^{-2\widehat{\phi}(z)},
\]
and hence
\[
\phi(z)=\widehat{\phi}-\frac{G(z)+\overline{G(z)}}{2}=\sum_{i,j=1}^n\frac{\p^2 \widehat{\phi}}{\p z^i\bar{z}^j}(0)z^i\bar{z}^j+O(|z|^3).
\]
Since $\frac{\p^2\widehat{\phi}}{\p z^i\bar{z}^j}(0)$ is again a hermitian matrix, we can unitarily diagonalize it. In other words, we can rechoose $(z_1,\dots,z_n)$ such that (i) and (ii) still hold and $\phi(z)$ becomes
\[
\phi(z)=\sum_{i=1}^n \lambda_{i,x}|z^i|^2+R(z,\bar{z}),\quad R(z,\bar{z})=O(|z|^3),
\]
in this case. Moreover, since $R^L(h)=2\p \db \phi=2\sum_{i,j=1}^n \frac{\p^2\phi}{\p z^i\p \bar{z}^j}dz^i\wedge \bar{z}^j$, we recognize that $4\lambda_{1,x},\dots,4\lambda_{n,x}$ are the eigenvalues of $\dot{R}^L(h)(x)$.
\end{proof}

\section{Gaffney Extension of Kodaira Laplacian}\label{Gaffney Extension of Kodaira Laplacian} Let  $(L,h)$ be a hermitian holomorphic line bundle over a hermitian manifold $(X,\omega)$. In section \ref{Differential Operators and Curvatures}, we have defined the differential operator $\db:=\db^{L^k}:A^{0,q}(X,L^k)\rightarrow A^{0,q+1}(X,L^k)$, for $q\in \{0,1\dots,n\}$. We denote  $\db^{*}:=\db^{L^k,*,\omega,k\phi}:A^{0,q}(X,L^k)\rightarrow A^{0,q-1}(X,L^k)$ by the formal adjoint of $\db$ with respect to $(\cdot|\cdot)_{\omega,k\phi}$ which is characterized by
\begin{equation}\label{formal adjoint}
(\db\alpha|\beta)_{\omega,k\phi}=(\alpha|\db^{*}\beta)_{\omega,k\phi},\quad \alpha\in A^{0,q}_c(X,L^k),\beta\in A^{0,q+1}(X,L^k).
\end{equation}
 \textbf{Kodaira Laplacian} for a hermitian holomorphic line bundle $(L,h)$ is defined by
\begin{equation}\label{Kodaira--Hodge Laplacian}
\square^{(q)}_{\omega,k\phi}:=\db\db^{*}+\db^{*}\db:A^{0,q}(X,L^k)\rightarrow A^{0,q}(X,L^k).
\end{equation}
To apply standard functional analysis formalism, we need to extend $\square^{(q)}_{\omega,k\phi}$ to an unbounded operator on $L^2_{p,q}(X,L^k)$. First, we consider an extension of $\db$.
\begin{defin}\label{maximal extension}
The \textbf{maximal extension} of $\db:A^{0,q}(X,L^k)\to A^{0,q+1}(X,L^k)$, denoted by $S_q:\dom S_q\subset L^2_{0,q}(X,L^k)\to L^2_{0,q+1}(X,L^k)$, is defined by
\[
\dom S_q:=\{u\in L^2_{0,q}(X,L^k):\db u\in L^2_{0,q+1}(X,L^k)\},\quad S_qu:=\db u,
\]
where $\db u$ is defined in the distribution sense by $(\db u|v)_{\omega,k\phi}=(u|\db^{*}v)_{\omega,k\phi}$, for any $v\in A^{p,q}_c(X,L^k)$.
\end{defin}
In other words, $\dom S_q$ consists of $L^2$-integrable $(0,q)$-form $u$ valued in $L^k$ so that the current $\db u$ is representable by a $L^2$-integrable $(0,q+1)$-form. Since $L^2_{0,q}(X,L^k)$ is the completion of $A^{0,q}_c(X,L^k)$ with respect to the $(\cdot|\cdot)_{\omega,k\phi}$ and thus $A^{0,q}_c(X,L^k)\subset \dom S_q$, it is a densely defined, closed extension of the Cauchy--Riemann operator (\cite[Lemma 3.1.1]{Mama}). We then define the Hilbert space adjoint $S_q^\dagger:\dom S_q^\dagger\subset L^2_{0,q+1}(X,L^k)\to L^2_{0,q}(X,L^k)$ of $S_q$, which is always a closed operator. (cf. \cite[Section 1.2]{Davies}). Since $A_c^{0,q+1}(X,L^k)\subset \dom (S_q^\dagger)$, $S_q^\dagger$ is also a densely defined, closed operator. \par

Following Gaffney \cite{Gaffney}, we define an extension for Kodaira Laplacian $\square^{(q)}_{\omega,k\phi}$ as following.
\begin{defin}\label{Gaffney extension}
\textbf{Gaffeney extension} \[
\square^{(q)}_k:\dom \square^{(q)}_k\subset L^2_{0,q}(X,L^k)\to L^2_{0,q}(X,L^k)
\]
of Kodaira Laplacian $\square^{(q)}_{\omega,k\phi}$ is defined by 
\[
\dom \square^{(q)}_k:=\{u\in \dom S_q\cap \dom S_{q-1}^\dagger:S_q u\in \dom S_q^\dagger,S_{q-1}^\dagger u\in \dom S_{q-1}\}
\]
and $\square_k^{(q)}u:=S_q^\dagger S_q u+S_{q-1}S_{q-1}^\dagger u$ for $u\in \dom \square^{(q)}_k$.
\end{defin} 
Notice that $A^{0,q}_c(X,L^k)\subset \dom S_q\cap \dom S_{q-1}^\dagger$ and for $u\in A_c^{0,q}(X,L^k)$,
\[
S_q=\db u\in A^{0,q+1}_c(X.L^k)\subset \dom S_q^\dagger,\quad S_{q-1}^\dagger u=\db^*u\in A^{0,q-1}(X,L^k)\subset \dom S_{q-1}.
\]
Hence, $A^{0,q}_c(X,L^k)\subset \dom \square^{(q)}_k$ and thus  $\square^{(q)}_k$ is a densely defined operator.

Such extension enjoys nice analytic properties due to the following fact, whose proof is standard and can be found in \cite[Proposition 3.1.2]{Mama} or \cite[Proposition 1.3.8]{FollandKohn}.
\begin{fact}[Gaffney]\label{Gaffney}
Gaffney extension $\square_k^{(q)}$ is a non-negative self-adjoint extension of the Kodaira Laplacian $\square_{\omega,k\phi}^{(q)}$.
\end{fact}
\begin{rmk}
In general, the Kodaira Laplacian is not essentially self-adjoint, and there are many possible self-adjoint extension of Kodaira Lapalcians if $X$ is non-compact or $\p X\neq \emptyset$. For instance, there exists another useful extension known as \textit{Friedrich extension}, see \cite[Theorem 4.4.5]{Davies}, \cite[Proposition 1.3.3]{FollandKohn}, \cite[Section A.2.3.]{Marinescu1}.

However, two extensions coincide if $X$ is complete (in this case, the Kodaira Laplacian is essentially self-adjoint, cf. \cite[Corollary A.14]{Marinescu1}). We refer to \cite[Appendix A.2]{Marinescu1} for proofs of statements above.
\end{rmk}
\begin{defin}
The \textbf{space  of $L^2$ harmonic $(0,q)$-form values in the bundle $L^k$} is defined by 
\[
\mathscr{H}^{q}(X,L^k):=\ker\square_k^{(q)}=\{u\in \dom \square^{(q)}_k:\square_k^{(q)}u=0\}.
\]
\end{defin}
By the definition of Gaffney extension, we have 
\begin{lem}
$\mathscr{H}^{q}(X,L^k)=\ker S_q\cap \ker S_{q-1}^\dagger$.
\end{lem}
\begin{proof} 
It is clear that $\ker S_q\cap \ker S_{q-1}^\dagger\subset \mathscr{H}^q(X,L^k)$. Conversely, if $u\in \dom \square^{(q)}_k$, then $u\in \dom S_q\cap \dom S_{q-1}^\dagger$ and
\[
0=(u|\square_k^{(q)}u)_{\omega,k\phi}=(u|S_{q}^\dagger S_q u)_{\omega,k\phi}+(u|S_{q-1} S_{q-1}^\dagger u)_{\omega,k\phi}=\|S_q u\|^2_{\omega,k\phi}+\|S_{q-1}^\dagger u\|^2_{\omega,k\phi}.
\]
Hence, we conclude that $S_q u=0$ and $S_{q-1}^\dagger u=0$.
\end{proof}
It is well--known that $\square_{\omega,k\phi}^{(q)}$ is an elliptic operator of second order. If $u\in \mathscr{H}^{q}(X,L^k)$, i.e., $\square^{(q)}_ku=0$, then elliptic regularity shows that $u\in A^q(X,L^k)$. Hence, we conclude that
\[
\mathscr{H}^{q}(X,L^k)=\{u\in \dom \square_k^{(q)}\cap A^q(X,L^k):\db u=\db^\dagger u=0\}
\] 
For the case $q=0$, $\square_k^{(0)}=S_{0}^\dagger S_0$. If $u\in\ker\square_k^{(0)}$, then 
\[
(u|S_0^\dagger S_0 u)_{\omega,k\phi}=\|S_0 u\|^2_{\omega,k\phi}=0.
\] 
Hence, $u\in \ker S_0$. This means that $u\in C^\infty(X,L^k)\cap L^2(X,L^k)$ satisfying $\db{u}=0$ in the distribution sense. Therefore, $u\in H^0(X,L^k)$.
\section{Bergman Kernel}\label{Bergman Kernel}
In the present section, we introduce the primary object of the our study, \textbf{Bergman kernel} for a hermitian holomorphic line bundle on a hermitian manifold. Before proceeding, let us first recall the Schwartz kernel Theorem.
\begin{fact}[Schwartz Kernel Theorem, cf. \cite{HormanderI}, Theorem 5.2.1, 5.2.6]\label{thm:Schwartz Kernel Theorem}~
\begin{itemize}
\item[(i)] Let $A:C^\infty_c(M,E)\to \D'(M,F)$ be a continuous linear operator. Then there exists a unique distribution $K\in \D'(M\times M, F\boxtimes E^*)$, called the \textbf{Schwartz kernel} such that
\[
(Au)(v)=K_A(v\otimes u), \quad u\in C^\infty_c(M,E), v\in C^\infty_c(M,F).
\] 
Here, $F\boxtimes E^*$ is the vector bundle on $M\times M$ whose fiber at $(x,y)\in M\times M$ is $\mathrm{Hom}(F_y,E_x)$, the space of linear transformations from $F_y$ to $E_x$.  Also, $v\otimes u\in C^\infty_c(M\times M,F\boxtimes E)$ is defined by $v\otimes u(x,y):=v(x)\otimes u(y)$. 
\item[(ii)] $A$ extends to a continuous linear map from $\E'(M,E)$ to $C^\infty(M,F)$ if and only if $K_A$ is represented by a smooth section $K_A(x,y)\in C^\infty(M\times M,F\boxtimes E^*)$ and
\[
(Au)(x):=u(K(x,\cdot)),\quad \forall u\in \E'(M,E).
\]
Moreover, with respect to any volume form $d\nu_M$, $A$ is given by
\[
Au(x)=\int_M K_A(x,y)(u(y))d\mu_M(y),\quad \forall u\in C^\infty_c(M,E).
\]
In such case, then we say $A$ is a \textbf{smoothing operator} and $K_A$ is a \textbf{smoothing kernel}. Equivalently, by Sobolev embedding, $A$ is a smoothing operator if and only if $A$ extends to a continuous linear operator from $W^{t}_{\comp}(M,E)$ to $W^{t+N}_{\loc}(M,F)$ for any $t\in \R$ any $N\in \N$.
\end{itemize}
\end{fact}

Let $(X,\omega)$ be a hermitian manifold, where $\omega$ is a hermitian form on $X$, $L\to X$ be a holomorphic line bundle with a hermitian fiber metric $h$ with local weight given by $\phi$. Since $\mathscr{H}^{q}(X,L^k):=\ker \square^{(q)}_k$ is a closed subspace of $L^2_{0,q}(X,L^k)$, one can consider the orthogonal projection $\Pi^{(q)}_{k}:L^2_{0,q}(X,L^k)\to \mathscr{H}^{q}(X,L^k)$, called the \textbf{Bergman projection} for $(0,q)$-forms. We then define
\begin{defin}\label{def:Bergman kernel}
The \textbf{Bergman kernel} $K_{k}^{(q)}(z,w)$ is the Schwartz kernel  representing the Bergman projection $\Pi_{k}^{(q)}:L^2_{0,q}(X,L^k)\to \mathscr{H}^q(X,L^k)$.
\end{defin}
A priori, $K_{k}^{(q)}\in \mathcal{D}'(X\times X, (L^k\otimes \bigwedge^{0,q}X)\boxtimes (L^k\otimes \bigwedge^{0,q}X)^*)$. We now prove
\begin{lem}\label{lem:Bergman kernel is smoothing} $\Pi_k$ is a smoothing operator.
\end{lem} 
\begin{proof}
We need to prove that $\Pi_k$ extends to a continuous linear operator from $W^{s}_{\comp}(X,\bigwedge^{0,q}X\otimes L^k)$ to $W^{s+N}_{\loc}(X,\bigwedge^{0,q}\otimes L^k)$, for any $s\in \R$ and $N\in \N$. By definition, it suffices to prove that for any $D_1,D_2\Subset X$, $\Pi_k:W^{-2m}_{\comp}(D_1,\bigwedge^{0,q}X\otimes L^k)\to W^{2m}_{\loc}(D_2,\bigwedge^{0,q}X\otimes L^k)$ is continuous, for any $m\in \N$.

For any $u\in A^{0,q}_c(D_1,L^k)$, $\square_k^{(q)}\Pi_ku=0$. Therefore, by elliptic estimate, for any $m\in \N$, any $\rho\in C^\infty_c(D_2,[0,1])$, there exists a constant $C=C(m,D_2)>0$ such that for any $u\in A^{0,q}_c(D_1,L^k)$, we have
\[
\|\rho\Pi_k^{(q)}u\|_{W^{2m},\omega,k\phi}\leq C( \|\widetilde{\rho}(\square_k^{(q)})^m \Pi_k^{(q)} u\|_{\omega,k\phi}+\|\widetilde{\rho}\Pi_k^{(q)} u\|_{\omega,k\phi})=C\|\widetilde{\rho}\Pi_k^{(q)}u\|_{\omega,k\phi},
\]
for some $\widetilde{\rho}\in C^\infty_c(D_2,[0,1])$ with $\supp \widetilde{\rho}\supset \supp \rho$. Also, since $\Pi_k^{(q)}$ is orthogonal projection and $0\leq \widetilde{\rho}\leq 1$, we have $\|\widetilde{\rho}\Pi_k^{(q)}u\|_{\omega,k\phi}\leq \|u\|_{\omega,k\phi}$. Therefore, by density argument, we see that 
\[
\Pi_k^{(q)}:W^0_{\comp}(D_1,\bigwedge\nolimits^{0,q}X\otimes L^k)\to W^{2m}_{\loc}(D_2,\bigwedge\nolimits^{0,q}X\otimes L^k)
\] 
is a continuous operator. On the other hand, $\Pi_k^{(q)}=(\Pi_k^{(q)})^*$ is also a continuous operator from $W^{-2m}_{\comp}(D_1,\bigwedge^{0,q}X\otimes L^k)$ to $W^0_{\loc}(D_2,\bigwedge^{0,q}X\otimes L^k)$, for $\Pi_k^{(q)}$ is self-adjoint and $D_1,D_2$ are arbitrary. Thus, for any $v\in W^{-2m}_{\comp}(D_1,\bigwedge^{0,q}X\otimes L^k)$, we choose $v_j\in A^{0,q}_c(D_1,L^k)$, for $j\in \N$, so that $v_j\to v$ in $W^{-2m}(D_1,\bigwedge^{0,q}X\otimes L^k)$. Applying above estimate, we get
\[
\|\rho\Pi_k^{(q)}(v_j-v_l)\|_{W^{2m},\omega,k\phi}\lesssim \|\widetilde{\rho} \Pi_k^{(q)}(v_j-v_l)\|_{\omega,k\phi}\lesssim \|\widetilde{\rho}(v_j-v_l)\|_{W^{-2m},\omega,k\phi}.
\]
Therefore, we see that $\{\rho\Pi_k^{(q)}v_j\}_{j=1}^\infty$ is a Cauchy sequence in $W^{2m}(D_2,\bigwedge^{0,q}X\otimes L^k)$. By completeness of Sobolev space, there exists $w\in W^{2m}(D_2,\bigwedge\nolimits^{0,q}X\otimes L^k)$, $\lim_{j\to\infty}\rho\Pi_k^{(q)}v_j=w$. In conclusion, 
\[
\Pi_k^{(q)}:W^{-2m}_{\comp}(D_1,\bigwedge\nolimits^{0,q}X\otimes L^k)\to W^{2m}_{\loc}(D_2,\bigwedge\nolimits^{0,q}X\otimes L^k)
\]
is continuous, for any $m\in \N$. As a result, $\Pi_k^{(q)}$ is a smoothing operator.
\end{proof}
In view of Fact \ref{thm:Schwartz Kernel Theorem}, $K_k(z,w)\in C^\infty(X\times X,(L^k\otimes \bigwedge^{0,q}X)\boxtimes (L^k\otimes \bigwedge^{0,q}X)^*)$. Hence, by Schwartz kernel Theorem, we obtain \textbf{reproducing property}:
\begin{equation}\label{reproducing property}
(\Pi_{k}^{(q)}u)(z)=\int_X K_{k}^{(q)}(z,w)(u(w))d\nu_X(w),\quad \forall u\in L^2_{0,q}(X,L^k).
\end{equation}
Now, we let $\{\Psi_j(z;k)\}_{j=1}^{d_k}$ be an orthonormal basis for $\mathscr{H}^q(X,L^k)$ with respect to $(\cdot|\cdot)_{\omega,k\phi}$ on $L^2_{0,q}(X,L^k)$, where $d_k\in \N_0\cup\{\infty\}$. By definition, for any $u\in \mathscr{H}^{q}(X,L^k)$, we can write
\begin{equation}\label{expansion into ONB}
u(z)=\sum_{j=1}^{d_k} (u(w)|\Psi_j(w;k))_{\omega,k\phi}\Psi_j(z;k)
\end{equation}
\begin{lem} Let $\Psi_j^*$ be the dual section given by $\Psi_j^*(u)=\la u|\Psi_j\ra_{\omega,k\phi}$. Then 
\begin{equation}\label{BK as expansion into ONB}
K^{(q)}_k(z,w)=\sum_{j=1}^{d_k} \Psi_j(z;k)\otimes \Psi_j^*(w;k).
\end{equation}
\end{lem}
\begin{proof} If $d_k<\infty$, then the sum of course converges and 
\[
\sum_{j=1}^{d_k} \Psi_j(z;k)\Psi_j^*(w;k)(u(w))=\sum_{j=1}^{d_k}\Psi_j(z;k)\la u(w)|\Psi_j(w;k)\ra_{\omega,k\phi}.
\]
Hence, we define a linear operator $\widehat{\Pi}_{k\phi}^{(q)}:L^2_{0,q}(X,L^k)\to \mathscr{H}^{q}(X,L^k)$ by
\begin{align*}
(\widehat{\Pi}_{k\phi}^{(q)}u)(z)&=\sum_{j=1}^{d_k}\Psi_j(z;k)\int_X\la u(w)|\Psi_j(w;k)\ra_{\omega,k\phi}\omega_n(w)\\
&=\sum_{j=1}^{d_k}(u(w)|\Psi_j(w;k))_{\omega,k\phi}\Psi_j(z;k).
\end{align*}
By \eqref{expansion into ONB}, we know that $(\widehat{\Pi}_{k\phi}^{(q)})^2=\widehat{\Pi}_{k\phi}^{(q)}$ and 
\[
(\widehat{\Pi}_{k}^{(q)}u|v)_{\omega,k\phi}=\sum_{j=1}^{d_k}(u|\Psi_j)_{\omega,k\phi}(\Psi_j|v)_{\omega,k\phi}=\sum_{j=1}^{d_k}(u|\Psi_j)_{\omega,k\phi}\overline{(v|\Psi_j)_{\omega,k\phi}}=(u|\widehat{\Pi}_{k}^{(q)}v)_{\omega,k\phi}.
\]
By \cite[Theorem 2, p.83]{Yosida} and the uniqueness of Scwhartz kernel, we know that $\Pi_{\omega,k\phi}^{(q)}=\widehat{\Pi}_{\omega,k\phi}^{(q)}$ and thus $K^{(q)}_k(z,w)=\sum_{j=1}^{d_k} \Psi_j(z;k)\otimes \Psi_j^*(w;k)$ in this case.

For $d_k=\infty$, we define a linear operator $\widehat{\Pi}_N:L^2_{0,q}(X,L^k)\to \mathscr{H}^{q}(X,L^k)$ by
\[
(\widehat{\Pi}_Nu)(z)=\sum_{j=1}^N \Psi_j(z) \int_X \langle u(w)|\Psi_j(w) \rangle_{\omega,k\phi} \omega_n(w)=\sum_{j=1}^N (u|\Psi_j)_{\omega,k\phi}\Psi_j(z).
\]
By \eqref{expansion into ONB}, for each $u\in \mathscr{H}^q(X,L^k)$, $\|u-\widehat{\Pi}_Nu\|_{\omega,k\phi}\to 0$ as $N\to \infty$. Now, for any compact subset $K\subset X$, by G\r{a}rding inequality and $u-\widehat{\Pi}_Nu\in \mathscr{H}^q(X,L^k)$, we know that
\[
\|u-\widehat{\Pi}_Nu\|_{K,W^{2l},\omega,k\phi}\leq C\|u-\widehat{\Pi}_Nu\|_{\omega,k\phi},
\]
for any $l\in \N$. Hence, combining with Sobolev embedding, we know that $\widehat{\Pi}_Nu\to u$ in $L^2$-norm implies $\widehat{\Pi}_Nu$ converges to $u$ locally uniformly in all derivatives. Thus, $\widehat{\Pi}_{k\phi}u:=\lim_{N\to \infty} \widehat{\Pi}_Nu=u=\Pi_{k}u$, for any $u\in \mathscr{H}^q(X,L^k)$. Hence, $\widehat{\Pi}_{k\phi}^2=\widehat{\Pi}_{k\phi}$. Furthermore, for each $N\in \N$,
\[
\|\widehat{\Pi}_N\|=\sup_{u\in \mathscr{H}^{q}(X,L^k),u\neq 0}\frac{\|\widehat{\Pi}_Nu\|_{\omega,k\phi}}{\|u\|_{\omega,k\phi}}=\sup_u\left(\frac{\sum_{j=1}^N |(u|\Psi_j)_{\omega,k\phi}|^2}{\sum_{j=1}^\infty |(u|\Psi_j)_{\omega,k\phi}|^2}\right)^{1/2}\leq 1
\]
Hence, $\widehat{\Pi}_{k\phi}$ is a bounded operator and $\|\widehat{\Pi}_{k\phi}\|\leq 1$.  By \cite[Theorem 3, p.84]{Yosida}, we know that $\widehat{\Pi}_{k\phi}=\Pi_{k}$. By uniqueness of Schwartz kernel, $K_{k}(z,w)=\sum_{j=1}^{\infty} \Psi_j(z;k)\otimes \Psi_j^*(w;k)$.
\end{proof}
 
\chapter{Asymptotic Expansion of Bergman Kernel}\label{chapter:Asymptotic Expansion of Bergman Kernel}
In this chapter, we prove the existence of pointwise asymptotic of Bergman kernel near a point in $X(0)$  under the local spectral gap condition.
\section{Notations and Set-Up}\label{Notation and Set-Up}
Let $X$ be a complex manifold with $\dim_\C X=n$, $\omega$ be a hermitian metric on $X$, and $L\to X$ be a holomorphic line bundle with hermitian metric $h$ on $L$. Let $\Pi_{k}^{(q)}:L^2_{0,q}(X,L^k)\to \mathscr{H}^q(X,L^k)$ be  Bergman projection and $K_{k}^{(q)}$ be Bergman kernel. Let $s$ be local holomorphic trivialization of $L$ on an open subsets of $X$ and $|s|_h^2=e^{-2\phi}$, where $\phi$ is the local weight of $h$ with respect to $s$.

Let $\eta$ be a $(0,1)$-form. We denote $\epsilon(\eta):=\eta\wedge\cdot:\bigwedge^{0,q}T^*_xX\to\bigwedge^{0,q+1}T^*_xX$ be the wedging $\eta$ from the left, and $\iota(\eta)$ be its adjoint with respect to $\langle\cdot|\cdot\ra_\omega$. Hence, for $\eta_1,\eta_2\in A^{0,1}(X)$,
\[
\epsilon(\eta_1)\iota(\eta_2)+\iota(\eta_2)\epsilon(\eta_1)=\la \eta_1|\eta_2\ra_{\omega} id.
\]
Let $e^1(z),\dots,e^n(z)$ be an orthonormal frame for $\bigwedge^{0,1}X$ over $U$, $Z_1,\dots,Z_n$ be its dual basis for $T^{0,1}X$. We can write Cauchy--Riemann operator $\db$ on $A^{0,q}(U,L^k)$ as
\begin{equation}
\db(s^k\otimes \alpha)=s^k\otimes \sum_{j=1}^n \left(\epsilon(e^j)Z_j+\epsilon(\db e^j)\iota(e^j)\right)\alpha.
\end{equation}
Thus, its formal adjoint with respect to the scalar product $(\cdot|\cdot)_{\omega,k\phi}$ is given by
\begin{equation}
\db^*(s^k\otimes \alpha)=s^k\otimes \sum_{j=1}^n\left(\iota(e^j)(Z_j^*+2k\overline{Z_j}(\phi))+\epsilon(e_j)\iota(\db e^j)\right)\alpha,
\end{equation}
where $Z_j^*$ is the formal adjoint of $Z_j$ with respect to the inner product $(\alpha|\beta)_\omega:=\int_X \la \alpha|\beta\ra_\omega \omega_n$ on $A^{0,q}_c(X)$. To put $\db$ and $\db^*$ in more symmetric form, we make the following identification.
\begin{equation}\label{unitary identification}
\begin{aligned}
A^{0,q}(U,L^k)\to A^{0,q}(U)&,\quad u=s^k\otimes \alpha\mapsto \widetilde{\alpha}:=\alpha e^{-k\phi}\\
A^{0,q}(U)\to A^{0,q}(U,L^k)&,\quad \beta\mapsto s^k\otimes e^{k\phi}\beta.
\end{aligned}
\end{equation}
This is a local unitary identification since for $u,v\in L^2_{0,q}(X,L^k)\cap \E'(X,\bigwedge^{0,q}X\otimes L^k)$,
\[
(u|v)_{\omega,k\phi}=\int_U \la\alpha|\beta\ra_{\omega} e^{-2k\phi}\omega_n=\int_U \la e^{-k\phi}\alpha|e^{-k\phi}\beta\ra_{\omega}=(\widetilde{\alpha}|\widetilde{\beta})_{\omega}, 
\]
where $u=s^k\otimes \alpha,v=s^k\otimes \beta$. Then under this unitary identification, we get 
\begin{equation}\label{localized dbar}
\db(s^k\otimes e^{k\phi}\alpha)=s^k\otimes e^{k\phi}\db_{k,s}\alpha,
\end{equation}
where $\alpha\in A^{0,q}(U)$ and 
\begin{equation}\label{localized dbar:defin}
\db_{k,s}=\sum_{j=1}^n \left(\epsilon(e^j)\otimes (Z_j+kZ_j(\phi))+\epsilon(\db e^j)\iota(e^j)\right)=\db+k\epsilon(\db \phi).
\end{equation}
The formal adjoint $\db_{k,s}^*$ with respect to the local scalar product $(\cdot|\cdot)_\omega$ is given by
\begin{equation}\label{localized dbar*:defin}
\db_{k,s}^*=\sum_{j=1}^n \left(\iota(e^j)\otimes (Z_j^*+k\overline{Z_j(\phi)})+\epsilon(e^j)\iota(\db e^j)\right)=\db^*+k\iota(\db \phi)
\end{equation}
and satisfies
\begin{equation}\label{localized dbar*}
\db^*(s^k\otimes e^{k\phi}\beta)=s^k\otimes e^{k\phi}\db^*_{k,s}\beta,\quad \beta\in A^{0,q+1}(U).
\end{equation}
We call $\db_{k,s}$ the \textbf{localized Cauchy-Riemann operator} with respect to $s$. Then the \textbf{localized Kodaira Laplacian} with respect to $s$ is then defined by
\begin{equation}\label{localized Laplacian}
\square_{k,s}^{(q)}:=\db_{k,s}^*\db_{k,s}+\db_{k,s}\db_{k,s}^*.
\end{equation}
Of course, from \eqref{localized dbar}, \eqref{localized dbar*}, we have
\begin{equation}
\square_{\omega,k\phi}^{(q)}(s^k\otimes e^{k\phi}\alpha)=s^k\otimes e^{k\phi}\square_{k,s}^{(q)}\alpha,\quad \alpha\in A^{0,q}(U).
\end{equation}
The \textbf{localized Bergman projection} $\Pi^{(q)}_{k,s}:L^2_{0,q}(U,\omega)\cap \mathcal{E'}(U,\bigwedge^{0,q}X)\to A^{0,q}(U)$ is given by
\begin{equation}\label{localized Bergman projection}
\Pi_{k,s}^{(q)}\alpha=e^{-k\phi}s^{-k}\Pi^{(q)}_k(e^{k\phi}\alpha\otimes s^k).
\end{equation}
Thus, we see that $\Pi_{k,s}^{(q)}:L^2_{0,q}(U,\omega)\cap \mathcal{E}'(U,\bigwedge^{0,q}X)\to \ker \square_{k,s}^{(q)}$. Let $K_{k,s}^{(q)}$ be the Schwartz kernel of $\Pi_{k,s}^{(q)}$, i.e.,
\begin{equation}\label{localized Bergman kernel}
(\Pi_{k,s}^{(q)}\alpha)(z)=\int_U K_{k,s}^{(q)}(z,w)(\alpha(w))\omega_n(w).
\end{equation}
The main result  in this chapter is the following 
\begin{thm}[=Theorem \ref{Pointwise Asymptotic of Bergman Kernel Function}]\label{thm:Pointwise Asymptotic of Bergman Kernel Function}
Suppose $X(0)\neq \emptyset$, say $x\in X(0)$. For any $D\subset X(0)$ of $x$ satisfying the spectral gap condition (cf. \eqref{defin:spectral gap condition}), there exists a trivializing open set $U\Subset D$ and a holomorphic coordinate $z$ on $U$ centered at $x$ so that on $U$,  we have 
\[
\rho(z)K_{k,s}(z,w)\chi_k(w)\in \h{S}^n_{\cl}(\C^n\times \C^n),
\]
where $\rho\in C^\infty_c(U)$, $\chi_k(z):=\chi(8k^{1/2-\epsilon}z)$, $\chi\in C^\infty_c(\C^n)$ satisfying
\[
\supp\chi\subset B_1(0),\quad \chi=1\text{ on }B_{1/2}(0),\quad \rho=1\text{ near } 0,
\]
and $\epsilon\in (0,\frac{1}{6})$.\end{thm}
To start, we give a local expression for localized Kodaira Laplacian.
\begin{prop}\label{local expression of localized Kodaira Laplacian}
The localized Kodaira Laplacian $\square_{k,s}^{(q)}$ satisfies
\begin{equation}
\begin{aligned}
&\square_{k,s}^{(q)}=\db_{k,s}^*\db_{k,s}+\db_{k,s}\db_{k,s}^*\\
=&\sum_{j=1}^n 1\otimes (Z_j^*+k\overline{Z_j(\phi)})(Z_j+kZ_j(\phi))\\
+&\sum_{j,l=1}^n \epsilon(e^j)\iota(e^l)\otimes \left[Z_j+kZ_j(\phi),Z_l^*+k\overline{Z_l(\phi)}\right]\\
+&O(1)(Z+kZ(\phi))+O(1)(Z^*+k\overline{Z(\phi)})+O(1),
\end{aligned}
\end{equation}
where $Z+kZ(\phi)$ indicates a remainder term of the form $\sum_{j=1}^n a_j(z)(Z_j+kZ_j(\phi))$ and $a_j(z)$ are some $k$-independent smooth function, and similarly for $Z^*+\overline{Z(\phi)}$. Also, $O(1)$ indicates some zero order differential operators which are independent of $k$.
\end{prop}
\begin{proof}
By direct computation,
\begin{equation}\label{Kodaira Laplacian expansion}
\begin{aligned}
\square_{k,s}^{(q)}=&\db_{k,s}^*\db_{k,s}+\db_{k,s}\db_{k,s}^*\\
=\sum_{j,l=1}^n (&\epsilon(e^j)\otimes (Z_j+kZ_j(\phi))(\iota(e^l)\otimes (Z_l^*+k\overline{Z_l(\phi)})))\\
+&(\iota(e^l)\otimes (Z_l+kZ_l(\phi)))(\epsilon(e^j)\otimes (Z_j+kZ_j(\phi)))\\
+&\epsilon(\db e^j)\iota(e^j)\iota(e^l)\otimes (Z_l^*+k\overline{Z_l(\phi)})\\
+&\iota(e^l)\epsilon(\db e^j)\iota(e^j)\otimes (Z_j+k\overline{Z_j(\phi)})\\
+&\epsilon(\db e^j)\iota(e^j)\epsilon(e^l)\iota(\db e^l)+\epsilon(e^l)\iota(\db e^l)\epsilon(\db e^j)\iota(e^j).
\end{aligned}
\end{equation}
Now, we combine the first two terms as
\begin{align*}
&(\epsilon(e^j)\iota(e^l)+\iota(e^l)\epsilon(e^j))((Z_j+kZ_j(\phi))(Z_l^*+k\overline{Z_l(\phi)}))\\
+&\epsilon(e^j)\iota(e^l)[Z_j+kZ_j(\phi),Z_l^*+k\overline{Z_l(\phi)}]\\
=&\la e^j|e^l\ra ((Z_j+kZ_j(\phi))(Z_l^*+k\overline{Z_l(\phi)})+\epsilon(e^j)\iota(e^l)[Z_j+kZ_j(\phi),Z_l^*+k\overline{Z_l(\phi)}]\\
=&\delta_{jl}(Z_j+kZ_j(\phi))(Z_l^*+k\overline{Z_l(\phi)})+\epsilon(e^j)\iota(e^l)[Z_j+kZ_j(\phi),Z_l^*+k\overline{Z_l(\phi)}].
\end{align*}
\end{proof}

To end this section, we define the concept of $k$-negligible kernel or $k$-negligible operators which we will use throughout the chapter.
\begin{defin}\label{defin:k-negligible}A $k$-dependent continuous linear operator $\mathcal{A}_k:L^2_{0,q}(X,L^k)\to L^2_{0,q}(X,L^k)$  is \textbf{$k$-negligible} if it is smoothing for sufficiently large $k$ and for any $\alpha,\beta\in \N_0^{2n}$, any $N\in \N$, there exists a $k$-independent constant $C_{\alpha,\beta,N,L}>0$ such that the smooth kernel $A_k(x,y)$ of $\mathcal{A}_k$ satisfying
\begin{equation}\label{k-negligible}
\left|\p_x^\alpha\p_x^\beta A_{k,s,t}(x,y)\right|\leq C_{\alpha,\beta,N,L}k^{-N}, \text{for } k\gg1,
\end{equation}
locally uniformly on any compact subset $L\subset U\times V$, where $s,t$ are local holomorphic trivialization of $L$ over $U,V$, respectively, and $A_k(x,y)=A_{k,s,t}(x,y)s^k(x)\otimes (t^k)^*(y)$, where $(t^k)^*$ is the metric dual of $t^k$. If so, we denote $\mathcal{A}_k\equiv 0\mod O(k^{-\infty})$ or $A_k\equiv 0\mod O(k^{-\infty})$.
\end{defin}
Notice that the condition of $k$-negligible is independent of the choice of local trivializations $s,t$ and local coordinates $x,y$. Also, by Sobolev embedding, $\mathcal{A}_k$ is $k$-negligible if and only if $\mathcal{A}_k$ extends to a operator from $W_{\comp}^{r}(U,\bigwedge^{0,q}X\otimes L^k)$ to $W_{\loc}^{r+M}(V,\bigwedge^{0,q}X\otimes L^k)$ with operator norm $O(k^{-N})$, for any $r\in \R$, $M,N\in \N$.
\section{Approximate Bergman Kernel and Semi-Classical \texorpdfstring{$L^2$}{}-Estimates}\label{sec:Approximate Bergman Kernel}
In this section, we define approximate Bergman kernel on $\C^n$ and develop $L^2$-estimates for this model case (cf. \cite{HormanderL2}, \cite{HormanderSCV}, and also \cite{Chen-Shaw}) but with semi-classical parameter $k$. Suppose $x\in X(0)$, let $U\subset X(0)$ be a local trivialization open set and $x\in U$. We choose a local coordinate $(U,z)$ centered at $x$ and a local holomorphic trivialization $s$ of $L$ on $U$ so that Lemma \ref{Canonical Coordinate} holds. That is,
\begin{equation}
\begin{aligned}
&\phi(z)=\phi_0(z)+O(|z|^3),\quad \phi_0(z)=\sum_{i=1}^n \lambda_{i,x}|z^i|^2,\quad \lambda_{i.x}>0,\\
&\omega(z)=\omega_0(z)+O(|z|),\quad \omega_0(z)=\frac{\ii}{2}\sum_{i=1}^n dz^i\wedge d\bar{z}^i,\quad z\in U.
\end{aligned}
\end{equation}
 Identifying $U$ as a bounded domain in $\C^n$, we extend $\phi$ and $\omega$ by $\widehat{\phi}$ and $\widehat{\omega}$ to whole $\C^n$ by
\begin{equation}\label{extension of local metric data}
\widehat{\phi}=\phi_0+\underbrace{\theta_k(\phi-\phi_0)}_{\phi_1},\quad \widehat{\omega}=\omega_0+\underbrace{\theta_k(\omega-\omega_0)}_{\omega_1},\text{ and}\quad \theta_k(z)=\theta(k^{1/2-\epsilon}z),
\end{equation}
where $\theta\in C^\infty_c(\C^n)$ is a cut-off function such that $\theta=1$ on $B_{1/2}(0)$ and $\supp\theta\subset B_1(0)$. Thus, $\phi_1\in C^\infty_c(U,\R)$, for any small $\epsilon>0$. Since $U\subset X(0)$, we know that $\phi$ is  strictly plruisubharmonic, i.e., there exists $C>0$ such that $\frac{\p^2\phi}{\p z^j\p\bar{z}^l}(z)\xi^j\overline{\xi}^l\geq C|\xi|^2$, for any $\xi\in \C^n\setminus \{0\}$, any $z\in U$. Thus, for sufficiently large $k$,  $\h{\phi}$ is strictly plurisubharmonic on $\C^n$. Let $\h{\omega}_n:=\h{\omega}^n/n!=\lambda(z)dm(z)$. Then $\lambda(z)=1$ outside $B_{k^{\epsilon-1/2}}(0)$. In other words, 
\begin{align*}
&\h{\omega}=\omega,\quad \h{\phi}=\phi \text{ on  }V_k:=B_{\frac{1}{2}k^{\epsilon-1/2}}(0)\\
&\h{\omega}=\omega_0,\quad \h{\phi}=\phi_0 \text{ on } |z|>k^{\epsilon-1/2}.
\end{align*}

We then consider $L^2$-space $L^2_{0,q}(\C^n,\h{\omega})$ which is the completion of $A^{0,q}_c(\C^n)$ with respect to the  $L^2$-norm given by 
\[
(f|g)_{\h{\omega}}:=\int_{\C^n}\la f|g\ra_{\h{\omega}}(z)\lambda(z)dm(z),\quad f,g\in A^{0,q}_c(\C^n).
\] 
\begin{rmk}
In the remaining of this section, unless otherwise stated, we will denote $(\cdot|\cdot)_{\h{\omega}}$ simply by $(\cdot|\cdot)$ for the sake of brevity.
\end{rmk}
We now define the \textbf{deformed Cauchy--Riemann operator} $\db_{k\h{\phi}}:=\db+k\epsilon(\db\h{\phi}):A^{0,q}(\C^n)\to A^{0,q+1}(\C^n)$ and its formal adjoint $\db_{k\h{\phi}}^*=\db^*+k\iota(\db\h{\phi})$ with respect to $(\cdot|\cdot)_{\h{\omega}}$. Hence, as before, the {\bf deformed Kodaira Laplacian} is then defined by
\begin{equation}
\triangle_{\h{\omega},\h{\phi}}^{(q)}=\db_{k\h{\phi}}\db_{k\h{\phi}}^*+\db^*_{k\h{\phi}}\db_{k\h{\phi}}.
\end{equation}
Note that the analogous formula \eqref{localized dbar:defin} and \eqref{localized dbar*:defin} still hold:
\begin{equation}\label{deformed dbar/*}
\begin{aligned}
&\db_{k\h{\phi}}=\sum_{j=1}^n \left(\epsilon(e^j)\otimes (Z_j+kZ_j(\h{\phi}))+\epsilon(\db e^j)\iota(e^j)\right)\\
&\db_{k\h{\phi}}^*=\sum_{j=1}^n \left(\iota(e^j)\otimes (Z_j^*+k\overline{Z_j(\h{\phi})})+\epsilon(e^j)\iota(\db e^j)\right)
\end{aligned}.
\end{equation}
Also, $\db_{\h{\omega},k\h{\phi}}=\db_{k,s}$, $\db_{\h{\omega},k\h{\phi}}^*=\db_{k,s}^*$, and $\triangle_{\h{\omega},k\h{\phi}}^{(q)}=\square_{k,s}^{(q)}$ on $V_k$. For $z\in \C^n\setminus U$, $\widehat{\omega}=\omega_0$, and thus $b_j^i(z)=\delta_j^i=c_j^i(z)$. 

From Proposition \ref{local expression of localized Kodaira Laplacian}, we see that for $q\geq 1$, $\alpha\in A^{0,q}_c(\C^n)$,
\begin{align*}
\triangle_{\h{\omega},k\h{\phi}}^{(q)}\alpha&=\sum_{j=1}^n 1\otimes (Z_j^*+k\overline{Z_j(\h{\phi})})(Z_j+kZ_j(\h{\phi}))\alpha\\
&+k\sum_{j,l=1}^n \epsilon(e^j)\iota(e^l)\otimes (Z_j\overline{Z_l}-Z_l^*Z_j)(\h{\phi})\alpha\\
&+O(1)(Z+kZ(\h{\phi}))\alpha+O(1)(Z^*+\overline{kZ(\h{\phi}}))\alpha+O(1)\alpha.
\end{align*}
Therefore, we get
\begin{align*}
(\triangle_{\h{\omega},k\h{\phi}}^{(q)}\alpha|\alpha)&=\left(\sum_{j=1}^n 1\otimes (Z_j^*+k\overline{Z_j(\h{\phi})})(Z_j+kZ_j(\h{\phi}))\alpha\middle|\alpha\right)\\
&+k\left(\sum_{j,l=1}^n \epsilon(e^j)\iota(e^l)\otimes (Z_j\overline{Z_l}-Z_l^*Z_j)(\h{\phi})\alpha\middle|\alpha\right)\\
&+\left(O(1)(Z+kZ(\h{\phi}))\alpha+O(1)\alpha\middle|\alpha\right)+\left(\alpha|O(1)(Z+kZ(\h{\phi}))\alpha\right)\\
&\geq \|(Z+kZ(\h{\phi}))\alpha\|^2+k\left(\sum_{j,l=1}^n \epsilon(e^j)\iota(e^l)\otimes (Z_j\overline{Z_l}-Z_l^*Z_j)(\h{\phi})\alpha\middle|\alpha\right)\\
&-\left|\left(O(1)(Z+kZ(\h{\phi}))\alpha|\alpha\right)\right|-\left|\left(\alpha|O(1)(Z+k(\h{\phi}))\alpha\right)\right|-\left|\left(O(1)\alpha\middle|\alpha\right)\right|,
\end{align*}
where $Z+kZ(\h{\phi})=\sum_{j=1}^n Z_j+kZ_j(\h{\phi})$. By Cauchy--Schwartz inequality,
\begin{align*}
&\left|\left(O(1)(Z+kZ(\h{\phi}))\alpha|\alpha\right)\right|\leq \frac{1}{2}\left(\epsilon \left\|O(1)(Z+kZ(\h{\phi}))\alpha\right\|^2+\frac{1}{\epsilon}\|\alpha\|^2\right),\\
&\left|\left(O(1)\alpha|\alpha\right)\right|\leq \frac{1}{2}\left(\|O(1)\alpha\|^2+\|\alpha\|^2\right).
\end{align*}
By the construction, these $O(1)$ terms are supported in a compact set and uniformly bounded in $k$. Hence, there exists a constant $C'>0$ independent of $k$ so that
\[
\|O(1)\alpha\|^2\leq C'\|\alpha\|^2,\quad \|O(1)(Z+kZ(\h{\phi}))\alpha\|^2\leq C'\|(Z+kZ(\h{\phi}))\alpha\|^2
\]
and we choose $\epsilon>0$ so that $\epsilon C'<1$. Also, since $\h{\phi}$ is strictly plurisubharmonic, there exists a $k$-independent constant $C_0>0$ so that
\[
\left(\sum_{j,l=1}^n \epsilon(e^j)\iota(e^l)\otimes (Z_j\overline{Z_l}-Z_l^*Z_j)(\h{\phi})\alpha\middle|\alpha\right)\geq C_0\|\alpha\|^2.
\]
Combining these estimates, we obtain
\begin{align*}
(\triangle_{\h{\omega},k\h{\phi}}^{(q)}\alpha|\alpha)&\geq \left(1-\frac{\epsilon}{C'}\right) \|(Z+kZ(\h{\phi}))\alpha\|^2\\
&+\left(C_0k-\frac{1}{\epsilon}-\frac{C'}{2}-1\right)\|\alpha\|^2,\quad \forall \alpha\in A^{0,q}_c(\C^n).
\end{align*}
Therefore, for sufficiently large $k$, there exists a constant $C$ independent of $k$ so that
\begin{equation}\label{spectral gap for deformed dbar}
(\triangle^{(q)}_{\h{\omega},k\h{\phi}}\alpha|\alpha)\geq Ck\|\alpha\|^2,\quad \forall \alpha\in A^{0,q}_c(\C^n),\quad q\geq 1.
\end{equation}
We next consider Gaffney extension of $\triangle^{(q)}_{\h{\omega},k\h{\phi}}$ and show that \eqref{spectral gap for deformed dbar} holds for any $\alpha\in \dom \triangle_{\h{\omega},k\h{\phi}}^{(q)}$, for $q\geq 1$. Similar to section \ref{Gaffney Extension of Kodaira Laplacian}, Gaffney extension for deformed Kodaira Laplacian on $(0,q)$-forms is given by
\begin{align*}
\dom \triangle_{\h{\omega},k\h{\phi}}^{(q)}&:=\{\alpha\in \dom(\h{S}_q)\cap \dom(\h{S}_{q-1}^\dagger):S_q\alpha\in \dom(\h{S}_q^\dagger),\h{S}_{q-1}^\dagger\alpha\in \dom \h{S}_{q-1}\}\\
\triangle^{(q)}_{\h{\omega},k\h{\phi}}\alpha&:=\h{S}_q^\dagger\h{S}_q\alpha+\h{S}_{q-1}\h{S}_{q-1}^\dagger\alpha,\quad \forall \alpha\in \dom \triangle_{\h{\omega},k\h{\phi}}^{(q)}.
\end{align*}
where $\h{S}_q$ is the maximal extension of $\db_{k\h{\phi}}:A^{0,q}(\C^n)\to A^{0,q+1}(\C^n)$ and $\h{S}_{q-1}^\dagger$ is the adjoint of $\h{S}_{q-1}$. As in section \ref{Gaffney Extension of Kodaira Laplacian}, we know that the Gaffney extension $\triangle_{\h{\omega},k\h{\phi}}^{(q)}$ is a self-adjoint, non-negative operator. Moreover, since $\triangle_{\h{\omega},k\h{\phi}}^{(q)}$ has the same principal symbol as usual Laplacian on $\C^n$, we know that it is elliptic. By the same argument as in section \ref{Bergman Kernel}, we know that the $L^2$-projection 
\begin{equation}\label{Approximate Bergman Kernel}
\mathcal{P}^{(q)}_{\h{\omega},k\h{\phi}}:L^2_{0,q}(\C^n,\h{\omega}_n)\to \ker \triangle_{\h{\omega},k\h{\phi}}^{(q)}
\end{equation}
 is a smooth operator and thus has smooth Schwartz kernel $P_{\h{\omega},k\h{\phi}}^{(q)}$. We call $P_{\h{\omega},k\h{\phi}}^{(q)}$ the \textbf{approximate Bergman kernel for $(0,q)$-forms}. Now, we show 
\begin{lem}[Approximation Lemma]\label{Approximation Lemma} Let $\alpha\in \dom \h{S}_q\cap \dom \h{S}^\dagger_{q-1}\subset L^2_{0,q}(\C^n,\h{\omega}_n)$. Then there exists a sequence $\{\alpha_j\}_{j=1}^\infty\subset A^{0,q}_c(\C^n)$ such that 
\[
\|\alpha_j-\alpha\|,\quad \|\db_{k\h{\phi}}\alpha_j-\h{S}_q\alpha\|,\quad\|\db_{k\h{\phi}}^*\alpha_j-\h{S}_{q-1}^\dagger\alpha\|\to 0,\quad \text{as }j\to\infty.
\]
\end{lem}
\begin{proof} First of all, we show that if $\chi\in C^\infty_c(\C^n,\R)$ is a test function, then for any $\alpha\in \dom(\h{S}_q)\cap \dom(\h{S}_{q-1}^\dagger)$, $\chi\alpha\in \dom(\h{S}_q)\cap \dom(\h{S}_{q-1}^\dagger)$. For $\alpha\in\dom(\h{S}_q)$, we first claim that the Leibniz rule $\db_{k\h{\phi}}(\chi \alpha)=\chi \db_{k\h{\phi}}\alpha+\db_{k\h{\phi}}\chi\wedge \alpha$ holds in the distribution sense. To see this, observe that for any $\beta\in A^{0,q+1}_c(\C^n)$, $\gamma\in A^{0,q}_c(\C^n)$,
\begin{align*}
&(\gamma|\db_{k\h{\phi}}^*(\chi \beta))=(\db_{k\h{\phi}}\gamma|\chi \beta)=(\chi \db_{k\h{\phi}}\gamma|\beta)\\
=&(\db_{k\h{\phi}}(\chi\gamma )-\db_{k\h{\phi}}\chi\wedge \gamma|\beta)=(\gamma|\chi\db_{k\h{\phi}}^*\beta)-(\gamma| \iota_{\db_{k\h{\phi}}\chi}\beta),
\end{align*}
where $\iota(\db_{k\h{\phi}}\chi)$ is the adjoint of $\epsilon(\db_{k\h{\phi}}\chi)$ with respect to $(\cdot|\cdot)$. Hence, $\db_{k\h{\phi}}^*(\chi\beta)=\chi\db_{k\h{\phi}}^*\beta-\iota_{\db_{k\h{\phi}}\chi}\beta$. From this, we deduce that
\begin{align*}
\db_{k\h{\phi}}(\chi \alpha)(\beta)&:=\chi\alpha(\db_{k\h{\phi}}^*\beta)=\alpha(\chi\db_{k\h{\phi}}^*\beta)\\
&=\alpha(\db_{k\h{\phi}}^*(\chi\beta))+\alpha(\iota_{\db_{k\h{\phi}}\chi}\beta)=\chi\db_{k\h{\phi}}\alpha(\beta)+\db_{k\h{\phi}}\chi\wedge\alpha(\beta).
\end{align*}
Hence, $\db_{k\h{\phi}}(\chi\alpha)=\chi\db_{k\h{\phi}}\alpha+\db_{k\h{\phi}}\chi\wedge \alpha$ holds indeed in the distribution sense. Therefore, since $\alpha\in \dom(\h{S}_q)$ and $\chi\in C^\infty_c(\C^n)$, $\db_{k\h{\phi}}\alpha\in L^2_{0,q}(\C^n)$ and $\db_{k\h{\phi}}\chi=\db\chi+k\chi \db\phi\in A^{0,q}_c(\C^n)$, we see that $\db_{k\h{\phi}}(\chi\beta)\in L^2_{0,q+1}(\C^n)$. Next, for $\alpha\in \dom(\h{S}_{q-1}^\dagger)$, we need to show that there exists constant $C>0$ such that
\[
\left|(\db_{k\h{\phi}}u|\chi\alpha)\right|\leq C\|u\|,\quad \forall u\in \dom(\h{S}_{q-1}).
\]
However, $\left|(\db_{k\h{\phi}}u|\chi\alpha)\right|\leq \sup_{z\in \C^n}|\chi(z)| \left|(\db_{k\h{\phi}}u|\alpha)\right|\leq C'\|u\|$. This shows that $\chi\alpha\in \dom(\h{S}_{q-1}^\dagger)$. 

Now, we can choose $\chi$ to belong to some partition of unity with compact supports and decompose $\alpha=\sum_{j=1}^n\chi_j\alpha$. It suffices to approximate each $\chi_j\alpha$ and thus we may assume that $\alpha$ has supports in some compact set $K$. Then we apply the standard regularization technique by convoluting the coefficients of $\alpha$ by the mollifiers $\rho_j(z):=j^{2n}\rho(jz)$, where $\rho$ is the standard modifier on $\C^n$. The result then follows from the classical Lemma of Friedrichs (cf. \cite{Demailly2} Chapter VII, Lemma 3.3).
\begin{fact}[Friedrichs]
Let $L:=\sum_{j=1}^m g_j\frac{\p}{\p x^k}+g$ be a first order differential operator on an open set $U\subset \R^m$ with coefifients $g_k\in C^1(U)$, $g\in C^0(U)$. Then for any for any $v\in L^2(\R^m)\cap \mathcal{E}'(U)$, 
$\lim_{j\to\infty}\|L(v*\rho_j)-(Lv)*\rho_j\|_{L^2}=0$.
\end{fact}
\end{proof}
From above approximation Lemma and \eqref{spectral gap for deformed dbar}, we deduce that for $q\geq 1$,
\[
\|\h{S}_q\alpha\|^2+\|\h{S}_{q-1}^\dagger\alpha\|^2\geq Ck\|\alpha\|^2,\quad \forall \alpha\in \dom \h{S}_q\cap \dom \h{S}_{q-1}^\dagger.
\]
For $\alpha\in \dom \triangle_{\h{\omega},k\h{\phi}}^{(1)}$, we also have $(\triangle_{\h{\omega},k\h{\phi}}^{(q)}\alpha|\alpha)=(\h{S}_q^\dagger\h{S}_q\alpha|\alpha)+(\h{S}_{q-1}\h{S}_{q-1}^\dagger\alpha|\alpha)=\|
\h{S}_q\alpha\|^2+\|\h{S}_{q-1}^\dagger\alpha\|^2$. Thus, this implies that \eqref{spectral gap for deformed dbar} extends to all elements in $\dom \triangle_{\h{\omega},k\h{\phi}}^{(q)}$, i.e.,
\begin{equation}\label{spectral gap for deformed dbar-3}
\left(\triangle_{\h{\omega},k\h{\phi}}^{(q)}\alpha\middle|\alpha\right)\geq Ck\|\alpha\|^2,\quad \forall \alpha\in \dom \triangle_{\h{\omega},k\h{\phi}}^{(q)},\quad q\geq 1.
\end{equation}
From this, we can prove
\begin{cor} The deformed Kodaira Laplacian $\triangle_{\h{\omega},k\h{\phi}}^{(q)}:\dom \triangle_{\h{\omega},k\h{\phi}}^{(q)}\to L^2_{0,q}(\C^n,\h{\omega})$ is bijective and has a bounded inverse.
\end{cor}
\begin{proof}
First, it is clear from \eqref{spectral gap for deformed dbar-3} that $\ker\triangle_{\h{\omega},k\h{\phi}}^{(q)}=0$. For surjectivity, given any $\beta\in L^2_{0,q}(\C^n,\h{\omega})$, we consider the a linear functional on $\im(\triangle_{\h{\omega},k\h{\phi}}^{(q)})$ given by 
\[
\ell_\beta(\triangle_{\h{\omega},k\h{\phi}}^{(q)}\alpha)=(\alpha|\beta),\quad \forall\alpha\in \dom\triangle_{\h{\omega},k\h{\phi}}^{(q)}.
\]
Injectivity of $\triangle_{\h{\omega},k\h{\phi}}^{(q)}$ implies that $\ell_\beta$ is well-defined. Also, \eqref{spectral gap for deformed dbar-3} implies $\|\ell\|\leq\frac{\|\beta\|}{Ck}$. By Hahn--Banach Theorem, $\ell_\beta$ extends to a bounded linear functional  on $\triangle_{\h{\omega},k\h{\phi}}^{(q)}$ with the same norm. By Riesz representation Theorem, there exists $\gamma\in L^2_{0,1}(\C^n,\h{\omega})$ with $\|\gamma\|\leq \frac{\|\beta\|}{Ck}$ such that
\[
(\alpha|\beta)=\ell_\beta(\triangle_{\h{\omega},k\h{\phi}}^{(q)}\alpha)=(\triangle_{\h{\omega},k\h{\phi}}^{(1)}\alpha|\gamma),\quad \forall \alpha\in \dom\triangle_{\h{\omega},k\h{\phi}}^{(q)}.
\]
In other words, $\gamma\in \dom(\triangle_{\h{\omega},k\h{\phi}}^{(q)})^\dagger=\dom\triangle_{\h{\omega},k\h{\phi}}^{(1)}$ and $\beta=\triangle_{\h{\omega},k\h{\phi}}^{(q)}\gamma$. This shows the first assertion.  The second assertion follows from \eqref{spectral gap for deformed dbar-3} that $\|(\triangle_{\h{\omega},k\h{\phi}}^{(q)})^{-1}\|\leq \frac{1}{Ck}$.
\end{proof}
\begin{rmk}
Since $\triangle_{\h{\omega},k\h{\phi}}^{(q)}$ is self-adjoint, $(\triangle_{\h{\omega},k\h{\phi}}^{(q)})^{-1}$ is also self-adjoint. \eqref{spectral gap for deformed dbar-3} implies that $\Spec \triangle_{\h{\omega},k\h{\phi}}^{(q)}\subset \{0\}\cup [Ck,\infty)$. Then we have $\Spec (\triangle_{\h{\omega},k\h{\phi}}^{(q)})^{-1}\subset [0,\frac{1}{Ck}]$.
\end{rmk}
Now, we turn to the $L^2$-existence Theorem for $\db_{k\h{\phi}}$ on $\C^n$. 
\begin{thm} \label{L^2-existence Theorem}
If $\alpha\in L^2_{0,1}(\C^n,\h{\omega})$ with $\db_{k\h{\phi}}\alpha=0$ in the distribution sense, then $u=\db_{k\h{\phi}}^*(\triangle_{\h{\omega},k\h{\phi}}^{(1)})^{-1}\alpha$ solves $\db_{k\h{\phi}}u=\alpha$ and we have the $L^2$-estimate
\[
\|u\|\leq \frac{1}{\sqrt{Ck}}\|\alpha\|.
\]
\end{thm}
\begin{proof}
Since $\beta:=(\triangle_{\h{\omega},k\h{\phi}}^{-1})^{-1}\alpha\in \dom(\triangle_{\h{\omega},k\h{\phi}})$, thus $u=\db_{k\h{\phi}}^*\beta\in \dom\db_{k\h{\phi}}\subset L^2(\C^n)$. The expression of $u$ is legitimate. Then we have
\[
\db_{k\h{\phi}}u=\triangle_{\h{\omega},k\h{\phi}}^{(1)}\beta-\db_{k\h{\phi}}^*\db_{k\h{\phi}}\beta=\alpha-\db_{k\h{\phi}}^*\db_{k\h{\phi}}\beta.
\]
Now, we claim that $\db_{k\h{\phi}}^*\db_{k\h{\phi}}\beta=0$. To see this, we first compute that
\[
\db_{k\h{\phi}}\db^*_{k\h{\phi}}\db_{k\h{\phi}}\beta=\db_{k\h{\phi}}\triangle_{\h{\omega},k\h{\phi}}^{(1)}\beta=\db_{k\h{\phi}}\alpha=0.
\]
Hence, $\db^*_{k\h{\phi}}\db_{k\h{\phi}}\beta\in \ker \db_{k\h{\phi}}$. Also, $\db_{k\h{\phi}}^*\db_{k\h{\phi}}\beta\in \ker \db_{k\h{\phi}}^*$ clearly. Therefore, $\db^*_{k\h{\phi}}\db_{k\h{\phi}}\beta\in \ker \triangle_{\h{\omega},k\h{\phi}}^{(1)}$. From \eqref{spectral gap for deformed dbar-3}, we know that $\ker\triangle_{\h{\omega},k\h{\phi}}^{(1)}=0$ and thus $\db_{k\h{\phi}}u=\alpha$.
As for the last statement, 
\begin{align*}
&\|u\|^2=(\db_{k\h{\phi}}^*(\triangle_{\h{\omega},k\h{\phi}}^{(1)})^{-1}\alpha|\db_{k\h{\phi}}^*(\triangle_{\h{\omega},k\h{\phi}}^{(1)})^{-1}\alpha)\\
=&\left((\triangle_{\h{\omega},k\h{\phi}}^{(1)})^{-1}\alpha|\db_{k\h{\phi}}\db_{k\h{\phi}}^*(\triangle_{\h{\omega},k\h{\phi}}^{(1)})^{-1}\alpha\right)=((\triangle_{\h{\omega},k\h{\phi}}^{(1)})^{-1}\alpha|\alpha)\leq \frac{1}{Ck}\|\alpha\|^2.
\end{align*}
\end{proof}
From above Theorem, we deduce the following ''Hodge decomposition'':
\begin{thm}\label{thm:alternative Hodge decomposition}
Let $\mathcal{P}_{\h{\omega},k\h{\phi}}:=\mathcal{P}^{(0)}_{\h{\omega},k\h{\phi}}:L^2(\C^n,\h{\omega})\to \ker\triangle_{\h{\omega},k\h{\phi}}^{(0)}=\ker(\db_{k\h{\phi}})$ be the orthogonal projection. Then it is given by
\begin{equation}\label{alternative Hodge decomposition}
\mathcal{P}_{\h{\omega},k\h{\phi}}=I-\db_{k\h{\phi}}^*(\triangle_{\h{\omega},k\h{\phi}}^{(1)})^{-1}\db_{k\h{\phi}}\quad \text{ on }C^\infty_c(\C^n).
\end{equation}
\end{thm}
\begin{proof}
First of all, for $u\in C^\infty_c(\C^n)$, we apply Theorem \ref{L^2-existence Theorem} to $v:=\db_{k\h{\phi}}u$ and thus $u_0=\db_{k\h{\phi}}^*(\triangle_{\h{\omega},k\h{\phi}})^{-1}\db_{k\h{\phi}}u$ solves $\db_{k\h{\phi}}u_0=\db_{k\h{\phi}}u$. Therefore, $u-u_0\in \ker\db_{k\h{\phi}}$. Also, for any $u\in \ker(\db_{k\h{\phi}})$, $(I-\db_{k\h{\phi}}^*(\triangle_{\h{\omega},k\h{\phi}}^{(1)})^{-1}\db_{k\h{\phi}})u=u=\mathcal{P}_{\h{\omega},k\h{\phi}}u$. This shows that 
\[
(I-\db_{k\h{\phi}}^*(\triangle_{\h{\omega},k\h{\phi}}^{(1)})^{-1}\db_{k\h{\phi}})^2u=u-\db_{k\h{\phi}}^*(\triangle_{\h{\omega},k\h{\phi}}^{(1)})^{-1}\db_{k\h{\phi}}u.
\] 
Finally, we compute
\begin{align*}
&(u-u_0|u_0)=(u|u_0)-(u_0|u_0)\\
=&(\db_{k\h{\phi}}^*(\triangle_{\h{\omega},k\h{\phi}}^{(1)})^{-1}\db_{k\h{\phi}}u|u)-(\db_{k\h{\phi}}^*(\triangle_{\h{\omega},k\h{\phi}}^{(1)})^{-1}\db_{k\h{\phi}}u|\db_{k\h{\phi}}^*(\triangle_{\h{\omega},k\h{\phi}}^{(1)})^{-1}\db_{k\h{\phi}}u)\\
=&((\triangle_{\h{\omega},k\h{\phi}}^{(1)})^{-1}\db_{k\h{\phi}}u|\db_{k\h{\phi}}u)-((\triangle_{\h{\omega},k\h{\phi}}^{1})^{-1}\db_{k\h{\phi}}u|\db_{k\h{\phi}}u)=0,
\end{align*}
where $\db_{k\h{\phi}}\db_{k\h{\phi}}^*(\triangle_{\h{\omega},k\h{\phi}}^{(1)})^{-1}\db_{k\h{\phi}}u=\db_{k\h{\phi}}u$ as in the proof of Theorem \ref{L^2-existence Theorem}. Therefore, $u-u_0\perp u_0$ and hence we conclude that $u-\db_{k\h{\phi}}^*(\triangle_{\h{\omega},k\h{\phi}}^{(1)})^{-1}\db_{k\h{\phi}}u=\mathcal{P}_{\h{\omega},k\h{\phi}}u$.
\end{proof}
Since $\mathcal{P}_{\h{\omega},k\h{\phi}}$ is a bounded operator on $L^2(\C^n,\h{\omega})$, we actually know that the \eqref{alternative Hodge decomposition} holds on $L^2(\C^n,\h{\omega})$ by density.
\begin{rmk}[Changing Metrics and Weights]\label{Changing Metrics and Weights} The results of this sections still hold if we replace $\widehat{\phi}$ and $\h{\omega}$ by even simpler standard metric $\phi_0=\sum_{j=1}^n \lambda_{j,x}|z^j|^2$ and $\omega_0:=\frac{\sqrt{-1}}{2}\sum_{j=1}^n dz^j\wedge d\bar{z}^j$.
\end{rmk}
\section{Symbolic Calculus and Asymptotic Sum}\label{symbolic calculus and Asymptotic Sum} To establish the asymptotic expansion for Bergman kernel, we develop the related symbol space and its  asymptotic sum in this section. We first define a space of functions which is rapidly decreasing off-diagonal (cf. \cite[section 3.1]{HsiaoNikhill}).
\begin{defin}\label{defin:Generalized Schwartz Space} The space $\widehat{S}(\R^d\times \R^d)$ consists of functions $a(x,y)\in C^\infty(\R^d\times \R^d)$ satisfying for any $(\alpha,\beta)\in \N_0^{2d}$. there exists $l=l(\alpha,\beta,a)\in \N$ such that for any $N>0$, there exists a constant $C=C_{\alpha,\beta,N}(a)>0$,
\begin{equation}\label{Generalized Schwartz Space}
\left|\p_x^\alpha\p_y^\beta a(x,y)\right|\leq C\frac{(1+|x|+|y|)^{l(\alpha,\beta)}}{(1+|x-y|)^N}, \quad \forall (x,y)\in  \R^d\times \R^d.
\end{equation}
\end{defin} 
Equivalently, \eqref{Generalized Schwartz Space} means that for any $\alpha,\beta\in \N_0^d$, there exists $l=l(\alpha,\beta)\in \N$ such that for any $N>0$,
\begin{equation}\label{Generalized Schwartz Space-2}
\sup_{x,y\in U} \frac{(1+|x-y|)^N\left|\p_x^\alpha\p_y^\beta a(x,y)\right|}{(1+|x|+|y|)^{l(\alpha,\beta)}}<\infty.
\end{equation}
Thus, if there exists $N_0>0$ such that \eqref{Generalized Schwartz Space-2} holds for $N>N_0$, then for $N\leq N_0$,
\[
\sup_{x,y\in U} \frac{(1+|x-y|)^N\left|\p_x^\alpha\p_y^\beta a(x,y)\right|}{(1+|x|+|y|)^{l(\alpha,\beta)}}\leq \sup_{x,y\in U} \frac{(1+|x-y|)^{N_0}\left|\p_x^\alpha\p_y^\beta a(x,y)\right|}{(1+|x|+|y|)^{l(\alpha,\beta)}}<\infty.
\] 
In other words, it suffices to show the condition \eqref{Generalized Schwartz Space} for sufficiently large $N$.
\begin{rmk}
One observes that if $a\in \widehat{S}(\R^d\times \R^d)$, then for each fixed $x,y\in U$, $a(x,\cdot)$, $a(\cdot,y)\in \mathscr{S}(\R^d)$, the Schwartz space of rapidly decreasing functions. If $a(x,y)\in \mathscr{S}(\R^d\times \R^d)$, then $a\in \widehat{S}(\R^d\times \R^d)$. However, the converse is not true. For instance, $e^{-|x-y|^2}\in \widehat{S}(\R^d\times \R^d)\setminus \mathscr{S}(\R^d\times \R^d)$. Hence, for any $\alpha,\beta\in \N_0^d$, any $a\in \widehat{S}(\R^d\times \R^d)$, for fixed $x,y\in U$, $\p_x^\alpha\p_y^\beta a(x,\cdot)$ and $\p_x^\alpha\p_y^\beta a(\cdot,y)$ are integrable in $x$ and $y$, repsectively.
\end{rmk}
Now, for a smooth function $a(x,y,k)$ with parameter $k$, recall that in Definition \ref{defin:semi-classical symbol space}, we have defined a kind a semi-classical symbol space. For $m\in \R$, a function $a(x,y,k)\in \widehat{S}^m(\R^d\times \R^d)$ if
\begin{itemize}
\item[(i)] $a(x,y,k)\in C^\infty(\R^d\times \R^d)$, for each $k\in \N$, and
\item[(ii)] for any$(\alpha,\beta)\in \N_0^{2d}$, there exists $l=l(\alpha,\beta,a)\in \N$ and $k_0\in \N$ such that for any $N>0$,  there exists a constant $C=C_{\alpha,\beta,N}(a)>0$,
\begin{equation}\label{semi-classical symbol space-2}
\left|\p_x^\alpha\p_y^\beta a(x,y,k)\right|\leq Ck^{m+\frac{|\alpha|+|\beta|}{2}}\frac{(1+|\sqrt{k}x|+|\sqrt{k}y|)^{l(\alpha,\beta)}}{(1+|\sqrt{k}(x-y)|)^N},\end{equation}
for any  $(x,y)\in \R^d\times \R^d$, any $k\geq k_0$. 
\end{itemize}
Similarly, one only needs to verify \eqref{semi-classical symbol space-2} for $N>N_0$, for some $N_0>0$, and \eqref{semi-classical symbol space-2} shows that for each fixed $x,y\in \R^d$, $k\in \N$, $a(x,\cdot,k), a(\cdot,y,k)\in \mathscr{S}(\R^d)$. 

The following Lemma asserts that such space is not vacuous.
\begin{lem}
If $a(x,y)\in \widehat{S}(\R^d\times \R^d)$, then $a(\sqrt{k}x,\sqrt{k}y)\in \widehat{S}^0(\R^d\times \R^d)$ and for any $m\in \R$, $k^ma(\sqrt{k}x,\sqrt{k}y)\in \widehat{S}^m(\R^d\times \R^d)$.
\end{lem}
\begin{proof}
The condition (i) in above definition is obviously satisfied by $a(\sqrt{k}x,\sqrt{k}y)$ and $k^ma(\sqrt{k}x,\sqrt{k}y)$, for any $a\in \widehat{S}(\R^d\times \R^d)$ and any $m\in \R$. For (ii),
for any $\alpha,\beta\in \N_0^d$, by \eqref{Generalized Schwartz Space}, we have the following estimates for $a(\sqrt{k}x,\sqrt{k},y)$:
\begin{align*}
&\left|\p_x^\alpha\p_y^\beta k^ma(\sqrt{k},\sqrt{k}y)\right|=k^{m+\frac{|\alpha|+|\beta|}{2}}|(\p_x^\alpha\p_x^\beta a)(\sqrt{k}x,\sqrt{k}y)|\\
\leq& Ck^{m+\frac{|\alpha|+|\beta|}{2}}\frac{(1+|\sqrt{k}x|+|\sqrt{k}y|)^{l(\alpha,\beta)}}{(1+|\sqrt{k}(x-y)|)^N}.
\end{align*}
\end{proof}
Now, we investigate some easy properties of symbol spaces.
\begin{lem}\label{basic properties for symbol spaces} For $m,m'\in \R$, $a\in \widehat{S}^m(\R^d\times \R^d)$ and $b\in \widehat{S}^{m'}(\R^d\times \R^d)$, then
\begin{itemize}
\item[(i)] $\widehat{S}^m(\R^d\times \R^d)\subset \widehat{S}^{m'}(\R^d\times \R^d)$ if $m<m'$.
\item[(ii)] For any $\alpha,\beta\in \N_0^d$, $\p_x^\alpha\p_y^\beta a\in \widehat{S}^{m+\frac{\alpha+\beta}{2}}(\R^d\times \R^d)$.
\item[(iii)] $ab\in \widehat{S}^{m+m'}(\R^d\times \R^d)$.
\end{itemize}
\end{lem}
\begin{proof}
(i) and (ii) are obvious from \eqref{semi-classical symbol space-2}. For (iii), for any $\alpha,\beta\in \N_0^d$, by assumption, there exists $l'(\alpha,\beta),l''(\alpha,\beta)\in \N$ such that for any $N\in \N$, there exists $C',C''>0$ such that
\begin{align*}
&\left|\p_x^\alpha\p_y^\beta a(x,y,k)\right|\leq C'k^{m+\frac{|\alpha|+|\beta|}{2}}\frac{(1+|\sqrt{k}x|+|\sqrt{k}y|)^{l'(\alpha,\beta)}}{(1+|\sqrt{k}(x-y)|)^{\frac{N}{2}}},\\
&\left|\p_x^\alpha\p_y^\beta b(x,y,k)\right|\leq C''k^{m'+\frac{|\alpha|+|\beta|}{2}}\frac{(1+|\sqrt{k}x|+|\sqrt{k}y|)^{l''(\alpha,\beta)}}{(1+|\sqrt{k}(x-y)|)^\frac{N}{2}}.
\end{align*}
Then for any $(x,y)\in \R^d\times \R^d$, any sufficiently large $k\in \N$, any $N\in \N$,
\begin{align*}
&\left|\p_x^\alpha\p_y^\beta (ab)(x,y,k)\right|\\
\leq &CC'k^{m+m'+\frac{|\alpha|+|\beta|}{2}}\sum_{\alpha'+\alpha''=\alpha,\beta'+\beta''=\beta} \frac{(1+|\sqrt{k}x|+|\sqrt{k}y|)^{l'(\alpha',\beta')+l''(\alpha'',\beta'')}}{(1+|\sqrt{k}(x-y)|)^N}
\end{align*}
One then takes $l(\alpha,\beta)=\max_{\alpha'+\alpha''=\alpha,\beta'+\beta''=\beta}\{l'(\alpha',\beta')+l''(\alpha'',\beta'')$.
\end{proof}
From the above Lemma (i), we then define \textbf{the space of symbols which rapidly decreasing in $k$} by $\widehat{S}^{-\infty}(\R^d\times \R^d):=\bigcap_{m\in \R}\widehat{S}^m(\R^d\times \R^d)$. To define the of asymptotic expansion in our case, we need to establish the notion of asymptotic sum.
\begin{thm}\label{thm:asymptotic sum for semi-classical symbol space} Given any sequence $\{m_j\}_{j=0}^\infty$ with  $m_j\searrow -\infty$ and $a_j(x,y,k)\in \widehat{S}^{m_j}(\R^d\times \R^d)$, there exists $a(x,y,k)\in \widehat{S}^{m_0}(\R^d\times \R^d)$ such that for any $q\in \N_0$, there exists $k_0(q)\in \N$ such that if $k\geq k_0$,
\begin{equation}\label{asymptotic sum for semi-classical symbol space}
a(x,y,k)-\sum_{j=0}^{q}a_j(x,y,k)\in \widehat{S}^{m_{q+1}}(\R^d\times \R^d).
\end{equation}

Moreover, such $a$ is uniquely modulo $\widehat{S}^{-\infty}(\R^d\times \R^d)$.
\end{thm}
\begin{proof}
For any positive sequence $\{\mu_j\}_{j=0}^\infty$ with $\lambda_j\nearrow \infty$, we define
\[
\tau_{j,k}=\mathbf{1}_{[0,1]}(\mu_j/k).
\]
For any $\alpha,\beta\in \N_0^d$ with $|\alpha|+|\beta|\leq j$, there exists $l_j\in \N$ such that for any $N>0$, there exists a constant $C=C_{N,j}>0$ satisfying
\[
|\p_x^\alpha\p_y^\beta a_j(x,y,k)|\leq C_{N,j}k^{m_j+\frac{j}{2}}\frac{1+|\sqrt{k}|x|+\sqrt{k}|y|)^{l_j}}{(1+|\sqrt{k}(x-y)|)^N},\quad \forall x,y\in \R^d.
\]
Let $\epsilon_j$ be a positive sequence tending to $0$ to be chosen later and
\begin{equation}\label{Borel construction}
a(x,y,k):=\sum_{j=0}^\infty A_j(x,y,k),\quad A_j(x,y,k):=a_j(x,y,k)\chi(k^{\frac{1}{2}-\epsilon_j}x,k^{\frac{1}{2}-\epsilon_j}y)\tau_{j,k},
\end{equation}
where $\chi\in C^\infty_c(\R^d\times \R^d)$ such that $0\leq \chi\leq 1$, $\chi(x,y)=1$ for $|x|,|y|\leq 1$, and $\chi(x,y)=0$ for $|x|,|y|\geq 2$

First of all,  $\tau_{j,k}\neq 0$ if and only if $\mu_j<k$. Since $\mu_j\nearrow \infty$, for each fixed $k$, there exists only finitely many $j$ with $\mu_j<k$. Hence, for each fixed $k\in \N$, \eqref{Borel construction} is a finite sum and hence $a(x,y,k)\in C^\infty(\R^d\times \R^d)$. Now, for sufficiently large $k$ so that $\tau_{j,k}=1$, we have
\[
A_j(x,y,k)-a_j(x,y,k)=\sum_{j=0}^q (\chi(k^{\frac{1}{2}-\epsilon_j}x,k^{\frac{1}{2}-\epsilon_j}y)-1)a_j(x,y,k).
\]
We claim that 
\begin{claim} For any $a\in \widehat{S}^m(\R^d\times \R^d)$ and $\epsilon>0$, $a(x,y,k)(\chi(k^{\frac{1}{2}-\epsilon}x,k^{\frac{1}{2}-\epsilon}y)-1)\in \widehat{S}^{-\infty}(\R^d\times \R^d)$, i.e., given any $r\in \N$, $a(\chi(k^{\frac{1}{2}-\epsilon}x,k^{\frac{1}{2}-\epsilon}y)-1)\in \widehat{S}^{m-r}(\R^d\times \R^d)$. 
\end{claim} 
\begin{proof}[Proof of Claim] The key is the following. For any $M\in \N$, any $\alpha,\beta\in \N_0^d$, there exists $C=C_{M,\alpha,\beta}>0$ such that
\begin{equation}\label{cut-off}
|\p_x^\alpha\p_y^\beta \chi(x,y)-1|\leq C_{M,\alpha,\beta}(|x|^M+|y|^M),\quad \forall x,y\in \C^n.
\end{equation}
For the proof of  \eqref{cut-off}, for $|x|,|y|<1$, $1-\chi(x,y)=0$ and for $|x|,|y|>2$, $1-\chi(x,y)=1$, the estimate holds obviously. For $1\leq |x|,|y|\leq 2$, we can expand
\[
\chi(x,y)-1=\sum_{0<|\alpha|+|\beta|\leq M} \frac{\p^{|\alpha|+|\beta|}\chi}{\p x^\alpha \p y^\beta}(0)x^\alpha y^\beta+\sum_{|\alpha|+|\beta|=M}R_{\alpha\beta}(x,y)x^\alpha y^\beta,
\]
where $R_{\alpha\beta}\in C^\infty(\R^d\times \R^d)$. Therefore, $|\p_x^\alpha\p_y^\beta \chi(x,y)-1|\leq C_{\alpha,\beta,M}(|x|^M+|y|^M)$.
Now, for any $\alpha,\beta\in \N_0^d$, any $N>0$,
\begin{align*}
&|\p_x^\alpha\p_y^\beta a(1-\chi(k^{1/2-\epsilon}x,k^{1/2-\epsilon}y))|\\
\leq &\sum_{\alpha'+\alpha''= \alpha,\beta'+\beta''=\beta}C'_{\alpha,\beta}k^{m+\frac{|\alpha'|+|\beta'|}{2}}\frac{(1+\sqrt{k}|x|+\sqrt{k}|y|)^{l(\alpha',\beta')}}{(1+\sqrt{k}|x-y|)^N}\\
\times &C''_{\alpha'',\beta'',M}(|k^{\frac{1}{2}-\epsilon}x|^{M+|\alpha''|+|\beta''}+|k^{\frac{1}{2}-\epsilon}y|^{M+|\alpha''|+|\beta''|})\\
\leq &C_3k^{m-\epsilon M+|\alpha|+|\beta|}\frac{(1+\sqrt{k}|x|+\sqrt{k}|y|)^{l_1(\alpha,\beta)+M}}{(1+\sqrt{k}|x-y|)^N},
\end{align*}
where $l_1(\alpha,\beta)=\max_{\alpha'+\alpha''=\alpha,\beta'+\beta''=\beta}\left(l(\alpha',\beta')+|\alpha''|+|\beta''|+M\right)$. Now, we choose $M$ so that $\epsilon M>r$ and thus 
\[
a(1-\chi(k^{\frac{1}{2}-\epsilon}x,k^{\frac{1}{2}-\epsilon}y))\in \widehat{S}^{m-\epsilon M}(\R^d\times \R^d)\subset \widehat{S}^{m-r}(\R^d\times \R^d).
\]
\end{proof}
As a result, $A_j(x,y,k)-a_j(x,y,k)\in \widehat{S}^{-\infty}(\R^d\times \R^d)$ for sufficiently large $k$. On the other hand, for any $j\in \N$, any $\alpha,\beta\in \N_0^d$ with $|\alpha|+|\beta|\leq j$, there exists $l_j\in \N$ such that 
\[
|\p_x^\alpha\p_y^\beta a_j|\leq C_jk^{m_j-\frac{j}{2}}\frac{(1+\sqrt{k}|x|+\sqrt{k}|y|)^{l_j}}{(1+\sqrt{k}|x-y|)^j}.
\]
Now, we apply above estimates to $A_j(x,y,k)=a_j(x,y,k)\chi(k^{\frac{1}{2}-\epsilon_j}x,k^{\frac{1}{2}-\epsilon_j}y)\tau_{j,k}$ and
\begin{align*}
&|\p_x^\alpha\p_y^\beta A_j(x,y,k)|\\
\leq &C_jk^{m_j-\frac{j}{2}-\frac{\epsilon_jj}{2}}\frac{(1+\sqrt{k}|x|+\sqrt{k}|y|)^{l_j}}{(1+\sqrt{k}|x-y|)^j}\sup_{|\alpha|+|\beta|\leq j}|(\p_x^\alpha\p_y^\beta\chi)(k^{\frac{1}{2}-\epsilon_j}x,k^{\frac{1}{2}-\epsilon_j}y)|.
\end{align*}
Since $\chi$ supports in $|x|,|y|\leq 2$, we see that $\sqrt{k}|x|,\sqrt{k}|y|\leq 2k^{\epsilon_j}$, we see that
\[
|\p_x^\alpha\p_y^\beta A_j(x,y,k)|\leq C'_jk^{m_j-\frac{j}{2}+\epsilon_jl_j}(1+\sqrt{k}|x-y|)^{-j}.
\]
We take $\epsilon_j$ so that $k^{\epsilon_jl_j-1}\leq \frac{1}{2^jC'_j}$ for sufficiently large $k$ and hence
\[
C'_jk^{m_j-\frac{j}{2}+\epsilon_jl_j}\leq 2^{-j}k^{m_j-\frac{j}{2}+1}.
\]
Given $\alpha,\beta\in \N_0^d$, $q\in \N$, we take $N\geq \max\{|\alpha|+|\beta|,q+1\}$ and $m_N+1\leq m_{q+1}$.
\begin{align*}
&\left|\p_x^\alpha\p_y^\beta (\sum_{j=N}^\infty A_j)\right|\leq\sum_{j=N}^\infty \frac{C'_jk^{m_j-\frac{j}{2}+\epsilon_jl_j}}{(1+\sqrt{k}|x-y|)^j}\\\leq &\sum_{j=N}^\infty \frac{k^{m_{q+1}-\frac{N}{2}}}{2^j(1+\sqrt{k}|x-y|)^N}\leq \frac{k^{m_{q+1}-\frac{|\alpha|+|\beta|}{2}}}{(1+\sqrt{k}|x-y|)^N}.
\end{align*}
Hence, 
\[
\left|\p_x^\alpha\p_y^\beta (a-\sum_{j=0}^qa_j)\right|\leq \left|\p_x^\alpha\p_y^\beta (\sum_{j=N}^{\infty} A_j)\right|+\left|\p_x^\alpha\p_y^\beta \sum_{j=0}^q (a_j-A_j)\right|+\left|\p_x^\alpha\p_y^\beta\sum_{j=q+1}^{N-1} A_j\right|
\]
Since $\chi\in \h{S}(\R^d\times \R^d)$, $A_j\in \h{S}^{m_{j}}(\R^d\times \R^d)$ (cf. Lemma \ref{lem:multiplication by k-cut-off}). Also, previous argument and above claim shows that $\sum_{j=N}^\infty A_j\in \h{S}^{m_{q+1}}(\R^d\times \R^d)$ and $\sum_{j=0}^q (a_j-A_j)\in \h{S}^{-\infty}(\R^d\times \R^d)$, we see that
\[
|\p_x^\alpha\p_y^\beta (a-\sum_{j=0}^qa_j)|\leq C_{\alpha,\beta,N}k^{m_{q+1}-\frac{|\alpha|+|\beta|}{2}}\frac{(1+|\sqrt{k}x|+|\sqrt{k}y|)^{l(\alpha,\beta)}}{(1+\sqrt{k}|x-y|)^N}.
\] 
In other words, $a-\sum_{j=0}^qa_j\in \h{S}^{m_{q+1}}$.
\end{proof}
If $a$ and $\{a_j\}$ satisfy the conclusion of the Theorem \ref{thm:asymptotic sum for semi-classical symbol space}, we then write $a\sim \sum_{j=0}^\infty a_j(x,y,k)$ and call $a$ the \textbf{asymptotic sum} for $\{a_j\}_{j=1}^\infty$. Moreover, we define:
\begin{defin}\label{defin:semi-classical classical symbol} The space $\widehat{S}^m_{\cl}(\R^d\times \R^d)$ of \textbf{classical symbol} of order $m$ consists of function $a(x,y,k)\in \widehat{S}^m(\R^d\times \R^d)$ such that there exists a sequence $a_j\in \widehat{S}(\R^d\times \R^d)$ for $j\in \N_0$ satisfying 
\begin{equation}\label{expansion of classical symbol}
a(x,y,k)\sim \sum_{j=0}^\infty k^{m-\frac{j}{2}}a_j(\sqrt{k}x,\sqrt{k}y)
\end{equation}
\end{defin}
Next, we define the \textbf{quantization} on symbol space $\widehat{S}^m(\R^d\times \R^d)$.

\begin{defin}\label{defin:quantization of semi-classical symbol} Given $a\in \widehat{S}^{m}(\R^d\times \R^d)$, we define a $k$-dependent continuous linear operator $Op_k(a)\in \widehat{L}^m(\R^d)$ by
\begin{equation}\label{quantization of semi-classical symbol}
Op_k(a)(u)(x)=\int_{\R^d}a(x,y,k)u(y)dm(y).
\end{equation}
A $k$-dependent continuous linear operator $A_k:C^\infty_c(\R
^n)\to D'(\R^d)$ is in the class $\widehat{L}^m(\R^d)$ if $A=Op_k(a)$ for some $a(x,y,k)\in \widehat{S}^m(\R^d\times \R^d)$.
\end{defin}
In particular, $A_k\in \widehat{L}^m(\R^d)$ implies the Schwartz kernel $K_{A_k}(x,y)=a(x,y,k)\in C^\infty(\R^d\times \R^d)$. Thus, $A_k\in \widehat{L}^m(\R^d)$ is a smoothing operator for any $k\in \N$. Moreover, if $a\in \widehat{S}^{-\infty}(\R^d\times \R^d)$, then $Op_k(a)$ is $k$-negligible in the sense of Definition \ref{defin:k-negligible} obviously. Hence, we may extend the Definition \ref{defin:quantization of semi-classical symbol} by $A_k\in \widehat{L}^m(\R^d)$ if there exists $a\in \widehat{S}^m(\R^d\times \R^d)$ such that $A_k-Op_k(a)= Op_k(a_1)$ with $a_1\in \widehat{S}^{-\infty}(\R^d\times \R^d)$. 

We also define the subclass $\widehat{L}^m_{\cl}(\R^d)\subset \widehat{L}^m(\R^d)$ by $A_k\in \widehat{L}^m_{\cl}(\R^d)$ if $A_k-Op_k(a)=Op_k(a_1)$, for some $a\in \widehat{S}^m_{\cl}(\R^d\times \R^d)$ and $a_1\in \widehat{S}^{-\infty}(\R^d\times \R^d)$. For $A_k=Op_k(a)\in \widehat{L}^m_{\cl}(\R^d)$, we then define the \textbf{principal symbol} $\sigma(A_k)$ by the leading term $a_0(x,y)\in \widehat{S}(\R^d\times \R^d)$ in the asymptotic sum \eqref{expansion of classical symbol}.
\begin{thm}\label{thm:composition and adjoint of operators} If $A_k=Op_k(a)\in \widehat{L}^m(\R^d)$, $B_k=Op_k(b)\in \widehat{L}^{m'}(\R^d)$, then
\begin{itemize}
\item[(i)] The formal adjoint $A^*_k\in \widehat{L}^m(\R^d)$ with $A^*_k=Op_k(a^*)$, where $a^*(x,y,k):=\overline{a(y,x,k)}$. 
\item[(ii)] $A_k\circ B_k\in \widehat{L}^{m+m'-\frac{d}{2}}(\R^d\times \R^d)$ with $A_k\circ B_k=Op_k(a\#b)$, where
\begin{equation}\label{composition of symbols}
(a\#b)(x,y,k):=\int_{\R^d}a(x,t,k)b(t,y,k)dm(t)\in \widehat{S}^{m+m'-\frac{d}{2}}(\R^d\times \R^d),
\end{equation}
i.e., if $a\in \widehat{S}^m_{\cl}(\R^d\times \R^d)$, $b\in \widehat{S}^{m'}_{\cl}(\R^d\times \R^d)$, then $a\#b\in \widehat{S}^{m+m'-\frac{d}{2}}(\R^d\times \R^d)$.
\end{itemize}
\end{thm}
\begin{proof}
For (i), for $u,v\in C^\infty_c(\R^d)$, the formal adjoint $A^*_k$ is given by $(A_k(u)|v)=(u|A_kv)$. Hence, we compute $A_k^*$ explicitly as
\begin{align*}
(A_ku|v)&=\int_{\R^d}\left(\int_{\R^d}a(x,y,k)u(y)dm(y)\right)\overline{v(x)}dm(x)\\
&=\int_{\R^d}\int_{\R^d} a(x,y,k)u(y)\overline{v(x)}dm(y)dm(x)\\
&=\int_{\R^d}u(y)\overline{\int_{\R^d}\overline{a(x,y,k)}v(x)dm(x)}dm(y)\\
&=\int_{\R^d}u(y)\overline{A_k^*(v)(y)}dy.
\end{align*}
This implies that $(A_k^*v)(x)=\int_{\R^d} \overline{a(y,x,k)}v(y)dm(y)$ and thus $A_k^*=Op_{a^*}$, where $a^*(x,y,k):=\overline{a(y,x,k)}$. Obviously, $a^*\in \widehat{S}^m_{\cl}(\R^d\times \R^d)$ if $a$ is.

For (ii), for each $f\in C^\infty_c(\R^d)$, $(B_kf)(t)=\int_{\R^d}b(t,y,k)f(y)dy$ and thus
\[
(A_k\circ B_k)(f)(x)=\int_{\R^d}a(x,t,k)\left(\int_{\R^d}b(t,y,k)f(y)dm(y)\right)dm(t)
\]
Observe that for fixed $x\in \R^d$, $k\in \N$, $a(x,t,k)b(t,y,k)f(y)$ is integrable in $t$ and $y$. Therefore, by Fubini-Tonelli Theorem, we then have
\begin{align*}
(A_k\circ B_k)(f)(x)&=\int_{\R^d}\left( \int_{\R^d}a(x,t,k)b(t,y,k)dm(t)\right)f(y)dm(y)\\
&=\int_{\R^d} (a\#b)(x,y,k)f(y)dm(y).
\end{align*}
We then see that $A_k\circ B_k=Op_{a\#b}$. Now, we show that $a\#b\in \widehat{S}^{m+m'-\frac{d}{2}}(\R^d\times \R^d)$. To see this, for any $\alpha,\beta\in \N_0^d$,
\begin{align*}
&\left|\p_x^\alpha\p_y^\beta (a\#b)\right|\leq \int_{U} \left|\p_x^\alpha a(x,t,k)\p_y^\beta b(t,y,k)\right|dm(t)\\
\leq &C_{\alpha,\beta,N}k^{m+m'+\frac{|\alpha|+|\beta|}{2}} \int_{\R^d} \frac{(1+\sqrt{k}|x|+\sqrt{k}|t|)^{l(\alpha,0)}}{(1+|\sqrt{k}(x-t)|)^N}\frac{(1+|\sqrt{k}t|+|\sqrt{k}y|)^{l'(0,\beta)}}{(1+|\sqrt{k}(t-y)|)^N}dm(t)
\end{align*}
We make the change of variable $s=\sqrt{k}t$ and thus $dm(s)=k^{\frac{d}{2}}dm(t)$. This implies
\begin{align*}
&\left|\p_x^\alpha\p_y^\beta (a\#b)\right|\\
\leq &Ck^{m+m'-\frac{d}{2}+\frac{|\alpha|+|\beta|}{2}}\int_{\R^d} \frac{(1+\sqrt{k}|x|+|s|)^{l(\alpha,0)}}{(1+|\sqrt{k}x-s|)^N}\frac{(1+|s|+|\sqrt{k}y|)^{l'(0,\beta)}}{(1+|s-\sqrt{k}y)|)^N}dm(s).
\end{align*}
Let $l(\alpha,\beta)=\max\{l(\alpha,0),l'(0,\beta)\}$. We observe that for any $M>0$,
\begin{align*}
&(1+|\sqrt{k}x-s|)^M(1+|\sqrt{k}y-s|)^M\geq(1+|\sqrt{k}(x-y)|+|\sqrt{k}x-s||\sqrt{k}y-s|)^M\\
=&(1+|\sqrt{k}(x-y)|)^M\left(1+\frac{|\sqrt{k}x-s||\sqrt{k}y-s|}{1+\sqrt{k}|x-y|}\right)^M\geq (1+\sqrt{k}|x-y|)^M
\end{align*}
Hence, by taking $M=N/2$ and $u=s-\sqrt{k}x$, we can write
\begin{align*}
&\left|\p_x^\alpha\p_y^\beta (a\#b)\right|\\
\leq &C\frac{k^{m+m'-\frac{d}{2}+\frac{|\alpha|+|\beta|}{2}}}{(1+\sqrt{k}|x-y|)^{N/2}}\int_{\R^d}\frac{(1+\sqrt{k}|x|+|s|)^{l(\alpha,\beta)}(1+\sqrt{k}|y|+|s|)^{l(\alpha,\beta)}}{(1+|\sqrt{k}x-s|)^{N/2}(1+|\sqrt{k}y-s|)^{N/2}}ds\\
\leq &C\frac{k^{m+m'-\frac{d}{2}+\frac{|\alpha|+|\beta|}{2}}}{(1+\sqrt{k}|x-y|)^{N/2}}\int_{\R^d} \frac{(1+2\sqrt{k}|x|+|u|)^{l(\alpha,\beta)}(1+\sqrt{k}|x|+\sqrt{k}|y|+|u|)^{l(\alpha,\beta)}}{(1+|u|)^{N/2}}du.
\end{align*}
Hence, there exists $N_0=N_0(\alpha,\beta,d)$ such that if $N>N_0$, the integral converges. By expanding the numerator of the integrand, we can find $l'(\alpha,\beta)$ so that for any $N>N_0$,   
\[
\left|\p_x^\alpha\p_y^\beta (a\#b)\right|\leq C_{\alpha,\beta,N}k^{m+m'-\frac{d}{2}+\frac{|\alpha|+|\beta|}{2}}\frac{(1+\sqrt{k}|x|+\sqrt{k}|y|)^{l'(\alpha,\beta)}}{(1+\sqrt{k}|x-y|)^{N/2}}.
\]

Thus, $a\#b\in \h{S}^{m+m'-\frac{d}{2}}(\R^d\times \R^d)$.
\end{proof}
\section{Asymptotic Expansion of Approximate Kernel}\label{Asymptotic Expansion of Approximate Kernel}
In this section, we establish the asymptotic expansion of approximate Bergman kernel $P_{\h{\omega},k\h{\phi}}^{(0)}$ through the symbolic calculus presented in section \ref{symbolic calculus and Asymptotic Sum}. 

Following remark \ref{Changing Metrics and Weights}, we first consider the case $P^{(q)}_{\omega_0,k\phi_0}$, where $\phi_0=\sum_{j=1}^n\lambda_{j,x}|z^j|^2$, $\omega_0=\frac{\sqrt{-1}}{2}\sum_{j=1}^n dz^j\wedge d\bar{z}^j$, which is the orthogonal projection $L^2_{0,q}(\C^n):=L^2_{0,q}(\C^n,\omega_0)\to \ker\triangle_{\omega_0,k\phi_0}$, where the analogous Laplacian $\triangle_{\omega_0,k\phi_0}^{(q)}$ is given by
\[
\triangle_{\omega_0,k\phi_0}^{(q)}:=\db_{k\phi_0}^{*,\omega_0}\db_{k\phi_0}+\db_{k\phi_0}\db_{k\phi_0}^{*,\omega_0},
\]
$\db_{k\phi_0}:=\db+k\epsilon(\db\phi_0)$, and $\db_{k\phi_0}^{*,\omega_0}$ is the formal adjoint with respect to the $L^2$-inner product
\[
(\alpha|\beta)_0:=\sum_{|I|=J}\nolimits^{\prime}\int_{\C^n}\alpha_I\overline{\beta_I}dm(z).
\]
Let $\delta_k(z)=\frac{z}{\sqrt{k}}$ be the scaling map on $\C^n$ with inverse $\delta_k^{-1}(z)=\sqrt{k}z$. Then for $u\in C^\infty(\C^n)$,
\begin{align*}
&\delta_k\db_{k\phi_0}\delta_k^{-1}u(z)=\delta_k\left(\sum_{i=1}^n\left( \sqrt{k}\frac{\p u}{\p \bar{z}^i}(\sqrt{k}z)+k\frac{\p \phi_0}{\p \bar{z}^i}(z)u(\sqrt{k}z)\right)d\bar{z}^i\right)\\
=&\delta_k\left(\sum_{i=1}^n \left(\sqrt{k}\frac{\p u}{\p \bar{z}^i}(\sqrt{k}z)+k\lambda_{i,x}z^iu(\sqrt{k}z)\right)d\bar{z}^i\right)\\
=&\delta_k\left(\sqrt{k}\sum_{i=1}^n\left( \frac{\p u}{\p \bar{z}^i}(\sqrt{k}z)+\sqrt{k}\lambda_{i,x}z^iu(\sqrt{k}z)\right)d\bar{z}^i\right)\\
=&\sqrt{k}\delta_k\delta_k^{-1}(\db_{\phi_0}u)(z)=\sqrt{k}\db_{\phi_0}u(z).
\end{align*}
Therefore, we can deduce that
\begin{equation}\label{scaling on deformed dbar and dbar*}
\db_{\phi_0}\delta_k=\frac{1}{\sqrt{k}}\delta_k\db_{k\phi_0},\quad \db_{\phi_0}^{*,\omega_0}\delta_k=\frac{1}{\sqrt{k}}\delta_k\db_{k\phi_0}^{*,\omega_0}.
\end{equation}
Hence, from \eqref{scaling on deformed dbar and dbar*}  we get
\begin{equation}\label{scaling on deformed Laplace}
\triangle_{\omega_0,\phi_0}\delta_k=\frac{1}{k}\delta_k\triangle_{\omega_0,k\phi_0}.
\end{equation}
Using \eqref{scaling on deformed Laplace}, if $\{\sigma_j(z)\}_{j=1}^d$, where $d\in \N_0\cup\{\infty\}$, is an orthonormal basis of $\ker \triangle_{\omega_0,k\phi_0}$ with respect to $(\cdot|\cdot)_{\omega_0}$, then $\delta_k\sigma_j$ satisfies $\triangle_{\omega_0,\phi_0}\delta_k\sigma_j=\frac{1}{k}\delta_k\triangle_{\omega_0,k\phi_0}\sigma_j=0$. Their inner product is given by
\[
(\delta_k{\sigma}_i|\delta_k{\sigma}_j)_{\omega_0}=\int_{\C^n}\delta_k{\sigma}_i(w)\overline{\delta_k{\sigma}_j(w)}dm(w)=k^n \int_{\C^n} \sigma_i(w)\overline{\sigma_j(w)}dm(w)=k^n\delta_{ij}.
\]
This shows that $\{k^{-\frac{n}{2}}\delta_k\sigma_j\}_{j=1}^d$ is an orthonormal basis for $\ker\triangle_{\omega_0,\phi_0}$. By \eqref{BK as expansion into ONB},
\begin{equation}
P_{\omega_0,\phi_0}(z,w)=k^{-n}\sum_{j=1}^d \sigma_j\left(\frac{z}{\sqrt{k}}\right)\overline{\sigma_j\left(\frac{z}{\sqrt{k}}\right)}=k^{-n}P_{\omega_0,k\phi_0}\left(\frac{z}{\sqrt{k}},\frac{w}{\sqrt{k}}\right).
\end{equation}
For $P_{\omega_0,\phi_0}$, we can compute it explicitly.
\begin{prop}\label{Schwartz kernel for simple model case} The approximate Bergman kernel $P_{\omega_0,\phi_0}^{(0)}(z,w)$ is given by
\[
P_{\omega_0,\phi_0}^{(0)}(z,w)=\frac{2^n\lambda_{1,x}\cdots \lambda_{n,x}}{\pi^n}e^{\sum_{j=1}^n\lambda_{j,x}\left(2z^j\overline{w}^j-|z^j|^2-|w^j|^2\right)}.
\]
\end{prop}
\begin{proof} First, we consider the trivial line bundle $L=\C\times \C^n$ over $\C^n$ with weight $|1|_{h^L}^2=e^{-2\phi_0}$. Its $L^2$-section can be identified as the weighted $L^2$-space $L^2(\C^n, e^{-2\phi_0}dm)$, and the subspace of holomorphic sections is identified the subspace $\mathcal{F}$, known as \textbf{Bargmann--Fock space}, consists of entire functions $f$ satisfying
\[
\|f\|_{\omega_0,\phi_0}^2=\int_{\C^n}|f(z)|^2e^{-2\sum_{j=1}^n\lambda_{j,x}|z^j|^2}dm(z)<\infty,\quad \db f=0,\quad \forall f\in \mathcal{F}.
\]
Then in section \ref{Bergman Kernel}, we denote $K_{\mathrm{BF}}(z,w)$ by the Schwartz kernel of the orthogonal projection $\Pi_{\mathrm{BF}}:L^2(\C^n,L)\to H^0(X,L)$, i.e.,
\[
(\Pi_{\mathrm{BF}}f)(z)=f(z)=\int_{\C^n} K_{\mathrm{BF}}(z,w)(f(w))dm(w),\quad \forall f\in H^0(X,L).
\]
\begin{claim}
$K_{\mathrm{BF}}(z,w)=\frac{2^n\lambda_{1,x}\cdots \lambda_{n,x}}{\pi^n}e^{2\sum_{j=1}^n\lambda_{j,x}(z^j\overline{w}^j-|w^j|^2)}$.
\end{claim}
\begin{proof}
For multi-index $\alpha=(\alpha_1,\dots,\alpha_n)\in \N_0^n$, we let $z^\alpha:=(z^1)^{\alpha_1}\dots (z^n)^{\alpha_n}$. Clearly, $\db(z^\alpha)=0$, for any $\alpha\in \N_0^n$. Hence, $z^\alpha\in \OO(\C^n)$. For $\alpha,\beta\in \N_0^n$, using polar coordinate $z^j=r_je^{-\ii\theta_j}$ and Fubini--Torelli Theorem,
\begin{align*}
&(z^\alpha|z^\beta)_{\omega_0,\phi_0}=\int_{\C^n} z^\alpha\bar{z}^\beta e^{-2\sum_{j=1}^n\lambda_{x,j}|z^j|^2}dm\\
=&\prod_{j=1}^n \left[\int_0^\infty\int_0^{2\pi} r_j^{\alpha_j+\beta_j+1}e^{\ii(\alpha_j-\beta_j)\theta_j}e^{-2\lambda_{x,j}r_j^2}d\theta_jdr_j\right].
\end{align*}
If $\alpha_j\neq \beta_j$, then $\int_0^{2\pi} e^{\ii(\alpha_j-\beta_j)\theta_j}d\theta_j=0$. Hence, $(z^\alpha,z^\beta)=0$ if $\alpha\neq \beta$. Now, observe that for $l\in \N$, 
\begin{equation}\label{freshman integral}
\int_0^\infty r^{2l+1}e^{-2\lambda r^2}dr=\frac{1}{2(2\lambda)^{l+1}}\int_0^\infty u^le^{-u}du=\frac{\Gamma(l+1)}{2(2\lambda)^{l+1}}=\frac{l!}{2(2\lambda)^{l+1}},
\end{equation} 
where $u=2\lambda r^2$. Therefore, the square of norm of $z^\alpha$ is given by
\[
\|z^\alpha\|^2_{\omega_0,\phi_0}=\prod_{j=1}^n \frac{2\pi \alpha_j!}{2(2\lambda_{j,x})^{\alpha_j+1}}=\frac{\pi^n\alpha!}{2^{|\alpha|+n}\lambda^{\alpha+1}},
\]
where $\lambda^{\alpha+1}:=\prod_{j=1}^n \lambda_{j,x}^{\alpha_j+1}$. As a result, $\left\{\Psi_\alpha:=\sqrt{\frac{2^{|\alpha|+n}\lambda^{\alpha+1}}{\pi^n\alpha!}}z^{\alpha}\right\}_{\alpha\in \N_0^n}$ is an orthonormal basis for $\mathcal{F}$. As in section \ref{Bergman Kernel}, $K_{\mathrm{BF}}(z,w)$ is given by 
\[
K_{\mathrm{BF}}(z,w)=\sum_{\alpha\in \N_0^n} \Psi_\alpha(z)\overline{\Psi_\alpha(w)}e^{-2\phi_0(w)},
\] 
for
\[
\int_{\C^n}\sum_{\alpha\in \N_0^n} \Psi_\alpha(z)\overline{\Psi_\alpha(w)}f(w)e^{-2\sum_{j=1}^n\lambda_{j,x}|w^j|^2}dm(w)=\sum_{\alpha\in \N_0^n}(f|\Psi_\alpha)_{\omega_0,\phi_0}\Psi_\alpha(z)=(\Pi_{\mathrm{BF}}f)(z).
\]
We then compute that
\begin{align*}
&K_{\mathrm{BF}}(z,w)=\sum_{\alpha\in \N_0^n} \frac{2^{|\alpha|+n}\lambda^{\alpha+1}}{\pi^n\alpha!}z^\alpha\overline{w}^\alpha e^{-2\sum_{j=1}^n\lambda_{j,x}|w^j|^2}\\
=&\sum_{d=1}^\infty \frac{2^{d+n}}{\pi^nd!} \sum_{|\alpha|=d}\frac{d!}{\alpha!} (\lambda z)^\alpha\overline{w}^\alpha e^{-2\sum_{j=1}^n\lambda_{j,x}|w^j|^2}
=\frac{2^n\lambda_{1,x}\cdots \lambda_{n,x}}{\pi^n}e^{2\sum_{j=1}^n\lambda_{j,x}(z^j\overline{w}^j-|w^j|^2)}.
\end{align*}
\end{proof}
Now, observe that $L^2(\C^n, e^{-2\sum_{j=1}^n\lambda_{j,x}|z^j|^2}dm)$ and $L^2(\C^n)$ is isometric via $u\mapsto ue^{\phi_0(z)}$. As in section \ref{Notation and Set-Up}, we know that
\[
\db (ue^{\sum_{j=1}^n\lambda_{j,x}|z^j|^2})=e^{\sum_{j=1}^n\lambda_{j,x}|z^j|^2}\db_{\phi_0}u,\quad \forall u\in C^\infty(\C^n).
\]
Therefore, $\mathcal{P}_{\omega_0,\phi_0}$ and $\Pi_{\mathrm{BF}}$ are related by $\mathcal{P}_{\omega_0,\phi_0}=e^{-\phi_0}\Pi_{\mathrm{BF}}e^{\phi_0}$, and their Schwartz kernel have the relation
\begin{equation}
\begin{aligned}
&P_{\omega_0,\phi_0}(z,w)=e^{-\phi_0(z)}K_{\mathrm{BF}}(z,w)e^{\phi_0(w)}\\
=&\frac{2^n\lambda_{1,x}\cdots \lambda_{n,x}}{\pi^n}e^{\sum_{j=1}^n\lambda_{j,x}\left(2z^j\overline{w}^j-|z^j|^2-|w^j|^2\right)}.
\end{aligned}
\end{equation}
\end{proof}
Our goal is to obtain asymptotic of $P_{\h{\omega},k\h{\phi}}$. Recall that in \eqref{extension of local metric data}, $\widehat{\phi}=\phi_0+\phi_1$ with $\phi_1\in C^\infty_c(\C^n)$. We consider $e^{-k(\h{\phi}-\phi_0)}u=e^{-k\phi_1}u$. Notice that $\phi_1\in C^\infty_c(\C^n)$ implies that
\[
\int_{\C^n} |u|^2e^{\pm2 k\phi_1}dm(z)<\infty,\quad \forall k\in \N,u\in L^2(\C^n).
\]
Hence, by similar argument as in Proposition \ref{Schwartz kernel for simple model case}, the map $u\mapsto ue^{-k\phi_1}$ defines an isometry on $L^2(\C^n)\to L^2(\C^n)$ which maps $\ker\db_{k\phi_0}$ bijectively onto $\ker\db_{k\h{\phi}}$ with inverse map $v\mapsto e^{k\phi_1}v$. On the other hand, we consider the change of base metric from $\omega_0$ to $\h{\omega}$. Observe that $\h{\omega}=\omega_0+\omega_1$ with $\omega_1$ supports in $B_{k^{\epsilon-\frac{1}{2}}}(0)$. This implies that the $L^2$-norm $\|\cdot\|_{\omega_0}$ and $\|\cdot\|_{\h{\omega}}$ are equivalent and thus $L^2(\C^n)=L^2(\C^n,\h{\omega})$. We may regard $\mathcal{P}_{\h{\omega},k\h{\phi}}:L^2(\C^n)\to \ker\db_{k\h{\phi}}$. 

We then define an intermediate operator $\widehat{\mathcal{P}}_{\omega_0,k\phi_0}:L^2(\C^n)\to L^2(\C^n)$ by
\begin{equation}
\widehat{\mathcal{P}}_{\omega_0,k\phi_0}=e^{-k\phi_1}\circ \mathcal{P}_{\omega_0,k\phi_0}\circ e^{k\phi_1}.
\end{equation}
By uniqueness of Schwartz kernel, we see that its Schwartz kernel $\widehat{P}_{\omega_0,k\phi_0}(z,w)$ is given by 
\begin{equation}\label{hat}
\widehat{P}_{\omega_0,k\phi_0}(z,w)=e^{-k\phi_1(z)}P_{\omega_0,k\phi_0}(z,w)e^{k\phi_1(w)}.
\end{equation}
Now, we observe that
\begin{lem}\label{lem:changing of metrics}
$\mathcal{P}_{\h{\omega},k\h{\phi}}=\widehat{\mathcal{P}}_{\omega_0,k\phi_0}\circ \mathcal{P}_{\h{\omega},k\h{\phi}}$ and $\widehat{\mathcal{P}}_{\omega_0,k\phi_0}=\mathcal{P}_{\h{\omega},k\h{\phi}}\circ\widehat{\mathcal{P}}_{\omega_0,k\phi_0}$.
\end{lem}
\begin{proof}
First, it is easy to see that the map $u\mapsto e^{-k\phi_1}u$ sends $\ker \db_{k\phi_0}$ to $\ker \db_{k\h{\phi}}$ and  $u\mapsto e^{k\phi_1}u$ sends $\ker \db_{k\h{\phi}}$ to $\ker \db_{k\phi_0}$. By \eqref{alternative Hodge decomposition} and above observation,
\begin{align*}
&\widehat{\mathcal{P}}_{\omega_0,k\phi_0}-\mathcal{P}_{\h{\omega},k\h{\phi}}\circ \widehat{\mathcal{P}}_{\omega_0,k\phi_0}=(I-\mathcal{P}_{\h{\omega},k\h{\phi}})\circ \widehat{\mathcal{P}}_{\omega_0,k\phi_0}\\
=&\db_{k\h{\phi}}^*(\triangle_{\h{\omega},k\h{\phi}}^{(1)})^{-1}\db_{k\h{\phi}}e^{-k\phi_1}\mathcal{P}_{\omega_0,k\phi_0}e^{k\phi_1}=0;\\
&\mathcal{P}_{\h{\omega},k\h{\phi}}-\widehat{\mathcal{P}}_{\omega_0,k\phi_0}\circ \mathcal{P}_{\h{\omega},k\h{\phi}}=e^{-k\phi_1}(I-\mathcal{P}_{\omega_0,k\phi_0})e^{k\phi_1}\circ \mathcal{P}_{\h{\omega},k\h{\phi}}\\
=&e^{-k\phi_1}\db_{k\phi_0}^{*,\omega_0}(\triangle_{\omega_0,k\phi_0}^{(1)})^{-1}\db_{k\phi_0}e^{k\phi_1}\mathcal{P}_{\h{\omega},k\h{\phi}}=0.
\end{align*}
\end{proof}
Moreover, let $\h{\mathcal{P}}^{*,\h{\omega}}_{\omega_0,k\phi_0}$ be the formal adjoint of $\widehat{\mathcal{P}}_{\omega_0,k\phi_0}$ with respect to the norm $(\cdot|\cdot)_{\h{\omega}}$. By direct computation, we know that its Schwartz kernel is given by
\begin{equation}\label{adjoint term}
\widehat{P}^{*,\h{\omega}}_{\omega_0,k\phi_0}(z,w)=\lambda^{-1}(z)e^{k\phi_1(z)}P_{k\phi_0}(z,w)e^{-k\phi_1(w)}\lambda(w),
\end{equation}
where $\lambda$ is the density of $\frac{\h{\omega}^n}{dm}$, i.e. $\h{\omega}_n=\lambda dm$. If we define $\mathcal{R}:=\widehat{\mathcal{P}}^{*,\h{\omega}}_{\omega_0,k\phi_0}-\widehat{\mathcal{P}}_{\omega_0,k\phi_0}$ to measure the extent which $\widehat{\mathcal{P}}_{\omega_0,k\phi_0}$ is not formally self-adjoint with respect to $(\cdot|\cdot)_{\h{\omega}}$, then its Schwartz kernel is given by
\begin{equation}\label{remainder term}
R(z,w)=P_{\omega_0,k\phi_0}(z,w)\left(\lambda^{-1}(z)\lambda(w)e^{k\phi_1(z)-k\phi_1(w)}-e^{k\phi_1(w)-k\phi_1(z)}\right).
\end{equation}
Now, if we take adjoint in the first formula in Lemma \ref{lem:changing of metrics}, we get
\[
\mathcal{P}_{\h{\omega},k\h{\phi}}=\mathcal{P}_{\h{\omega},k\h{\phi}}\circ \widehat{\mathcal{P}}^{*,\h{\omega}}_{\omega_0,k\phi_0}=\mathcal{P}_{\h{\omega},k\h{\phi}}\circ(\widehat{\mathcal{P}}_{\omega_0,k\phi_0}+\mathcal{R})=\widehat{\mathcal{P}}_{\omega_0,k\phi_0}+\mathcal{P}_{\h{\omega},k\h{\phi}}\circ \mathcal{R},
\]
where we use the second formula in Lemma \ref{lem:changing of metrics} in the last line. We then get
\begin{equation}\label{Neumann series-0}
\mathcal{P}_{\h{\omega},k\h{\phi}}(I-\mathcal{R})=\widehat{\mathcal{P}}_{\omega_0,k\phi_0}.
\end{equation}
 Now, for any $M\in \N$, if we multiply $(I+\mathcal{R}+\mathcal{R}^2+\cdots+\mathcal{R}^{M-1})$ from the right on the both sides of \eqref{Neumann series-0}, then we obtain
\begin{equation}\label{Neumann series}
\widehat{\mathcal{P}}_{\omega_0,k\phi_0}+\widehat{\mathcal{P}}_{\omega_0,k\phi_0}\circ\mathcal{R}+\cdots+ \widehat{\mathcal{P}}_{\omega_0,k\phi_0}\circ \mathcal{R}^{M-1}+\mathcal{P}_{\h{\omega},k\h{\phi}}\circ \mathcal{R}^M=\mathcal{P}_{\h{\omega},k\h{\phi}}.
\end{equation}
\eqref{Neumann series} is the key observation for establishing asymptotic expansion for $P_{\h{\omega},k\h{\phi}}$ near $(0,0)$. We will now employ the symbolic calculus developed in previous section to \eqref{Neumann series} to achieve this. First of all, from Proposition \ref{Schwartz kernel for simple model case} and
\[
2z^j\overline{w}^j-|z^j|^2-|w^j|^2=-|z^j-w^j|^2+2\ii \Imp z^j\overline{w}^j, \quad \forall  z,w\in \C^n,
\]
we know that
\[
P_{\omega_0,\phi_0}(z,w)=\frac{2^n\lambda_{1,x}\cdots \lambda_{n,x}}{\pi^n}e^{-\sum_{j=1}^n\lambda_{j,x}|z^j-w^j|^2+2\sqrt{-1}\Imp z^j\overline{w}^j}\in \h{S}(\C^n\times \C^n).
\]
Therefore, $P_{\omega_0,k\phi_0}(z,w)=k^nP_{\omega_0,\phi_0}(\sqrt{k}z,\sqrt{k}w)\in \h{S}^n_{\cl}(\C^n\times \C^n)$. Now, we show \begin{lem}\label{models in symbol space} For $\epsilon\in [0,1/6)$, we have
\[
\widehat{P}_{\omega_0,k\phi_0}(z,w)\in \h{S}^n_{\cl}(\C^n\times \C^n),\quad R(z,w)\in \h{S}^{n-\frac{1}{2}}_{\cl}(\C^n\times \C^n).
\]
\end{lem}
\begin{proof} By our choice of $\phi_1$ and $\omega_1$ as in \eqref{extension of local metric data}, we know that for $|z|,|w|>k^{\epsilon-1/2}$, $\widehat{P}_{\omega_0,k\phi_0}(z,w)=P_{\omega_0,k\phi_0}(z,w)\in \h{S}^n_{\cl}(\C^n\times \C^n)$. For $|z|,|w|<k^{\epsilon-1/2}$, since $\phi_1(z)=O(|z|^3)$, $|k\phi_1(z)|\leq Ck|z|^3 \leq Ck^{-1/2}|\sqrt{k}z|^3$. Since $|z|<k^{\epsilon-1/2}$, we see that $|k\phi_1(z)|\leq Ck^{3\epsilon-1/2}$. This shows that 
\begin{equation}\label{remainder term is lower degree}
\left|e^{k\phi_1(x)-k\phi_1(y)}-1\right|\leq Ce^{k^{3\epsilon-\frac{1}{2}}}\sup_{|x|,|y|<k^{-\frac{1}{2}+\epsilon}}(k|x|^3+k|y|^3)\leq Ck^{3\epsilon-1/2}e^{k^{3\epsilon-\frac{1}{2}}},
\end{equation}
where $x,y$ are the underlying real coordinates for $z$ and $w$. Hence, if $\epsilon\in [0,\frac{1}{6})$, then $3\epsilon-\frac{1}{2}<0$. Therefore, $\h{P}_{\omega_0,k\phi_0}\in \h{S}^n(\C^n\times \C^n)$. Furthermore,
\begin{equation}\label{Taylor expansion for kphi_1}
\begin{aligned}
&\left|e^{k\phi_1(x)-k\phi_1(y)}-\sum_{l=1}^N \frac{(k\phi_1(x)-k\phi_1(y))^l}{l!}\right|\\
\lesssim &\int_0^1 |k\phi_1(x)-k\phi_1(y)|^N e^{k\phi_1(tx)-k\phi_1(ty)}dt\\
\lesssim &k^{(N+1)(3\epsilon-1/2)}e^{k^{3\epsilon-1/2}}.
\end{aligned}
\end{equation}
This shows that $\h{P}_{\omega_0,k\phi_0}\in \h{S}^n_{\cl}(\C^n\times \C^n)$. Now, if we expand $\lambda(x)\lambda^{-1}(y)$ in Taylor expansion:
\begin{align*}
&\lambda^{-1}(x)\lambda(y)=1+\sum_{j=1}^{N-1}\sum_{|\alpha|+|\beta|=j} \frac{\p_x^\alpha \p_y^\beta (\lambda^{-1}(x)\lambda(y))(0,0)}{(\alpha+\beta)!}x^\alpha y^\beta\\
+&N\int_0^1 (1-t)^{N-1}\sum_{|\alpha|+|\beta|=N} \frac{\p_x^\alpha \p_y^\beta (\lambda^{-1}(x)\lambda(y))(tx,ty)x^\alpha y^\beta}{(\alpha+\beta)!}dt.
\end{align*}
This shows that
\begin{equation}\label{Taylor expansion for lambda}
\begin{aligned}
&\frac{\lambda(y)}{\lambda(x)}-1-\sum_{j=1}^{N-1}k^{-\frac{j}{2}}\sum_{|\alpha|+|\beta|=j} \frac{\p_x^\alpha \p_y^\beta (\lambda^{-1}(x)\lambda(y))(0,0)}{(\alpha+\beta)!}(\sqrt{k}x)^\alpha (\sqrt{k}y)^\beta\\
=&O(k^{-\frac{N}{2}}).
\end{aligned}
\end{equation}
Hence,  we get
\begin{align*}
&|R(z,w)|\\
\leq&|P_{\omega_0,k\phi_0}(z,w)|\left(\left|\lambda^{-1}(z)\lambda(w)-1\right|\left|e^{k\phi_1(z)-k\phi_1(w)}\right|+\left|e^{k\phi_1(z)-k\phi_1(w)}-e^{k\phi_1(w)-k\phi_1(z)}\right|\right)\\
\leq &|P_{\omega_0,k\phi_0}(z,w)|\left( \left|\lambda^{-1}(z)\lambda(w)-1\right|(1+Ck^{3\epsilon-1/2}e^{k^{3\epsilon-\frac{1}{2}}})+2Ck^{3\epsilon-1/2}e^{k^{3\epsilon-\frac{1}{2}}}\right).
\end{align*}
The derivative estimate of $R$ follows similarly as above. This shows that $R\in \h{S}^{n-\frac{1}{2}}(\C^n\times \C^n)$ if $\epsilon\in [0,\frac{1}{6})$. Moreover, \eqref{Taylor expansion for kphi_1} and \eqref{Taylor expansion for lambda} shows that $R\in \h{S}^{n-\frac{1}{2}}_{\cl}(\C^n\times \C^n)$.
\end{proof}
From Theorem \ref{thm:composition and adjoint of operators}, we know that for any $j\in \N$ and $R_j:=R^{\#j}$,
\[
\widehat{P}_{\omega_0,k\phi_0}\#R_j\in \h{S}^{n-j/2}_{\cl}(\C^n\times \C^n).
\]
Before proving our main result for this section, we need to first show that the remainder kernel $P_{\h{\omega},k\h{\phi}}\#R^{\#j}$ for $\mathcal{P}_{\h{\omega},k\h{\phi}}\circ \mathcal{R}^j$ in \eqref{Neumann series} is well-defined. 
\begin{lem}\label{lem:finiteness of remainder kernel} Let $R_j:=R^{\#j}$ be the Schwartz kernel of $\mathcal{R}^j$. Then $\mathcal{P}_{\h{\omega},k\h{\phi}}\circ \mathcal{R}^j$ is well-defined as a smoothing operator with smoothing kernel $P_{\h{\omega},k\h{\phi}}\#R_j$, for any $j\in \N$.
\end{lem}
\begin{proof}
For any $\alpha,\beta\in \N_0^{2n}$ any $x_0,y_0\in \C^n$, by Cauchy--Schwartz inequality,
\begin{align*}
&\left|\p_x^{\alpha}\p_y^\beta\int_{\C^n} P_{\h{\omega},k\h{\phi}}(x_0,u)R_j(u,y_0)dm(u)\right|
\leq  \int_{\C^n}|\p_x^\alpha P_{\h{\omega},k\h{\phi}}(x_0,u)||\p_y^\beta R_j(u,y_0)|dm(u)\\
\leq &\left(\int_{\C^n} |\p_x^\alpha P_{\h{\omega},k\h{\phi}}(x_0,u)|^2dm(u)\right)^{1/2}\left(\int_{\C^n} |\p_y^\beta R_j(u,y_0)|^2dm(u)\right)^{1/2}.
\end{align*}  
First of all, since $R_j\in \h{S}^{n-{j/2}}(\C^n\times \C^n)$, we know that $R_j(\cdot,y)\in \mathscr{S}(\R^{2n})$ for fixed $y_0$ and thus $(\p_y^\beta R_j)(\cdot,y_0)\in L^2(\C^n)$, for any $\beta\in \N_0^{2n}$. On the other hand, since $\|\cdot\|_{\omega_0}$ and $\|\cdot\|_{\h{\omega}}$ are equivalent, there exists a constant $C>0$ such that
\begin{align*}
&\int_{\C^n} |\p_x^\alpha P_{\h{\omega},k\h{\phi}}(x_0,u)|^2dm(u)\leq C\int_{\C^n} |\p_x^\alpha P_{\h{\omega},k\h{\phi}}(x_0,u)|^2\h{\omega}_n(u)\\
=&C\p_x^\alpha \p_y^\alpha P_{\h{\omega},k\h{\phi}}(x_0,x_0)<\infty.
\end{align*}
It is clear that $P_{\h{\omega},k\h{\phi}}\#R_j$ is the Schwartz kernel of $\mathcal{P}_{\h{\omega},k\h{\phi}}\circ \mathcal{R}^j$ and thus $\mathcal{P}_{\h{\omega},k\h{\phi}}\circ \mathcal{R}^j$ is a smoothing operator, for $P_{\h{\omega},k\h{\phi}}\#R_j$ is smooth.
\end{proof}
Hence, the kernel version of \eqref{Neumann series} is well-defined: 
\begin{equation}\label{Neumann series kernel}
\widehat{P}_{\omega_0,k\phi_0}+\h{P}_{\omega_0,k\phi_0}\# R+\cdots+ \widehat{P}_{\omega_0,k\phi_0}\# R_{M-1}+P_{\h{\omega},k\h{\phi}}\# R_M=P_{\h{\omega},k\h{\phi}}, \forall M\in \N.
\end{equation}
Also, we need the following simple observation. Let $\chi,\widetilde{\chi}\in C^\infty_c(\C^n)$ with
\[
\supp\chi\subset B_1(0),\quad \supp\widetilde{\chi}\subset B_2(0),\quad \widetilde{\chi}=1\text{ on }\supp\chi, \quad \chi=1\text{ on }B_{1/2}(0).
\]
and set $\chi_k(z):=\chi(8k^{1/2-\epsilon}z)$ and $\widetilde{\chi}_k(z):=\widetilde{\chi}(8k^{1/2-\epsilon}z)$.
\begin{lem}\label{lem:multiplication by k-cut-off} For any $a\in \h{S}^m(\C^n\times \C^n)$, $\widetilde{\chi}_k(x)a(x,y,k)\chi_k(y)\in \h{S}^m(\C^n\times \C^n)$.
\end{lem}
\begin{proof}
For any $\alpha,\beta\in \N_0^n$, we estimate 
\[
|\p_x^\alpha\p_y^\beta a(x,y,k)\widetilde{\chi}_k(x)\chi_k(y)\leq C_{\alpha,\beta} \sum_{\alpha'\leq \alpha,\beta'\leq \beta} |\p_x^{\alpha'}\widetilde{\chi}_k||\p_y^{\beta'}\chi_k||\p_x^{\alpha-\alpha'}\p_y^{\beta-\beta'}a|.
\]
Since $a\in \h{S}^m(\C^n\times \C^n)$, for each $\alpha-\alpha'$, $\beta-\beta'$, there exists $l(\alpha-\alpha',\beta-\beta')\in \N$ such that for any $N\in \N$, for any $x,y\in \C^n$, we have
\[
|\p_x^{\alpha-\alpha'}\p_y^{\beta-\beta'}a(x,y,k)|\leq C_{\alpha-\alpha',\beta-\beta',N}k^{m+\frac{|\alpha-\alpha'|+|\beta-\beta'|}{2}}\frac{(1+\sqrt{k}|x|+\sqrt{k}|y|)^{l(\alpha-\alpha',\beta-\beta')}}{(1+\sqrt{k}|x-y|)^N}.
\]
On the other hand, we have $|\p_x^{\alpha'}\widetilde{\chi}_k|\leq C_{\alpha'} k^{|\alpha'|(1/2-\epsilon)}$ and $ |\p_y^\beta\chi_k|\leq C_{\beta'}k^{|\beta'|(1/2-\epsilon)}$. Hence, we conclude that
\[
|\p_x^\alpha\p_y^\beta a(x,y,k)\widetilde{\chi}_k(x)\chi_k(y)|\leq C_{\alpha,\beta,N}k^{m+\frac{|\alpha|+|\beta|}{2}}\frac{(1+\sqrt{k}|x|+\sqrt{k}|y|)^{l(\alpha,\beta)}}{(1+\sqrt{k}|x-y|)^N}, \quad \forall x,y\in \C^n,
\]
where $l(\alpha,\beta):=\max\{l(\alpha-\alpha',\beta-\beta'):\alpha'\leq \alpha,\beta'\leq \beta\}$.
\end{proof}
Now, we are ready to establish the asymptotic expansion of $P_{\h{\omega},k\h{\phi}}$ near $(0,0)$.
\begin{thm}\label{Kernels in Symbol Space} For $\epsilon\in [0,1/6)$, we have
\[
\widetilde{\chi}_k(x)P_{\h{\omega},k\h{\phi}}(x,y)\chi_k(y)\in \h{S}^n_{\cl}(\C^n\times \C^n),
\]
where $\chi_k$, $\widetilde{\chi}_k$ as above.
\end{thm}
\begin{proof} We first show that $\widetilde{\chi}_kP_{\h{\omega},k\h{\phi}}\chi_k\in \h{S}^n(\C^n\times \C^n)$. For $z,w$ in a compact set $K$ of $0\in \C^n$, by standard scaling technique (cf. \cite[Theorem 4.3]{HsiaoMa}), one can prove that for any $\alpha\in \N_0^{2n}$, there exists a constant $C=C_K>0$ such that any $u\in \ker\db_{k\h{\phi}}$,
\[
|(\p_x^\alpha u)(z)|\leq C_{\alpha,K}k^{\frac{n+|\alpha|}{2}}\|u\|_{\h{\omega}},\quad \forall z\in K. 
\]
Let $\{\Psi_j\}_{j=1}^{d_k}$ be an orthonormal basis of $\ker \triangle_{\h{\omega},k\h{\phi}}$ with respect to $(\cdot|\cdot)_{\h{\omega}}$. Fix $\alpha\in \N_0^{2n}$ and $x_0\in K$, we may assume that $\sum_{j=1}^{d_k} |\p_x^\alpha \Psi_j(x_0)|^2\neq 0$. We then set
\[
u(z):=\frac{\sum_{j=1}^{d_k}\overline{(\p_x^\alpha \Psi_j)(x_0)}\Psi_j(z)}{\left(\sum_{j=1}^{d_k}|\p_x^\alpha \Psi_j(x_0)|^2\right)^{\frac{1}{2}}}.
\]
Since $P_{\h{\omega},k\h{\phi}}(z,w)=\sum_{j=1}^{d_k}\Psi_j(z)\overline{\Psi_j(w)}$ is smooth, the sum $\sum_{j=1}^{d_k} |\p_x^\alpha \Psi_j(x_0)|^2$ converges, and thus $u\in \ker\triangle_{\omega_0,k\h{\phi}}$ and $\|u\|^2_0=1$. By above argument, there exists a constant $C_{\alpha,K}$ so that 
\[
\left(\sum_{j=1}^{d_k} |\p_x^\alpha \Psi_j(x_0)|^2\right)^{\frac{1}{2}}=|(\p_x^\alpha u)(x_0)|\leq C_{\alpha,K}k^{\frac{n+|\alpha|}{2}}.
\]
Since $|\p_x^\alpha\p_y^\beta (P_{\h{\omega},k\h{\phi}})(x_0,x_0)|$ is dominated by $\left(\sum_j|\p_x^\alpha \Psi_j(x_0)|^2\right)^{1/2}\left(\sum_j |\p_y^\beta \Psi_j(x_0)|^2\right)^{1/2}$,
\begin{equation}\label{on-diagonal estimate}
|\p_x^\alpha \p_y^\beta(P_{\h{\omega},k\h{\phi}})(x_0,x_0)|\leq C_{\alpha,\beta,K}k^{n+\frac{|\alpha|+|\beta}{2}},
\end{equation}
and the same estimates holds for any $z\in K$ with the same constant $C_{\alpha,K}$. 

Now, for off-diagonal estimates, we notice that $|z|<\frac{1}{4}k^{-1/2+\epsilon}$, $|w|<\frac{1}{8}k^{-1/2+\epsilon}$. Therefore, $|z-w|<\frac{5}{8}k^{-1/2+\epsilon}$. For any $M\in \N$, we now multiply $\widetilde{\chi}_k(x)\chi_k(y)$ on the both sides of \eqref{Neumann series kernel}:
\begin{align*}
&\widetilde{\chi}_k(x)P_{\h{\omega},k\h{\phi}}(x,y)\chi_k(y)\\
=&\widetilde{\chi}_k(x)(P_{\h{\omega},k\h{\phi}}\#R_M)(x,y)\chi_k(y)+\sum_{j=0}^{M-1} \widetilde{\chi}_k(x)\h{P}_{\omega_0,k\phi_0}\#R_j(x,y)\chi_k(y).
\end{align*}
By Lemma \ref{lem:multiplication by k-cut-off}, we know that $\widetilde{\chi}_k\h{P}_{\omega_0,k\phi_0}\#R_j\chi_k\in \h{S}^{n-j/2}(\C^n\times \C^n)$. Hence, given any $N\in \N$, to estimate $(1+\sqrt{k}|z-w|)^N|\widetilde{\chi}_kP_{\h{\omega},k\h{\phi}}\chi_k|$, tt remains to estimate $(1+\sqrt{k}|z-w|)^N|\widetilde{\chi}_kP_{\h{\omega},k\h{\phi}}\#R_M\chi_k|$. We observe that \begin{align*}
&(1+\sqrt{k}|z-w|)^N\widetilde{\chi}_k(z)\left|\int_{\C^n} P_{\h{\omega},k\h{\phi}}(z,u)R_M(u,w)dm(u)\right|\chi_k(w)\\
\leq &\left(1+\frac{5}{8}k^{\epsilon}\right)^N C_Lk^{n-M/2} \int_{\C^n} |P_{\h{\omega},k\h{\phi}}(z,u)|\frac{(1+\sqrt{k}|u|+\sqrt{k}|w|)^l}{(1+\sqrt{k}|u-w|)^L}dm(u)\\
\leq &\left(1+\frac{5}{8}k^{\epsilon}\right)^N C_Lk^{n-M/2} |P_{\h{\omega},k\h{\phi}}(z,z)|^{1/2}\left(\int_{\C^n}\frac{(1+\sqrt{k}|u|+\sqrt{k}|w|)^{2l}}{(1+\sqrt{k}|u-w|)^{2L}}dm(u)\right)^{1/2}.
\end{align*}
By above on-diagonal estimates, we know that $|P_{\h{\omega},k\h{\phi}}(z,z)|^{1/2}\lesssim k^{n/2}$. Now, as in the proof of Theorem \ref{thm:composition and adjoint of operators}, we make the change of variable $t=\sqrt{k}u-\sqrt{k}w$:
\[
\int_{\C^n}\frac{(1+\sqrt{k}|u|+\sqrt{k}|w|)^{2l}}{(1+\sqrt{k}|u-w|)^{2L}}dm(u)\leq k^{-n}\int_{\C^n}\frac{(1+|t|+2\sqrt{k}|w|)^{2l}}{(1+|t|)^{2L}}dm(t).
\]
By choosing $L=l+2n$,  we get
\begin{align*}
&(1+\sqrt{k}|z-w|)^N\widetilde{\chi}_k(z)\left|(P_{\h{\omega},k\h{\phi}}\#R_M)(z,w)\right|\chi_k(w)\\
\leq& C_{n,l}(1+\frac{5}{8}k^{\epsilon})^Nk^{n-M/2}(1+\sqrt{k}|w|)^l.
\end{align*}
If we choose $M>2\epsilon N$, then $k^{-M/2}(1+\frac{5}{8}k^{\epsilon})^N\leq 2^N$. Similar estimate works for $(1+\sqrt{k}|z-w|)^N|\p_x^\alpha\p_y^\beta(\widetilde{\chi}_k(P_{\h{\omega},k\h{\phi}}\#R_M)\chi_k)|$ with the same $M$ but now $l$ may depends on $\alpha,\beta\in \N_0^{2n}$.  Hence, we conclude that $\widetilde{\chi}_kP_{\h{\omega},k\h{\phi}}\chi_k\in \h{S}^n(\C^n\times \C^n)$.

Finally, we show that $\widetilde{\chi}_kP_{\h{\omega},k\h{\phi}}\chi_k\in \h{S}^n_{\cl}(\C^n\times \C^n)$. To see this, for any $\alpha,\beta\in \N_0^{2n}$, any $M,N\in \N$, above argument shows that one can find $M'>M$ so that 
\begin{align*}
&(1+\sqrt{k}|z-w|)^N\left|\p_x^\alpha\p_y^\beta \widetilde{\chi}_k(x)(P_{\h{\omega},k\h{\phi}}\#R_{M'})(x,y)\chi_k(y)\right|\\
\leq &C_{\alpha,\beta,N}k^{n-\frac{M}{2}}(1+\sqrt{k}|z|+\sqrt{k}|w|)^{l(\alpha,\beta)}.
\end{align*}
 Therefore, by Lemma \ref{lem:multiplication by k-cut-off} and above estimate,
\begin{align*}
&\widetilde{\chi}_k(x)P_{\h{\omega},k\h{\phi}}(x,y)\chi_k(y)-\sum_{j=0}^M  \widetilde{\chi}_k(x)\h{P}_{\omega_0,k\phi_0}\#R_j(x,y)\chi_k(y)\\
=&\chi_k(x)(P_{\h{\omega},k\h{\phi}}\#R_{M'})(x,y)\chi_k(y)+\sum_{j=M+1}^{M'-1} \widetilde{\chi}_k(x)\h{P}_{\omega_0,k\phi_0}\#R_j(x,y)\chi_k(y)
\end{align*}
is in $\h{S}^{n-M/2}(\C^n\times \C^n)$. In view of Theorem \ref{thm:asymptotic sum for semi-classical symbol space}, we conclude that
\[
\widetilde{\chi}_k(x)P_{\h{\omega},k\h{\phi}}(x,y)\chi_k(y)\sim \sum_{j=0}^\infty  \widetilde{\chi}_k(x)\h{P}_{\omega_0,k\phi_0}\#R_j(x,y)\chi_k(y)
\]
and thus $\widetilde{\chi}_k(x)P_{\h{\omega},k\h{\phi}}(x,y)\chi_k(y)\in \h{S}^n_{\cl}(\C^n\times \C^n)$.
\end{proof}
\section{Localization of Global Bergman Kernel} \label{sec:Localization of Bergman Kernel via Spectral Gap}                                                                                                                                                                                                                
In this section, we complete the proof of Theorem  \ref{thm:Pointwise Asymptotic of Bergman Kernel Function} by localizing global Bergman kernel to the approximate Bergman kernel whose asymptotic expansion is already established in section \ref{Asymptotic Expansion of Approximate Kernel}. 

Let us recall the setting of Theorem \ref{thm:Pointwise Asymptotic of Bergman Kernel Function}. Let $(L,h)\to (X,\omega)$ be a a hermitian holomorphic line bundle over a complex manifold. We assume that $X(0)\neq \emptyset$. Given such open set $U$ of $x\in X(0)$, we choose a holomorphic coordinate $z=(z^1,\dots,z^n)$ and a holomorphic trivialization $s$ for $L$ on the open set $U\Subset X$ as in Lemma \ref{Canonical Coordinate}. Hence, we may identify the  data $(L,h)\to (X,\omega)$ near a point $x$ to an open neighborhood $U$ of $0\in\C^n$ with hermitian metric $\omega=\frac{\ii}{2}\sum_{i,j=1}^n H_{ij}(z)dz^i\wedge d\bar{z}^j$ and trivial line bundle $U\times \C$ with bundle metric given by $|1|^2_{\phi}=e^{-2\phi}$. In \eqref{extension of local metric data}, we consider a cut-off function $\theta_k\in C^\infty_c(U)$ with $\theta_k=1$ on $V_k:=\{|z|<\frac{1}{2}k^{\epsilon-\frac{1}{2}}\}\Subset U$ and extends $\omega$ and $\phi$ to $\h{\omega}$ and $\h{\phi}$ to whole $\C^n$ by
\[
\h{\omega}=\omega_0+\theta_k(\omega-\omega_0),\quad \h{\phi}=\phi_0+\theta_k(\phi-\phi_0).
\]
As in section \ref{sec:Approximate Bergman Kernel}, we have defined the deformed Cauchy--Riemann operator by
\[
\db_{k\h{\phi}}:=\db+k\epsilon(\db \h{\phi}),
\]
and the approximate Bergman kernel $\mathcal{P}_{\h{\omega},k\h{\phi}}:L^2(\C^n,\h{\omega})\to \ker\db_{k\h{\phi}}$. On the other hand, in section \ref{Notation and Set-Up}, we  have defined the \textbf{localized Bergman projection} $\Pi_{k,s}:L^2_{\comp}(U,\omega)\to \ker\db_{k,s}^{(0)}\subset L^2(U,\omega)$ by
\begin{equation}\label{localized Bergman projection}
\Pi_{k,s}u=e^{-k\phi}s^{-k}\Pi_k(e^{k\phi}us^k),\quad u\in L^2_{\comp}(U,\omega),
\end{equation}
where $\Pi_k:L^2(X)\to \ker \square_{\omega,k\phi}^{(0)}$ is the global Bergman kernel and $\db_{k,s}^{(0)}$ is the localized Cauchy--Riemann operator defined in \eqref{localized dbar:defin}. By definition of cut-off, we know that $\h{\omega}=\omega$ and $\h{\phi}=\phi$ on $V_k$. Thus, $\h{\omega}_n=\frac{\h{\omega}^n}{n!}=\frac{\omega^n}{n!}=\omega_n$. This implies that
\[
(u|v)_{\omega}=\int_{V_k} u\overline{v}\omega_n=\int_{V_k} u\overline{v}\h{\omega}_n=(u|v),\quad \forall u,v\in L^2(V_k,\omega)=L^2(V_k,\h{\omega}).
\]
and $\db_{k,s}u=\db u+k\epsilon(\db\phi)=\db u+k\epsilon(\db\h{\phi})$, for $u\in C^\infty(V_k)$.

Our goal is to establish the relation between $\Pi_{k,s}$ and $P_{\h{\omega},k\h{\phi}}$. To achieve this, we need to modify approximate Bergman kernel to a kernel defined on $U$. First, we consider a sequence of bump functions $\{\psi_i\}_{i=1}^\infty\subset C^\infty_c(U,[0,1])$ such that for any compact set $K\subset U$, $K\cap \supp \psi_i\neq \emptyset$, for only finitely many $i$, and $\sum_{i=1}^\infty \psi_i=1$ on $U$, and we define
\begin{equation}\label{properly support function}
\eta(z,w):=\sum_{\supp \psi_i\cap \supp \psi_j\neq \emptyset}\psi_i(z)\psi_j(w)
\end{equation}
\begin{lem}\label{lem:properly support function} $\eta\in C^\infty(U\times U)$ and $\eta\equiv1$ on a neighborhood $\Omega$ of the diagonal $\Delta_U\subset U\times U$. Furthermore, the projection $\supp\eta\to U$ on both $z$ and $w$ directions are proper maps.
\end{lem}
\begin{proof} Clearly, it suffices to prove that $\eta$ is smooth on a neighborhood of any $(x_0,y_0)\in U\times U$. For any $(x_0,y_0)\in U\times U$, any neighborhoods $W,W'\Subset U$ of $x_0$ and $y_0$, respectively. By construction of $\{\psi_i\}$, we know that there exist only finitely many $i,j\in \N$ such that $\supp \psi_i\cap W\neq\emptyset$ and $\supp \psi_j\cap W'\neq \emptyset$. Therefore, the sum in \eqref{properly support function} is a finite sum on $W\times W'$ and thus $\eta\in C^\infty(W\times W')$. For the second assertion, observe that 
\begin{align*}
&1-\eta(z,w)=\sum_{i=1}^\infty \psi_i(z)-\chi(z,w)=\sum_{i=1}^\infty \psi_i(z)\left(\sum_{j=1}^\infty \psi_j(w)\right)-\chi(z,w)\\
=&\sum_{i,j=1}^\infty \psi_i(z)\psi_j(w)-\sum_{\supp \psi_i\cap \supp \psi_j\neq \emptyset}\psi_i(z)\psi_j(w)\\
=&\sum_{\supp \psi_i\cap \supp \psi_j=\emptyset}\psi_i(z)\psi_j(w).
\end{align*}
If $1-\eta(z,w)\neq 0$, then $\psi_i(z)\neq 0$ and $\psi_j(w)\neq 0$, for some pair $(i,j)$ with $\supp \psi_i\cap \supp\psi_j=\emptyset$. We know that for such pair $(i,j)$, $z\in \supp \psi_i$ and $w\in \supp \psi_j$ and thus $(z,w)\in \supp\psi_i\times \supp \psi_j$. This shows that
\[
\supp(1-\eta)=\bigcup_{\supp\psi_i\cap \supp\psi_j=\emptyset}\supp\psi_i\times \supp\psi_j,
\] 
and thus the intersection of $\supp(1-\eta)$ with the diagonal of $U\times U$ is  empty. Furthermore, for each $z_0\in V$, if $z_0\in \supp \psi_{i_0}$, then there exists only finitely many $j_1,\dots,j_N$ such that $\supp \psi_{i_0}\cap \supp\psi_{j_l}\neq \emptyset$, for $l=1,\dots,N$. Thus, we can pick a neighborhood $W$ of $z_0$ such that $W\cap \bigcup_{l=1}^N\supp \psi_{j_l}=\emptyset$. This shows that $(W\times W)\cap \supp(1-\eta)=\emptyset$. As a result, $\eta\equiv 1$ on an open neighborhood $\Omega$ of the diagonal.

Finally, for each compact set $K\subset U$, the pre-image of it under the first projection is then given by $(K\times U)\cap \supp \eta$ Since there exists only finitely many index $i_1,\dots,i_N$  such that $\supp\psi_{i_l}\cap K\neq \emptyset$, for $l=1,\dots,N$, and for each $l=1,\dots,N$, there exists only finitely many $j_{l,m}$, for $m=1,\dots,M$, such that $\supp \psi_{i_l}\cap \supp\psi_{j_{l,m}}\neq \emptyset$. This shows that
\begin{align*}
\supp\eta\cap (K\times U)&=(K\times U)\cap \bigcup_{\supp\psi_i\cap \supp\psi_j\neq\emptyset }\supp\psi_i\times \supp\psi_j\\
&=\bigcup_{l=1}^N\supp\psi_{i_l}\times \bigcup_{m=1}^M\supp\psi_{j_{l,m}},
\end{align*}
and thus $\supp\eta\cap (K\times U)$ is compact. The proof for second projection is the same.
 \end{proof}
We define \textbf{localized approximate projection} $\h{\Pi}_k:L^2_{\comp}(U,\omega)\to L^2(U,\omega)$ by
\begin{equation}
(\h{\Pi}_ku)(z)=\int_U P_{\h{\omega},k\h{\phi}}(z,w)\eta(z,w)u(w)\omega_n(w)
\end{equation}
whose Schwartz kernel is given by $\h{K}_k(z,w)=P_{\h{\omega},k\h{\phi}}(z,w)\eta(z,w)$. By Lemma \ref{lem:properly support function}, we know that the projections $(z,w)\in \supp\eta \to w\in U$ is proper, and thus $\eta(z,w)u(w)$ also has compact support in $z$. This shows that $\h{\Pi}_k:L^2_{\comp}(U,\omega)\to L^2_{\comp}(U,\omega)$. On the other hand, by the properness of $(z,w)\in \supp \eta\to w\in U$, we know that for $u\in L^2(U,\omega)$, any $\tau\in C^\infty_c(U)$, $\tau(z)(\h{\Pi}_ku)(z)\in L^2(U,\omega)$. This implies that $\h{\Pi}_k:L^2(U,\omega)\to L^2_{\loc}(U,\omega)$.

On the other hand, we also define \textbf{localized approximate projection $\widetilde{\Pi}_k:L^2(U,\omega)\to L^2_{\comp}(U,\omega)$ concentrated near origin} by 
\begin{equation}
(\widetilde{\Pi}_ku)(z):=\widetilde{\chi}_k(z)(\h{\Pi}_k(\chi_ku))(z),
\end{equation}
where $\widetilde{\chi}_k(z)=\widetilde{\chi}(8k^{1/2-\epsilon}z)$, $\chi_k(z)=\chi(8k^{1/2-\epsilon}z)$ and
\[
\supp\chi_k\subset B_1(0),\quad \supp\widetilde{\chi}\subset B_2(0), \quad \widetilde{\chi}=1\text{ on }\supp\chi,\quad \chi=1\text{ on }B_{1/2}(0).
\]
By construction, we know that $\supp\widetilde{\chi}_k,\supp\chi_k\subset V_k$. We denote $\widetilde{K}_k(z,w)$ by the Schwartz kernel of $\widetilde{\Pi}_k$. Therefore, we have
\[
\widetilde{K}_k(z,w)=\widetilde{\chi}_k(z)\h{K}_k(z,w)\chi_k(w)=\widetilde{\chi}_k(z)P_{\h{\omega},k\h{\phi}}(z,w)\eta(z,w)\chi_k(w).
\]
We first prove the a crucial result which is important in our later arguments.
\begin{thm}\label{thm:off-diagonal k-negligible}
For $\epsilon\in (0,1/6)$, $(1-\widetilde{\chi}_k)\mathcal{P}_{\h{\omega},k\h{\phi}}\chi_k$ is a $k$-negligible operator in the sense of Definition \ref{defin:k-negligible}.
\end{thm}
In view of Definition \ref{defin:k-negligible}, it suffices to prove that for any $l,N\in \N$, any compact set $K\subset \C^n\times \C^n$, $\|(1-\widetilde{\chi}_k(x))P_{\h{\omega},k\h{\phi}}(x,y)\chi_k(y)\|_{C^l(K)}\leq C_{N,l}k^{-N}$, for some constant $C_{N,l}>0$ independent of $k$. We first prove a lemma.
\begin{lem} For any $m\in \R$, $a\in \h{S}^m(\C^n\times \C^n)$, 
\[
(1-\widetilde{\chi}_k(z))P_{\h{\omega},k\h{\phi}}(z,w)\chi_k(w)\equiv 0\mod O(k^{-\infty}).
\]
\end{lem}
\begin{proof}
Notice that $(1-\widetilde{\chi}_k)a\chi_k$ supports in $(\supp \chi_k)^c\times \supp \chi_k$. Since $\delta:=d(\supp(1-\widetilde{\chi}),\supp\chi)>0$, we only need to consider
\[
|w|<\frac{1}{8}k^{-1/2+\epsilon},\quad |z|>\frac{1}{8}k^{-1/2+\epsilon},\quad \text{ and }|z-w|\geq \frac{1}{8}k^{-\frac{1}{2}+\epsilon}\delta.
\]
Now, given any $\alpha,\beta\in \N_0^{2n}$, any $L\in \N$, since $a\in \h{S}^m(\C^n\times \C^n)$, we have the following esimtate
\[
\left|\p_x^\alpha\p_y^\beta (1-\widetilde{\chi}_k)a\chi_k\right|\leq C_{\alpha,\beta,L}k^{m+\frac{|\alpha|+|\beta|}{2}-\epsilon(|\alpha|+|\beta|)}\frac{(1+\sqrt{k}|z|+\sqrt{k}|w|)^{l(\alpha,\beta)}}{(1+\sqrt{k}|z-w|)^L}.
\]
Now, given any compact set $K\subset \C^n\times \C^n$ with $(\supp \chi_k)^c\times \supp \chi_k$, we may assume $ \frac{1}{8}k^{-1/2+\epsilon}<|z|\leq R$, for some $R:=R_K>0$ depending on $K$.  Hence, we get
\begin{align*}
&\left|\p_x^\alpha\p_y^\beta (1-\widetilde{\chi}_k)a\chi_k\right|\leq C_{\alpha,\beta,L}k^{m+\frac{|\alpha|+|\beta|}{2}}\frac{(1+|k|^\epsilon+\sqrt{k}R)^{l(\alpha,\beta)}}{(1+8k^\epsilon\delta)^L}\\
\leq &C_{\alpha,\beta,N}k^{m+\frac{|\alpha|+|\beta|+l(\alpha,\beta)}{2}}k^{-L\epsilon}(8\delta)^{-L}R^{l(\alpha,\beta)}\leq C_{\alpha,\beta,L,K}(8\delta)^{-L}k^{m+\frac{|\alpha|+|\beta|+l(\alpha,\beta)}{2}-\epsilon L}.
\end{align*}
For any $N\in \N$, we choose $L>\frac{m+\frac{|\alpha|+|\beta|+l(\alpha,\beta)}{2}+N}{2\epsilon}$. Therefore, we see that
\[
\left|\p_x^\alpha\p_y^\beta (1-\widetilde{\chi}_k)a\chi_k\right|\leq C_{\alpha,\beta,m,N,\epsilon,K}k^{-N}.
\]
\end{proof}

\begin{proof}[Proof of Theorem \ref{thm:off-diagonal k-negligible}]
Now, if we multiply \eqref{Neumann series kernel} by $(1-\widetilde{\chi}_k)(z)\chi_k(w)$, then by above Lemma,  
\[
(1-\widetilde{\chi}_k)(z)\h{P}_{\omega_0,k\phi_0}\#R_j\chi_k\equiv 0 \mod O(k^{-\infty}),
\] 
for any $j\in \N$. Similar to Lemma \ref{lem:finiteness of remainder kernel}, we estimate
\begin{align*}
&|(1-\widetilde{\chi}_k(z))\left|\int_{\C^n} P_{\h{\omega},k\h{\phi}}(z,u)R_M(u,w)dm(u)\right|\chi_k(w)\\
\leq & |1-\widetilde{\chi}_k(z)|P_{\h{\omega},k\h{\phi}}(z,z)|^{1/2}C_Lk^{n-M/2}\left(\int_{\C^n} \frac{(1+|t|+2\sqrt{k}|w|)^{2l}}{(1+|t|)^{2L}}\right)^{1/2}|\chi_k(w)|\\
\leq &C_{L,K}k^{n-\frac{M}{2}+\frac{l}{2}}\leq C_{n,l,K}k^{-N},
\end{align*}
where we choose $L=L_0(n,l)$ so that the integral converges and $M>2n+l+2N$, for any $N\in \N$. The derivatives estimates proceeds in similar fashion but $l$ may depends on the degree of differentiation. In conclusion, for any $N\in \N$, any $\alpha,\beta\in \N_0^{2n}$, there exists $M>2n+l(\alpha,\beta)+2N\in \N$ and $C_{n,\alpha,\beta,K}>0$ independent of $k$ such that
\[
\sup_K \left|\p_x^\alpha\p_y^\beta (1-\widetilde{\chi}_k)P_{\h{\omega},k\h{\phi}}\#R_M\chi_k\right|\leq C_{n,\alpha,\beta,K}k^{-N}.
\]
Hence, given any $\alpha,\beta\in \N_0^{2n}$, any $N\in \N$, we choose $M>2n+l(\alpha,\beta,2N)$ so that
\begin{align*}
&\sup_K \left|\p_x^\alpha\p_y^\beta (1-\widetilde{\chi}_k(z))P_{\h{\omega},k\h{\phi}}(z,w)\chi_k(w)\right|\\
\leq& \sum_{j=0}^{M-1}\sup_K \left|\p_x^\alpha\p_y^\beta (1-\widetilde{\chi}_k)(z)(\widehat{P}_{\omega_0,k\phi_0}\#R_j)(z,w)\chi_k(w)\right|\\
+&\sup_K \left|\p_x^\alpha\p_y^\beta (1-\widetilde{\chi}_k)(P_{\h{\omega},k\h{\phi}}\#R_M)(z,w)\chi_k(w)\right|\\
\leq& C_{\alpha,\beta,K,N}k^{-N}.
\end{align*}

\end{proof}
\begin{rmk}\label{rmk:off-diagonal k-negligible}
Notice that the condition on above Lemma can be relaxed. We actually proved that $\chi_kP_{\h{\omega},k\h{\phi}}\tau_k\equiv 0 \mod O(k^{-\infty})$, for $\chi_k\in C^\infty_c(\C^n,[0,1])$, $\tau_k\in C^\infty(\C^n,[0,1])$ with $d(\supp \chi_k,\supp \tau_k)>\frac{1}{8}k^{-1/2+\epsilon}\delta$, for some $\delta>0$ independent of $k$. Particularly, we can exchange the role of $1-\widetilde{\chi}_k$ and $\chi_k$. Also, by Theorem \ref{thm:alternative Hodge decomposition}, we see that
\begin{equation}\label{consequence of lemma}
(1-\widetilde{\chi}_k)\mathcal{P}_{\h{\omega},k\h{\phi}}\chi_k=-(1-\widetilde{\chi}_k)\db_{k\h{\phi}}^*(\triangle_{\h{\omega},k\h{\phi}}^{(1)})^{-1}\db_{k\h{\phi}}\chi_k\equiv 0\mod O(k^{-\infty}).
\end{equation}
\end{rmk}
\begin{thm}\label{approximate projection}
$\square^{(0)}_{k,s}\widetilde{\Pi}_k$ is $k$-negligible, i.e., $\square^{(0)}_{k,s}\widetilde{\Pi}_k\equiv 0\mod O(k^{-\infty})$ on $U$.
\end{thm}
\begin{proof} It suffices to prove $\db_{k,s}\widetilde{\Pi}_k\equiv 0\mod O(k^{-\infty})$ on $L^2(U,\omega)$. In view of Definition \ref{defin:k-negligible}, it suffices to prove that
\[
\db_{k,s}\left(\widetilde{\chi}_k(z)P_{\h{\omega},k\h{\phi}}(z,w)\eta(z,w)\chi_k(w)\right)\equiv 0\mod O(k^{-\infty}),\quad \text{ on } U.
\]
To see this, we write
\begin{align*}
&\widetilde{\chi}_k(z)P_{\h{\omega},k\h{\phi}}(z,w)\eta(z,w)\chi_k(w)\\
=&\widetilde{\chi}_k(z)P_{\h{\omega},k\h{\phi}}(z,w)\chi_k(w)+(1-\eta(z,w))\widetilde{\chi}_k(z)P_{\h{\omega},k\h{\phi}}(z,w)\chi_k(w)
\end{align*}
For the latter term, from the proof of Lemma \ref{lem:properly support function}, we see that
\begin{align*}
&(1-\eta(z,w))\widetilde{\chi}_k(z)P_{\h{\omega},k\h{\phi}}(z,w)\chi_k(w)\\
=&\sum_{\supp\psi_i\cap \supp\psi_j=\emptyset} \psi_i(z)\widetilde{\chi}_k(z)P_{\h{\omega},k\h{\phi}}(z,w)\chi_k(w)\psi_j(w).
\end{align*}
Now, notice that the proof of Theorem \ref{thm:off-diagonal k-negligible} works for $\psi_i(z)\widetilde{\chi}_k(z)$ and $\psi_j(w)\chi_k(w)$ (cf. Remark \ref{rmk:off-diagonal k-negligible}). Hence, we see that
\[
(1-\eta(z,w))\widetilde{\chi}_k(z)P_{\h{\omega},k\h{\phi}}(z,w)\chi_k(w)\equiv 0\mod O(k^{-\infty}).
\]
For the first term,
\begin{align*}
&\db_{k,s}(\widetilde{\chi}_k(z)P_{\h{\omega},k\h{\phi}}(z,w)\chi_k(w))=(\db+k\db\phi)(\widetilde{\chi}_k(z)P_{\h{\omega},k\h{\phi}}(z,w)\chi_k(w))\\
=&(\db \widetilde{\chi}_k)(z)P_{\h{\omega},k\h{\phi}}(z,w)\chi_k(w)+\widetilde{\chi}_k(z)\db_{k,s}(P_{\h{\omega},k\h{\phi}}(z,w))\chi_k(w).
\end{align*}
 Since $\chi_k(w)$ supports in $V_k$ and $\omega=\h{\omega}$, $\h{\phi}=\phi$ on $V_k$, we see that $\db_{k,s}(P_{\h{\omega},k\h{\phi}}(z,w))=0$. On the other hand, Theorem \ref{thm:off-diagonal k-negligible} (cf. again Remark \ref{rmk:off-diagonal k-negligible}) shows that
  \[
 (\db\widetilde{\chi}_k)(z)P_{\h{\omega},k\h{\phi}}(z,w)\chi_k(w)\equiv 0\mod O(k^{-\infty}).
 \]
\end{proof}

Since $\Pi_{k,s}:L^2_{\comp}(U,\omega)\to L^2(U,\omega)$ and $\widetilde{\Pi}_k:L^2(U,\omega)\to L^2_{\comp}(U,\omega)$, the composition $\Pi_k\circ \widetilde{\Pi}_{k,s}:L^2(U)\to L^2(U)$ makes sense. Recall that if the local spectral gap condition \eqref{spectral gap condition} holds on an open set $U\subset X(0)$. By local unitary identification in section \ref{Notation and Set-Up}, we have
\begin{equation}\label{spectral gap condition-local form}
\|(I-\Pi_{k,s})u\|_{\omega}\leq \frac{1}{Ck^d} \|\square^{(0)}_{k,s}u\|_\omega,\quad u\in C^\infty_c(U).
\end{equation}
Hence, we can prove
\begin{thm}\label{approximate projection-1} If the local spectral gap condition \eqref{spectral gap condition} holds on an open set $U\subset X(0)$, then the operator
$\Pi_{k,s}\widetilde{\Pi}_k-\widetilde{\Pi}_{k}$ is $k$-negligible on $U$.
\end{thm}
\begin{proof} For any $u\in C^\infty_c(U)$, we have the following estimate for $L^2$-norm.
\[
\|(\Pi_{k,s}\widetilde{\Pi}_{k}-\widetilde{\Pi}_{k})u\|_{\omega}=\|(\Pi_{k,s}-I)\widetilde{\Pi}_ku\|_\omega\leq C^{-1}k^{-d}\|\square_{k,s}^{(0)}\widetilde{\Pi}_ku\|_{\omega}.
\]
Using Theorem \ref{approximate projection}, $\square^{(0)}_{k,s}\widetilde{\Pi}_k\equiv 0\mod O(k^{-\infty})$. Thus, for any $N>0$, there exists a $k$-independent constant $C:=C_{M,N}>0$ so that $\|\square_{k,s}^{(0)}\widetilde{\Pi}_ku\|_{\omega}.
\leq Ck^{-N}\|u\|_{W^{-M},\omega}$, for any $M>0$. Also, notice that $
\square_{k,s}^{(0)}(\widetilde{\Pi}_{k}-\Pi_{k,s}\widetilde{\Pi}_{k})=\square_{k,s}^{(0)}\widetilde{\Pi}_k\equiv 0 \mod O(k^{-\infty})$.
By elliptic estimate, for any $u\in C^\infty_c(U)$, we know that there exists a constant $C>0$ independent of $k$ and $l>0$ so that for any $m\in \N$, any $N\in \N$,
\begin{align*}
&\|(\Pi_{k,s}\widetilde{\Pi}_{k}-\widetilde{\Pi}_{k})u\|_{W^{2m},\omega}\\
\leq &C k^{lm}(\|(\square_{k,s}^{(0)})^m(\Pi_{k,s}\widetilde{\Pi}_{k}-\widetilde{\Pi}_{k})u\|_{\omega}+\|(\Pi_{k,s}\widetilde{\Pi}_{k}-\widetilde{\Pi}_{k})u\|_{\omega})\\
\lesssim& k^{lm}(\|(\square_{k,s}^{(0)})^m\widetilde{\Pi}_ku\|_{\omega}+k^{-d}\|\square_{k,s}^{(0)}\widetilde{\Pi}_ku\|_{\omega})\lesssim k^{-N}\|u\|_{W^{-2m},\omega},\end{align*}
for coefficient of $\square_{k,s}^{(0)}$ has at most polynomial growth in $k$. By density argument, above estimate hold for any $u\in L^2(U,\omega)$ and thus
\[
\Pi_{k,s}\widetilde{\Pi}_{k}-\widetilde{\Pi}_{k}:W^{-2m}(U,\omega)\to W^{2m}(U,\omega)
\]
has operator norm $O(k^{-N})$, for any $m\in \N$ and $N\in \N$. We conclude that $\Pi_{k,s}\widetilde{\Pi}_k\equiv \widetilde{\Pi}_{k}\mod O(k^{-\infty})$.
\end{proof}
On the other hand, we prove
\begin{thm}\label{thm:almost analytic}
The operator $\chi_k\h{\Pi}_k\widetilde{\chi}_k\Pi_{k,s}-\chi_k\Pi_{k,s}$ is $k$-negligible on $U$.
\end{thm}
Before proving Theorem \ref{thm:almost analytic}, we need an asymptotic upper bound for $\Pi_{k,s}$.
\begin{lem}\label{upper bound for global kernel} For any $\alpha,\beta\in \N_0^{2n}$, 
there exists a constant $C_{\alpha,\beta,U}>0$ so that
\[
\left|\p_x^\alpha \p_y^\beta(K_{k,s})(x',y')\right|\leq C_{\alpha,\beta,U}k^{n+\frac{|\alpha|+|\beta}{2}},\quad \forall x',y'\in U.
\]
\end{lem}
\begin{proof}
This proceeds similar to the proof of \eqref{on-diagonal estimate} in the first part of the proof of Theorem \ref{Kernels in Symbol Space}. 
\end{proof}
\begin{proof}[Proof of Theorem \ref{thm:almost analytic}]
Since $\chi_k,\widetilde{\chi}_k$ supports in $V_k$ and $\h{\omega}=\omega$, $\h{\phi}=\phi$ on $V_k$, by Thoerem \ref{thm:alternative Hodge decomposition}, for any $u\in L^2_{\comp}(U,\omega)$, we can write
\begin{align*}
&(\chi_k\h{\Pi}_k\widetilde{\chi}_k\Pi_{k,s}u)(z)=\sum_{\supp \psi_i\cap \supp \psi_j\neq \emptyset}\chi_k(z)\psi_i(z)(\mathcal{P}_{\h{\omega},k\h{\phi}}(\psi_j\widetilde{\chi}_k\Pi_{k,s}u))(z)\\
=&-\sum_{\supp\psi_i\cap \supp\psi_j= \emptyset}\chi_k(z)\psi_i(z)(\mathcal{P}_{\h{\omega},k\h{\phi}}(\psi_j\widetilde{\chi}_k\Pi_{k,s}u))(z)+\chi_k(z)\mathcal{P}_{\h{\omega},k\h{\phi}}(\widetilde{\chi}_k\Pi_{k,s}u)(z)\\
=&-\sum_{\supp\psi_i\cap \supp\psi_j= \emptyset}\chi_k(z)\psi_i(z)(\mathcal{P}_{\h{\omega},k\h{\phi}}(\psi_j\widetilde{\chi}_k\Pi_{k,s}u))(z)\\
+&\chi_k(z)(I-\db_{k\h{\phi}}^*(\triangle_{\h{\omega},k\h{\phi}}^{(1)})^{-1}\db_{k\h{\phi}})(\widetilde{\chi}_k\Pi_{k,s}u)(z)\\
=&\chi_k(z)(\Pi_{k,s}u)(z)-\sum_{\supp\psi_i\cap \supp\psi_j=\emptyset}\chi_k(z)\psi_i(z)(\mathcal{P}_{\h{\omega},k\h{\phi}}(\psi_j\widetilde{\chi}_k\Pi_{k,s}u))(z)\\
-&\chi_k(z)\db_{k\h{\phi}}^*(\triangle_{\h{\omega},k\h{\phi}}^{(1)})^{-1}(\db\widetilde{\chi}_k)\Pi_{k,s}u)(z)-\chi_k(z)\db_{k\h{\phi}}^*(\triangle_{\h{\omega},k\h{\phi}}^{(1)})^{-1}\widetilde{\chi}_k(\db_{k\h{\phi}}\Pi_{k,s}u)(z).
\end{align*} 
Applying Theorem \ref{thm:off-diagonal k-negligible} and Remark \ref{rmk:off-diagonal k-negligible} to the last two terms, we have
\begin{align*}
\sum_{\supp\psi_i\cap \supp\psi_j=\emptyset}\psi_i\mathcal{P}_{\h{\omega},k\h{\phi}}\psi_j&\equiv 0 \mod O(k^{-\infty}),\\
\chi_k\db_{k\h{\phi}}^*(\triangle_{\h{\omega},k\h{\phi}}^{(1)})^{-1}(\db\widetilde{\chi}_k)&\equiv 0\mod O(k^{-\infty}).
\end{align*}
Also, $\h{\phi}=\phi$ on $V_k$, $\chi_k(z)\db_{k\h{\phi}}^*(\triangle_{\h{\omega},k\h{\phi}}^{(1)})^{-1}\widetilde{\chi}_k(\db_{k\phi}\Pi_{k,s}u)(z)=0$.

Combining with Lemma \ref{upper bound for global kernel} which shows that derivatives of $K_{k,s}$ is at most polynomials in $k$, we see that $\chi_k(z)\Pi_{k,s}\equiv \chi_k\widehat{\Pi}_k\widetilde{\chi}_k\Pi_{k,s}\mod O(k^{-\infty})$.
\end{proof}
\begin{proof}[Proof of Theorem \ref{thm:Pointwise Asymptotic of Bergman Kernel Function}] If we take adjoint in Theorem \ref{thm:almost analytic}, we get
\[
\Pi_{k,s}\chi_k\equiv \Pi_{k,s}\widetilde{\chi_k}\h{\Pi}_k^{*,\omega}\chi_k\mod O(k^{-\infty}),
\]
where $\h{\Pi}_k^{*,\omega}$ means the adjoint with respect to $\omega$. As in the proof of Theorem \ref{approximate projection-1} , since $\chi_k,\widetilde{\chi}_k$ supports in $V_k$ and $\h{\omega}=\omega$, $\h{\phi}=\phi$ on $V_k$, we see that 
\[
\widetilde{\chi_k}\h{\Pi}_k^{*,\omega}\chi_k=\widetilde{\chi}_k\h{\Pi}_k^{*,\h{\omega}}\chi_k=\widetilde{\chi_k}\h{\Pi}_k\chi_k.
\] 
The last identity follows from $P_{\h{\omega},k\h{\phi}}$ is self-adjoint (with respect to $\h{\omega}$), $\eta(z,w)=\overline{\eta(w,z)}=\eta(w,z)$, and  $\h{K}_k(z,w)=P_{\h{\omega},k\h{\phi}}(z,w)\eta(z,w)$. By Theorem \ref{approximate projection-1} and the assumption that $U$ satisfies local spectral gap, we conclude that
\[
\Pi_{k,s}\chi_k\equiv \Pi_{k,s}\widetilde{\Pi}_k\equiv \widetilde{\Pi}_k \mod O(k^{-\infty}).
\]
In terms of kernels, this shows that
\[
K_{k,s}(z,w)\chi_k(w)\equiv \widetilde{\chi}_k(z)P_{\h{\omega},k\h{\phi}}(z,w)\eta(z,w)\chi_k(w)\mod O(k^{-\infty}).
\]
Since $\eta(z,w)$ is $k$-independent, by multiplying $\eta$ to \eqref{Neumann series kernel}, we see that
\[
\widetilde{\chi}_k(z)P_{\h{\omega},k\h{\phi}}(z,w)\eta(z,w)\chi_k(w)\in \h{S}^n_{\cl}(\C^n\times \C^n) 
\]
and $\widetilde{\chi}_k(z)P_{\h{\omega},k\h{\phi}}(z,w)\eta(z,w)\chi_k(w)\sim \sum_{j=0}^\infty \widetilde{\chi}_k(z)(\h{P}_{\omega_0,k\phi_0}\#R^{\#j})(z,w)\eta(z,w)\chi_k(w)$. 

Finally, by multiplying $\rho$ on above, we have
\[
\rho(z)\widetilde{\chi}_k(z)P_{\h{\omega},k\h{\phi}}(z,w)\eta(z,w)\chi_k(w)\equiv \rho(z)K_{k,s}(z,w)\chi_k(w) \mod O(k^{-\infty}).
\]
We see that $ \rho(z)K_{k,s}(z,w)\chi_k(w)\in \h{S}^n_{\cl}(\C^n\times \C^n)$ as it is supported in $U\times V_k$.
\end{proof}
We also deduce Theorem \ref{the first coefficient}.
\begin{proof}[Proof of Theorem \ref{the first coefficient}]
From the proof of Theorem \ref{Kernels in Symbol Space}, we know that
\[
\widetilde{\chi}_k(z)P_{\h{\omega},k\h{\phi}}(z,w)\chi_k(w)\sim \sum_{j=0}^\infty \widetilde{\chi}_k(z)(\h{P}_{\omega_0,k\phi_0}\#R^{\#j})(z,w)\chi_k(w)
\] 
and by Theorem \ref{thm:off-diagonal k-negligible}, $\widetilde{\chi}_k(z)P_{\h{\omega},k\h{\phi}}(z,w)\eta(z,w)\chi_k(w)\equiv P_{\h{\omega},k\h{\phi}}(z,w)\chi_k(w)$. Hence, 
\[
\rho(z)K_{k,s}(z,w)\chi_k(w)\sim \rho(z) \widetilde{\chi}_k(z)P_{\h{\omega},k\h{\phi}}(z,w)\chi_k(w).
\]
The first coefficient in the asymptotic sum is given by
\[
\widehat{P}_{\omega_0,k\phi_0}(z,w)=\frac{2^nk^n\lambda_{1,x}\cdots \lambda_{n,x}}{\pi^n}e^{k\sum_{j=1}^n\lambda_{j,x}(2z^j\overline{w}^j-|z^j|^2-|w^j|^2)-k(\phi_1(z)-\phi_1(w))}.
\]
By \eqref{remainder term is lower degree}, we know that $1-e^{-k(\phi_1(z)-\phi_1(w))}$ is of lower degree in $k$, we get
\[
a_0(z,w)=\frac{2^n\lambda_{1,x}\cdots \lambda_{n,x}}{\pi^n}e^{\sum_{j=1}^n\lambda_{j,x}(2z^j\overline{w}^j-|z^j|^2-|w^j|^2)}.
\]
\end{proof}
\backmatter
\bibliographystyle{alpha}
\bibliography{Bergman}

\begin{thebibliography}{HKSX16}

\bibitem[BBS08]{BBS}
Robert Berman, Bo~Berndtsson, and Johannes Sj\"{o}strand.
\newblock A direct approach to {B}ergman kernel asymptotics for positive line
  bundles.
\newblock {\em Ark. Mat.}, 46(2):197--217, 2008.

\bibitem[BdMS76]{BMS}
L.~Boutet~de Monvel and J.~Sj\"{o}strand.
\newblock Sur la singularit\'{e} des noyaux de {B}ergman et de {S}zego.
\newblock In {\em Journ\'{e}es: \'{E}quations aux {D}\'{e}riv\'{e}es
  {P}artielles de {R}ennes (1975)}, Ast\'{e}risque, No. 34--35, pages 123--164.
  Soci\'{e}t\'{e} Math\'{e}matique de France, Paris, 1976.

\bibitem[Ber22]{Bergman}
Stefan Bergmann.
\newblock \"{U}ber die {E}ntwicklung der harmonischen {F}unktionen der {E}bene
  und des {R}aumes nach {O}rthogonalfunktionen.
\newblock {\em Math. Ann.}, 86(3-4):238--271, 1922.

\bibitem[BMS94]{Schlichenmaier1}
Martin Bordemann, Eckhard Meinrenken, and Martin Schlichenmaier.
\newblock Toeplitz quantization of {K}\"{a}hler manifolds and
  {${\mathrm{gl}}(N)$}, {$N\to\infty$} limits.
\newblock {\em Comm. Math. Phys.}, 165(2):281--296, 1994.

\bibitem[Bou90]{Bouche3}
Thierry Bouche.
\newblock Convergence de la m\'{e}trique de {F}ubini-{S}tudy d'un fibr\'{e}
  lin\'{e}aire positif.
\newblock {\em Ann. Inst. Fourier (Grenoble)}, 40(1):117--130, 1990.

\bibitem[Bou95]{Bouche2}
Thierry Bouche.
\newblock Two vanishing theorems for holomorphic vector bundles of mixed sign.
\newblock {\em Math. Z.}, 218(4):519--526, 1995.

\bibitem[BS07]{BS}
Robert Berman and Johannes Sj\"{o}strand.
\newblock Asymptotics for {B}ergman-{H}odge kernels for high powers of complex
  line bundles.
\newblock {\em Ann. Fac. Sci. Toulouse Math. (6)}, 16(4):719--771, 2007.

\bibitem[Cat99]{Catlin}
David Catlin.
\newblock The {B}ergman kernel and a theorem of {T}ian.
\newblock In {\em Analysis and geometry in several complex variables ({K}atata,
  1997)}, Trends Math., pages 1--23. Birkh\"{a}user Boston, Boston, MA, 1999.

\bibitem[CDS15a]{CDS1}
Xiuxiong Chen, Simon Donaldson, and Song Sun.
\newblock K\"{a}hler-{E}instein metrics on {F}ano manifolds. {I}:
  {A}pproximation of metrics with cone singularities.
\newblock {\em J. Amer. Math. Soc.}, 28(1):183--197, 2015.

\bibitem[CDS15b]{CDS3}
Xiuxiong Chen, Simon Donaldson, and Song Sun.
\newblock K\"{a}hler-{E}instein metrics on {F}ano manifolds. {III}: {L}imits as
  cone angle approaches {$2\pi$} and completion of the main proof.
\newblock {\em J. Amer. Math. Soc.}, 28(1):235--278, 2015.

\bibitem[CS01]{Chen-Shaw}
So-Chin Chen and Mei-Chi Shaw.
\newblock {\em Partial differential equations in several complex variables},
  volume~19 of {\em AMS/IP Studies in Advanced Mathematics}.
\newblock American Mathematical Society, Providence, RI; International Press,
  Boston, MA, 2001.

\bibitem[Dav95]{Davies}
E.~B. Davies.
\newblock {\em Spectral theory and differential operators}, volume~42 of {\em
  Cambridge Studies in Advanced Mathematics}.
\newblock Cambridge University Press, Cambridge, 1995.

\bibitem[Dem82]{Demailly3}
Jean-Pierre Demailly.
\newblock Estimations {$L^{2}$} pour l'op\'{e}rateur {$\bar \partial $} d'un
  fibr\'{e} vectoriel holomorphe semi-positif au-dessus d'une vari\'{e}t\'{e}
  k\"{a}hl\'{e}rienne compl\`ete.
\newblock {\em Ann. Sci. \'{E}cole Norm. Sup. (4)}, 15(3):457--511, 1982.

\bibitem[Dem12]{Demailly2}
Jean-Pierre Demailly.
\newblock Complex analytic and algebraic geometry.
\newblock https://www-fourier.ujf-grenoble.fr/~demailly/manuscripts/agbook.pdf,
  2012.

\bibitem[DK10]{Klevtsov}
Michael~R. Douglas and Semyon Klevtsov.
\newblock Bergman kernel from path integral.
\newblock {\em Comm. Math. Phys.}, 293(1):205--230, 2010.

\bibitem[DLM06]{DLM}
Xianzhe Dai, Kefeng Liu, and Xiaonan Ma.
\newblock On the asymptotic expansion of {B}ergman kernel.
\newblock {\em J. Differential Geom.}, 72(1):1--41, 2006.

\bibitem[Don01]{Donaldson2001}
S.~K. Donaldson.
\newblock Scalar curvature and projective embeddings. {I}.
\newblock {\em J. Differential Geom.}, 59(3):479--522, 2001.

\bibitem[Don03]{Donnelly}
Harold Donnelly.
\newblock Spectral theory for tensor products of {H}ermitian holomorphic line
  bundles.
\newblock {\em Math. Z.}, 245(1):31--35, 2003.

\bibitem[DS14]{DS}
Simon Donaldson and Song Sun.
\newblock Gromov-{H}ausdorff limits of {K}\"{a}hler manifolds and algebraic
  geometry.
\newblock {\em Acta Math.}, 213(1):63--106, 2014.

\bibitem[FK72]{FollandKohn}
G.~B. Folland and J.~J. Kohn.
\newblock {\em The {N}eumann problem for the {C}auchy-{R}iemann complex}.
\newblock Annals of Mathematics Studies, No. 75. Princeton University Press,
  Princeton, N.J.; University of Tokyo Press, Tokyo, 1972.

\bibitem[Gaf55]{Gaffney}
Matthew~P. Gaffney.
\newblock Hilbert space methods in the theory of harmonic integrals.
\newblock {\em Trans. Amer. Math. Soc.}, 78:426--444, 1955.

\bibitem[GH78]{GH}
Phillip Griffiths and Joseph Harris.
\newblock {\em Principles of algebraic geometry}.
\newblock Pure and Applied Mathematics. Wiley-Interscience [John Wiley \&
  Sons], New York, 1978.

\bibitem[H\"65]{HormanderL2}
Lars H\"{o}rmander.
\newblock {$L^{2}$} estimates and existence theorems for the {$\bar \partial $}
  operator.
\newblock {\em Acta Math.}, 113:89--152, 1965.

\bibitem[H\"66]{HormanderSCV}
Lars H\"{o}rmander.
\newblock {\em An introduction to complex analysis in several variables}.
\newblock D. Van Nostrand Co., Inc., Princeton, N.J.-Toronto, Ont.-London,
  1966.

\bibitem[H\"03]{HormanderI}
Lars H\"{o}rmander.
\newblock {\em The analysis of linear partial differential operators. {I}}.
\newblock Classics in Mathematics. Springer-Verlag, Berlin, 2003.
\newblock Distribution theory and Fourier analysis, Reprint of the second
  (1990) edition [Springer, Berlin; MR1065993 (91m:35001a)].

\bibitem[HKSX16]{Seto2}
Hamid Hezari, Casey Kelleher, Shoo Seto, and Hang Xu.
\newblock Asymptotic expansion of the {B}ergman kernel via perturbation of the
  {B}argmann-{F}ock model.
\newblock {\em J. Geom. Anal.}, 26(4):2602--2638, 2016.

\bibitem[HM14]{HsiaoMa}
Chin-Yu Hsiao and George Marinescu.
\newblock Asymptotics of spectral function of lower energy forms and {B}ergman
  kernel of semi-positive and big line bundles.
\newblock {\em Comm. Anal. Geom.}, 22(1):1--108, 2014.

\bibitem[HS20]{HsiaoNikhill}
Chin-Yu Hsiao and Nikhil Savale.
\newblock Bergman-szeg\"{o} kernel asymptotics in weakly pseudoconvex finite
  type cases.
\newblock arXiv:2009.07159, Sep 2020.

\bibitem[Hsi15]{Hsiao}
Chin-Yu Hsiao.
\newblock Bergman kernel asymptotics and a pure analytic proof of the {K}odaira
  embedding theorem.
\newblock In {\em Complex analysis and geometry}, volume 144 of {\em Springer
  Proc. Math. Stat.}, pages 161--173. Springer, Tokyo, 2015.

\bibitem[Mar05]{Marinescu1}
George Marinescu.
\newblock The laplace operator on high tensor powers of line bundles.
\newblock Preprint, http://www. mi. uni-koeln. de/gmarines/PREPRINTS/habil.
  pdf, 2005.

\bibitem[MM06]{MM2}
Xiaonan Ma and George Marinescu.
\newblock The first coefficients of the asymptotic expansion of the {B}ergman
  kernel of the {$\mathrm{Spin}^c$} {D}irac operator.
\newblock {\em Internat. J. Math.}, 17(6):737--759, 2006.

\bibitem[MM07]{Mama}
Xiaonan Ma and George Marinescu.
\newblock {\em Holomorphic {M}orse inequalities and {B}ergman kernels}, volume
  254 of {\em Progress in Mathematics}.
\newblock Birkh\"{a}user Verlag, Basel, 2007.

\bibitem[MM08]{MM}
Xiaonan Ma and George Marinescu.
\newblock Generalized {B}ergman kernels on symplectic manifolds.
\newblock {\em Adv. Math.}, 217(4):1756--1815, 2008.

\bibitem[MM12]{MamaBT}
Xiaonan Ma and George Marinescu.
\newblock Berezin-{T}oeplitz quantization on {K}\"{a}hler manifolds.
\newblock {\em J. Reine Angew. Math.}, 662:1--56, 2012.

\bibitem[RSN20]{Siorstrand2}
Oph\'{e}lie Rouby, Johannes Sj\"{o}strand, and San V\~{u} Ng\d{o}c.
\newblock Analytic {B}ergman operators in the semiclassical limit.
\newblock {\em Duke Math. J.}, 169(16):3033--3097, 2020.

\bibitem[Rua98]{Ruan}
Wei-Dong Ruan.
\newblock Canonical coordinates and {B}ergman metrics.
\newblock {\em Comm. Anal. Geom.}, 6(3):589--631, 1998.

\bibitem[Sch10]{Schlichenmaier2}
Martin Schlichenmaier.
\newblock Berezin-{T}oeplitz quantization for compact {K}\"{a}hler manifolds.
  {A} review of results.
\newblock {\em Adv. Math. Phys.}, pages Art. ID 927280, 38, 2010.

\bibitem[Set15]{Seto}
Shoo Seto.
\newblock {\em On the Asymptotic Expansion of the Bergman Kernel}.
\newblock PhD thesis, UC Irvine, 2015.

\bibitem[Siu84]{Siu}
Yum~Tong Siu.
\newblock A vanishing theorem for semipositive line bundles over
  non-{K}\"{a}hler manifolds.
\newblock {\em J. Differential Geom.}, 19(2):431--452, 1984.

\bibitem[Tia90]{Tian}
Gang Tian.
\newblock On a set of polarized {K}\"{a}hler metrics on algebraic manifolds.
\newblock {\em J. Differential Geom.}, 32(1):99--130, 1990.

\bibitem[Yos95]{Yosida}
Kosaku Yosida.
\newblock {\em Functional analysis}.
\newblock Classics in Mathematics. Springer-Verlag, Berlin, 1995.
\newblock Reprint of the sixth (1980) edition.

\bibitem[Zel98]{Zelditch}
Steve Zelditch.
\newblock Szego kernels and a theorem of {T}ian.
\newblock {\em Internat. Math. Res. Notices}, (6):317--331, 1998.

\end{thebibliography}
\end{document}